%  The latest version
%  July 31 version

%  quantum R-matrices

%  This is an amstex file.
\mag 1200
\input amstex
\input amsppt.sty
\NoRunningHeads
\overfullrule 0pt

\hsize 6.25truein
\vsize 9.03truein

\let\al\alpha
\let\bt\beta
\let\gm\gamma 
\let\dl\delta 
\let\epe\epsilon \let\eps\varepsilon \let\epsilon\eps

\let\la\lambda 
\let\om\omega 
 \let\phi\varphi

\def\O{\Cal O}
\def\o{\otimes}

\def\d{\partial}
\def\<{\langle}
\def\>{\rangle}

\def\wo{\widetilde\o}

\def\C{\Bbb C}

\def\Z{\Bbb Z}

\def\F{\Cal F}

\def\g{\frak g}

  \def\Fwh{\,\wh{\!\F}}

\def\A{\roman A}

\def\E{\roman E}

\let\Fun\Fwh

% my notations

\def\A{\Cal A} 

\def\Fun/{ {\text{Fun}}} 
\def\End/{ {\text{End}}} 
\def\Hom/{ {\text{Hom}}} 
\def\Rm/{\^{$R$-}matrix} 

\def\bo{\bar\o}
\def\Vh{{\V_\h}}
\def\B{\Cal B}

\def\M/{\Cal M } 
\def\A/{$A_{\tau,\eta}(sl_2)$} 
\def\Am/{$A_{\tau,\eta,\mu}(sl_2)$} 

\def\E/{$E_{\tau,\eta}(sl_2)$} 

\def\Ex/{ \widetilde {  E   }}
\def\EE/{ $\widetilde {  E   }_{\tau,\eta}(sl_2)$}

\def\eqg/{elliptic quantum group}
\def\Em/{$E_{\tau,\eta,\mu}(sl_2)$} 
\def\E/{$E_{\tau,\eta}(sl_2)$}

\def\Ce/{\C/{1\over \eta}\Z }
\def\V{\Cal V }

\def \g{\frak g}

\def \h{\frak h}

\def\End {{\text{End}}}
\def\Id  {{\text{Id}}}

\def\coth  {{\text{cotanh}}}

\def\sin {{\text{sin}}}

\topmatter

\title
Solutions of the quantum dynamical 
Yang-Baxter equation and dynamical quantum groups\endtitle

\author
Pavel Etingof and Alexander Varchenko
\endauthor

\thanks
The authors were supported in part by an 
NSF postdoctoral fellowship and
NSF grant
DMS-9501290
\endthanks

\date
August, 1997
\enddate

\abstract
{The quantum dynamical Yang-Baxter (QDYB) equation is a useful generalization 
of the quantum Yang-Baxter (QYB) equation. This generalization was
 introduced by 
Gervais, Neveu, and Felder. Unlike the QYB equation, the QDYB equation is 
not an algebraic but a difference
equation, with respect to a matrix function rather than a matrix. 
The QDYB  equation and its quasiclassical analogue (the classical dynamical 
Yang-Baxter equation) arise in several areas of mathematics and 
mathematical physics (conformal field theory, integrable systems, 
representation theory). The most interesting solution 
of the QDYB equation is the elliptic solution, discovered by Felder.  

In this paper, we prove the first classification results 
for solutions of the QDYB equation. These results are parallel
to the classification of solutions of the classical dynamical Yang-Baxter 
equation, obtained in our previous paper. All solutions we found 
can be obtained from Felder's elliptic solution by a limiting process 
and gauge transformations. 

Fifteen years ago the quantum Yang-Baxter equation gave rise to the theory 
of quantum groups. Namely, it turned out that the language of quantum groups 
(Hopf algebras) is the adequate algebraic language to talk about solutions 
of the quantum Yang-Baxter equation. 

In this paper we propose a similar language, originating from 
Felder's ideas, which we found to be adequate  
for the dynamical Yang-Baxter equation. This is the language of dynamical 
quantum groups (or $\h$-Hopf algebroids), which is the quantum
counterpart of the language of dynamical Poisson groupoids,
introduced in our previous paper.  
}
\endabstract

 \leftheadtext{ P.Etingof and A.Varchenko }
  \rightheadtext{ Solutions of the quantum dynamical 
Yang-Baxter equation and dynamical quantum groups}
\endtopmatter

\document

\head Introduction\endhead

This paper is 
devoted to the quantum dynamical Yang-Baxter equation, its solutions, 
and the related algebraic structures (quantum groupoids,
Hopf algebroids); abusing language, we will call these structures by the
collective name {\bf ``dynamical quantum groups''}. 

Let $\h$ be a finite dimensional commutative Lie algebra over $\C$, 
$V$ a semisimple finite dimensional $\h$-module, and $\gamma$ a complex 
number.  
The quantum dynamical Yang-Baxter (QDYB) equation is the equation
$$
\align
R^{12}(\la - \gm h^{(3)})\,
R^{13}(\la )\, &
R^{23}(\la - \gm h^{(1)})\,
\tag 1
\\
=
R^{23}(\la )\,&
R^{13}(\la - \gm h^{(2)})\,
R^{12}(\la)\,
\endalign
$$
with respect to a meromorphic function $R:\h^*\to \End (V\o V)$, 
where by definition $R^{12}(\la-\gamma h^{(3)})(v_1\o v_2\o v_3):=
(R^{12}(\la-\gamma \mu)(v_1\o v_2))\o v_3$ if $v_3$ has weight $\mu$, and  
$R^{13}(\la-\gamma h^{(2)}), R^{23}(\la-\gamma h^{(1)})$ are defined 
analogously. 

It is also useful to consider the quantum dynamical Yang-Baxter equation 
with spectral parameter, with respect to a meromorphic function
$R:\C\times \h^*\to \End (V\o V)$. By definition, 
the QDYB equation with spectral parameter is just equation (1), with $R^{ij}(*)$ 
replaced by $R^{ij}(z_i-z_j,*)$, where $z_1,z_2,z_3\in \C$. 

Solutions of the QDYB equation which are invariant under $\h$ are called 
quantum dynamical R-matrices.
 
A brief history of the QDYB equation is as follows. 
The QDYB equation 
was proposed by Felder 
\cite{F2} as a quantization of the classical dynamical Yang-Baxter 
equation \cite{F1}, but it also appeared earlier in 
physical literature \cite{GN}. Examples of dynamical R-matrices 
appeared in \cite{Fad1,AF}). 
As Felder showed \cite{F2}, the QDYB equation 
 is equivalent to the star-trangle relation 
in statistical mechanics. The most interesting known solution of the 
QDYB equation with spectral parameter 
is the elliptic solution given in \cite{F1,F2}. As was shown 
in \cite{TV}, this solution 
arises when one studies monodromies of the quantum KZ equation 
introduced in \cite{FR}, see also \cite{FTV1-2}. 
The algebraic structure corresponding to this 
solution was described in \cite{F1,F2, FV1-3} 
and called ``the elliptic quantum group''. 
Although the elliptic quantum group is not a Hopf algebra, it is 
very similar to a Hopf algebra in many respects. For example, its category 
of representations, with a suitable definition of the tensor product,
 is a tensor category, which was studied in \cite{FV1,FV2}. 

This paper has two goals.  

1. To classify quantum dynamical R-matrices in the case 
when $\h\subset \End (V)$ is the algebra of all diagonal operators in some 
basis. 

2. To describe the axiomatics of the algebraic structure corresponding 
to a quantum dynamical R-matrix. 

The first goal is partially attained in Chapters 1 and 2. 

In Chapter 1, we study dynamical R-matrices without spectral parameter. We
define the notion of a dynamical R-matrix of Hecke type which is a dynamical 
R-matrix satisfying a generalized unitarity condition. Then we define 
gauge transformations, which map the set of such 
dynamical R-matrices to itself. After this, we classify 
dynamical R-matrices of Hecke type, with $\h$ as above. 
The answer turns out to be completely parallel to the 
classical case (\cite{EV}, Chapter 3). In particular, any 
classical dynamical r-matrix from \cite{EV} without spectral parameter
(for the Lie algebra $\frak{gl}_N$) can be quantized.

In Chapter 2, we study dynamical R-matrices with spectral parameter, 
satisfying the unitarity condition. As in Chapter 1, we define 
gauge transformations, which map the set of such 
dynamical R-matrices to itself. After this, we list 
all known examples, and give a partial classification result 
(for R-matrices given by a power series in $\gamma$, 
which are quantizations of elliptic r-matrices from \cite{EV}, Chapter 4).
As before, the results are parallel to the 
classical case. In particular, any 
classical dynamical r-matrix from \cite{EV} with spectral parameter
(for the Lie algebra $\frak{gl}_N$) can be quantized.

{\bf Remark.} 
We were not able to obtain a nice classification result for 
dynamical R-matrices with spectral parameter and numerical 
$\gamma$, since we do not understand what is the correct analogue of 
the residue condition in \cite{EV}. However, we expect that such a result 
can be obtained along the same lines as in Chapter 4 of \cite{EV}, 
and Chapter 
1 of this paper. 

The second goal is attained in Chapters 3-6. 

In Chapter 3, we explain the connection between dynamical R-matrices 
and monoidal categories. We introduce the tensor category of $\h$-vector
spaces, and show that a tensor functor from a braided monoidal category 
to the category of $\h$-vector spaces gives a dynamical R-matrix, 
in the same way as a tensor functor from a braided monoidal category 
to the category of vector spaces gives a usual R-matrix. 
We also attach to every dynamical R-matrix 
a tensor category of its representations, following the ideas of 
\cite{F1,F2,FV1,FV2}.   
This category is nontrivial (for example, it contains 
the basic representation), has natural notions of the left and right
dual objects, and is equipped with a canonical tensor functor to $\h$-vector 
spaces.  

In Chapter 4 we introduce the notions of an $\h$-algebra, 
$\h$-bialgebroid, and $\h$-Hopf algebroid, which are generalizations 
of the notions of an algebra, bialgebra, and Hopf algebra.
We define the notion of a dynamical representation of an $\h$-algebra, and 
show that the category 
of dynamical representations $\text{Rep}(A)$ of an $\h$-bialgebroid $A$
is a tensor category
with a natural tensor functor to $\h$-vector spaces. 
If $A$ is an $\h$-Hopf algebroid, this category in addition 
has natural notions of the left and right 
dual representation. 
 
Using a generalization of the Faddeev-Reshetikhin-Sklyanin-Takhtajan
formalism \cite{FRT},\cite{FT} 
which assigns a Hopf algebra to any R-matrix, 
we assign an $\h$-bialgebroid $A_R$
to any dynamical R-matrix  $R$.
If $R$ has an additional rigidity property, 
then $A_R$ is an $\h$-Hopf algebroid.     
We call the bialgebroid 
$A_R$ the 
dynamical quantum group associated to $R$. 
We show that the category of representations of $R$ is equivalent 
to the category $\text{Rep}(A_R)$ as a tensor category with duality
and with a functor to $\h$-vector spaces. 

In Chapter 5, we define quantum counterparts of the quasiclassical
objects defined in \cite{EV} (in the setting of perturbation theory).
 More specifically, we  
define the notions of a biequivariant algebra
(biequivariant quantum space), 
a biequivariant Hopf algebroid (biequvariant quantum groupoid), 
a dynamical Hopf
algebroid (dynamical quantum groupoid), which are the quantum analogues 
of the notions of a biequivariant Poisson algebra (biequivariant 
Poisson manifold),
a biequivariant Poisson-Hopf algebroid (biequivariant 
Poisson groupoid), a dynamical 
Poisson-Hopf algebroid (dynamical Poisson groupoid), introduced in \cite{EV}. 
We introduce the notion of quantization for biequivariant and 
dynamical objects, and conjecture that any 
dynamical Poisson groupoid can be quantized. 

This material is a generalization of the material of Chapter 4, because,
as we explain in Section 5.5, 
the notion of an $\h$-algebra ($\h$-bialgebroid, $\h$-Hopf algebroid) is 
essentially a special 
case of the notion of a biequivariant algebra 
(bialgebroid, Hopf algebroid). 

{\bf Remark.} The general notion of a Hopf algebroid was introduced by J.H.Lu 
\cite{Lu}. It is easy to check that bieqivariant and dynamical 
Hopf algebroids as defined in Chapter 5 of our paper are Hopf algebroids 
in the sense of Lu. However, the notion considered in \cite{Lu} is more
general than the one considered in the paper.  

In Chapter 6, we study $\h$-bialgebroids associated 
to dynamical R-matrices of strong Hecke type. 
Using the semisimplicity of the Hecke algebra for a generic value
of the parameter, we prove a Poincare-Birkhoff-Witt theorem 
for such bialgebroids. This result explains the meaning of the Hecke type 
condition, which was artificially introduced in Chapter 1.
Using the same method, we show that the $\h$-Hopf algebroid
associated to a dynamical R-matrix of Hecke type 
of the form $R=1-\gm r+..$ is a 
flat deformation (quantization) of the Poisson-Hopf algebroid 
corresponding to $r$. 
 
In the next papers, we plan to
 develop the theory of dynamical quantum groups.
We plan to  describe the infinite-dimensional dynamical quantum groups 
associated to dynamical R-matrices with spectral parameter, and dynamical 
quantum groups (both finite and infinite dimensional) associated 
to Lie groups other than $GL_N$.  
We plan to  develop the representation theory of dynamical 
quantum groups, and 
explain its connection
with exchange 
(Zamolodchikov) algebras, Kazhdan-Lusztig functors, 
KZ and quantum KZ equations. 

\head 1. Classification of Quantum Dynamical R-matrices
without spectral parameter
\endhead

\subhead 1.1. Quantum dynamical R-matrix
\endsubhead

Let $\h$ be an abelian finite dimensional Lie algebra.
A finite dimensional diagonalizable $\h$-module is a complex
finite dimensional vector space $V$ with a weight decomposition 
$V=\oplus_{\mu\in\h^*}V[\mu]$, such that $\h$ acts on $V[\mu]$
by $xv =\mu(x)v$,  where $x \in \h$, $v \in V[\mu]$.

Let $V_i, \, i=1, 2, 3,$ be finite dimensional
diagonalizable  $\h$ modules,
$$
R_{V_iV_j}\,:\, \h^* \,\to\, \End (V_i \otimes V_j), \qquad 1\leq i < j \leq 3,
$$
meromorphic functions, $\gm $ a nonzero complex number. 
The  equation in $\End (V_1 \otimes V_2 \otimes V_3)$,
$$
\align
R_{V_1V_2}^{12}(\la - \gm h^{(3)})\,
R_{V_1V_3}^{13}(\la )\,&
R_{V_2V_3}^{23}(\la - \gm h^{(1)})\,
\tag 1.1.1
\\
=\,R_{V_2V_3}^{23}(\la )\,&
R_{V_1V_3}^{13}(\la - \gm h^{(2)})\,
R_{V_1V_2}^{12}(\la)\,
\endalign
$$
is called  {\it the quantum dynamical Yang-Baxter
equation with step} $ \gm$ (QDYB equation).

Here we use the following notation.
 If $X\in\End(V_i)$, then
we denote by $X^{(i)}\in\End(V_1\otimes\dots\otimes V_n)$
the operator $\cdots\otimes\Id\otimes X\otimes\Id\otimes\cdots$, acting
non-trivially on the $i$-th factor of a tensor product of vector spaces,
and
if $X=\sum X_k\otimes Y_k\in\End(V_i\otimes V_j)$, then we set
$X^{ij}=\sum X_k^{(i)}Y_k^{(j)}$.   
The shift of $\la$ by $\gm h^{(i)}$ is defined in the standard way.
For instance,
$R_{V_1V_2}^{12}(\la - \gm h^{(3)})$ acts on 
a tensor
$v_1\otimes v_2\otimes v_3$ as 
$R_{V_1V_2}^{12}(\la - \gm \mu_3) \otimes \Id$
if $v_3$ has weight $\mu_3$.

A function
$R_{V_iV_j} \, : \h^* \to \End (V_i \otimes V_j)$ is called  {\it a
function of zero
weight} if
$$
[R_{V_iV_j}(\la), h \otimes 1 + 1 \otimes h ] \, = \, 0
\tag 1.1.2
$$
for all $h \in \h$, $\la \in \h^*$.
A solution $\{ R_{V_iV_j} \}_{1\leq i< j\leq 3}$  of the QDYB equation is called
a solution of zero weight if
each of the functions is of zero weight.

If all the spaces $V_i$ are equal to a space $V$,  then consider the QDYB equation on
one function $R \, : \h^* \to \End (V \otimes V)$,
$$
\align
R^{12}(\la - \gm h^{(3)})\,
R^{13}(\la )\, &
R^{23}(\la - \gm h^{(1)})\,
\tag 1.1.3
\\
=
R^{23}(\la )\,&
R^{13}(\la - \gm h^{(2)})\,
R^{12}(\la)\,.
\endalign
$$
 An invertible function $R$ of zero weight
satisfying
the QDYB equation (1.1.3) is called  {\it a quantum dynamical R-matrix}.

\subhead 1.2. Quantization and quasiclassical limit
\endsubhead

Let $x_1,...,x_N$ be a basis in $\h$. The basis defines a linear system of coordinates
on $\h^*$. For any $\la \in \h^*$, set $\la_i = x_i (\la)$, $i=1,..., N$.

Let $R_\gm: \h^* \to \End (V\otimes V)$ 
be a smooth family of solutions to the QDYB equation with step $\gm$
such that
$$
R_\gm(\la)\,=\, 1\,-\,\gm\,r(\la)\,+\, O(\gm^2).
\tag 1.2.1
$$
Then the function $r : \h^* \to \End (V\otimes V)$
satisfies the classical dynamical Yang-Baxter equation (CDYB),
$$
\align
 \, \sum _{i=1}^N \, x_i^{(1)} {\partial r^{23} 
\over \partial x_i}&\, + \, \sum_{i=1}^N \, x_i^{(2)} {\partial r^{31} 
\over \partial x_i}\, + \, \sum _{i=1}^N \, x_i^{(3)} {\partial r^{12 }
\over \partial x_i}\, +
\tag 1.2.2
\\
& \, [r^{12},r^{13}]\,+ \, [r^{12},r^{23}]\,
+\,[r^{13},r^{23}] \,= \,0\, .
\endalign
$$
A function $r$ of zero weight satisfying the CDYB 
equation is called  {\it a classical dynamical
r-matrix}.
The function $r$ in (1.2.1) is called  {\it the quasiclassical limit of} $R$,
and the function $R$ is called  {\it a quantization of} $r$.

Let $U\subset \h^*$ be an open set, and let
$R:U\to \End (V\o V)$ be a
zero weight meromorphic function on $U$. We will say that 
$R$ is a quantum dynamical R-matrix on $U$ if the QDYB equation is satisfied 
for $R$ whenever it makes sense. 

{\bf Remark.} If $U$ is a bounded set, this notion is only interesting 
for small $\gamma$, so that the QDYB equation 
makes sense on a nonempty open set $U'\subset U$.

A classical dynamical r-matrix $r(\la)$ on $U$ is called {\it quantizable}
if there exists a power series in $\gm$,
$$
R_\gm(\la)\,=\, 1\,-\,\gm\,r(\la)\,+\, \sum_{n=2}^\infty \gm^n r_n(\la),
\tag 1.2.3
$$
convergent for small 
$|\gamma|$ for any fixed $\la\in U$ and
such that $R_\gm(\la)$ is a quantum dynamical R-matrix on $U$ with step 
$\gamma$.

\subhead 1.3. Quantum dynamical R-matrices of  Hecke type
\endsubhead

Let $\h$ be an abelian Lie algebra of dimension $N$. Let $V$ be a diagonalizable
$\h$-module of the same dimension $N$ such that
its weights  $\omega_1, ... , \omega_N$
form a basis in $\h^*$. Let $x_1,...,x_N$ be the dual basis of $\h$.
Let $v_1, ..., v_N$ be an eigenbasis for $\h$ in $V$ such that $x_i v_j = \dl_{ij} v_j$.
Then the $\h$-module $V\otimes V$ has the weight decomposition,
$$
V\otimes V \,=\, \oplus_{a=1}^N V_{aa} \, \oplus
\oplus_{a < b} V_{ab}\,,
\tag 1.3.1
$$
where $V_{aa}\, =\, \C\,v_a\otimes v_a$ and $V_{ab}\, =\, \C\,v_a\otimes v_b \oplus
\C\, v_b \otimes v_a $ .

Introduce a basis  $E_{ij}$ in $ \End (V)$ by $E_{ij}v_k = \dl_{jk} v_i$.

A quantum dynamical R-matrix  $R : \h^* \to \End (V\otimes V)$
for these $\h$ and $V$ will be called an R-matrix of $gl_N$ type.

The zero weight condition implies that the R-matrix preserves the weight decomposition (1.3.1)
and has the form
$$
R(\la)\,=\, \sum_{a,b = 1}^N \,\al_{ab}(\la)\,E_{aa}\otimes E_{bb}\,+
\,\sum_{a\neq b } \,\bt_{ab}(\la)\,E_{ba}\otimes E_{ab}\,
\tag 1.3.2
$$
where $\al_{ab}, \bt_{ab} \,:\, \h^* \to \C$ 
are suitable meromorphic functions.

Let $P 
%\,=\, \sum_{a, b=1 }^N \,E_{ab}\otimes E_{ba}\, 
\in \End (V\otimes V)$
be the permutation of factors. Set $R^\vee = P R$.

Let $ p, q $ be nonzero complex numbers, $p \neq -q$. 
A function  $R : \h^* \to \End (V\otimes V)$ will be called  {\it a function of
Hecke type with parameters $p$, $q$} if
\item{1.3.3.}{  The function preserves the weight decomposition (1.3.1).
}
\item{1.3.4.}{ For any $a = 1, ..., N$ and
$\la \in \h^*$, we have  $R^\vee (\la) v_a\otimes v_a \,= \,p\, v_a\otimes v_a$.}
\item{1.3.5.}{ For any $a \neq b$ and
$\la \in \h^*$, the operator $R^\vee (\la)$ restricted to  the two dimensional space
$V_{ab}$ has eigenvalues $p$ and $ - q$.}

A function  $R : \h^* \to \End (V\otimes V)$ will be called  {\it a function of 
weak Hecke type with parameters $p$, $q$} if
it preserves the weight decomposition (1.3.1) and for any $\la \in \h^*$
satisfies the equation
$$
(R^\vee(\la) - p)\, (R^\vee(\la) + q)\,=\, 0\,.
\tag 1.3.6
$$

A relation between Hecke types is given by the following simple observation.
Let $R_t : \h^* \to \End (V\otimes V), t \in [0,1]$, be a continuous
 family of meromorphic 
functions, which is analytic when $t \in (0,1)$. Assume that
 for any $t$ the function $R_t$ is
of weak Hecke type  and   $R_{t=0} 
= \Id$. Then $R_t$ is of Hecke type
for any $t$. In fact, the matrix
$R^{\vee}_{t=0} = P$ satisfies (1.3.4-5) and hence
$R^{\vee}_{t}$ satisfies (1.3.4-5) for any $t$.

In the following sections we
 classify quantum dynamical R-matrices of $gl_N$ Hecke type.

\subhead 1.4. Gauge transformations and multiplicative closed 2-forms
\endsubhead 

In this subsection we introduce gauge transformations of
 quantum dynamical R-matrices of  Hecke type.
We shall use the notion of a multiplicative  form.

{\it A multiplicative  $k$-form} on a vector space
with a linear coordinate system $\la_1, ..., \la_N$ is a collection, 
$$
\phi \,=\,\{ \phi_{a_1, ... , a_k}(\la_1, ... ,
 \la_N)\}\,,
$$
of meromorphic
functions , where $a_1,..., a_k$ run through
all $k$ element subsets of $\{ 1, ..., N\}$, such that for any
subset $a_1,..., a_k$  and any $i,\, 1\leq i < k$, we have
$$
\phi_{a_1, ..., a_{i+1}, a_i, ...  , a_k}(\la_1, ..., \la_N)
\, \phi_{a_1, ... , a_k}(\la_1, ..., \la_N) \, = \, 1 \,.
$$
Let $\Omega^k$ be the set of all multiplicative $k$-forms.

If $\phi$
and $ \psi$
are  multiplicative $k$-forms, then
$\{ \phi_{a_1, ... , a_k}(\la_1, ..., \la_N) \, \cdot \,
\psi_{a_1, ... , a_k}(\la_1, ..., \la_N)\}$
and 
$\{ \phi_{a_1, ... , a_k}(\la_1, ..., \la_N)\,/\,
 \psi_{a_1, ... , a_k}(\la_1, ..., \la_N)\}$
are multiplicative $k$-forms. This gives an abelian group structure on $\Omega^k$.
The zero element in $\Omega^k$ is the form 
$\{ \phi_{a_1, ... , a_k}(\la_1, ..., \la_N) \equiv 1\}$.

Fix a nonzero complex number $\gm$. For any $a = 1, ...., N$, introduce an operator
$\dl_a$ 
on the space of meromorphic functions $f(\la_1, ..., \la_N)$ by
$$
\dl_a \,:\, f(\la_1, ..., \la_N) \, \mapsto \,
f(\la_1, ..., \la_N)\,/\, f(\la_1, ..., \la_a - \gm, ..., \la_N)\,
$$
and an operator $d_\gm : \Omega^{k} \to \Omega^{k+1},\,
\phi \mapsto d_\gm\phi,$ by
$$
(d_\gm \phi)_{a_1, ... , a_{k+1}}(\la_1, ..., \la_N)\,=\,
\prod_{i=1}^{k+1} \,(\dl_{a_i}\phi_{a_1, ... , a_{i-1},a_{i+1}, ...,a_{k+1}}
(\la_1, ..., \la_N))^{(-1)^{i+1}}.
$$
We have $d_\gm^2\,=\, 0$.
A form $\phi$ will be called {\it $\gm$-closed} if $d_\gm\phi = 0$.

Let $\phi(\gm) =\{ \phi_{a_1, ... , a_k}(\la_1, ... ,
 \la_N,\gm)\}$ be a smooth family of multiplicative $k$-forms
such that for all $a_1, ... , a_k$,
$$
\phi_{a_1, ... , a_k}(\la, \gm)\,=\, 1\,-\, \gm\,
C_{a_1, ... , a_k}(\la) + O(\gm^2)\,
$$
for suitable functions $C_{a_1, ... , a_k}(\la)$.
Then the functions
$\{ C_{a_1, ... , a_k}(\la) \}$ are skew-symmetric with respect 
to permutation of the indices, so it is natural to consider 
a differential form
$C =\sum_{a_1<...<a_k}\,C_{a_1, ... , a_k}(\la)\,dx_{a_1}
\wedge...\wedge dx_{a_k}$.  
The differential form $C$ is called the 
{\it quasiclassical limit} of the multiplicative 
form $\phi(\gm)$ and the multiplicative form
$\phi(\gm)$ is called a {\it quantization} of the differential form $C$.
It is easy to see that if $\phi(\gm)$ is $\gm$-closed, then $C$ is closed. 

Let $U\subset \C^N$ be an open set, and let
$\phi$ be a multiplicative 
meromorphic $k$-form on $U$. We will say that $\phi$ is 
$\gamma$-closed if the equation $d_\gamma\phi=0$ is satisfied 
whenever it makes sense. 

A closed differential form $\{C_{a_1,...,a_k}(\la)\}$  
is called {\it quantizable}
if there exists a power series in $\gm$, 
$$
\phi_{a_1,...,a_k}(\la, \gm)
\,=\, 1\,-\,\gm\,C_{a_1,...,a_k}(\la)\,+\,\sum_{n=2}^\infty \gm^n 
C_{n;\,a_1,...,a_k}(\la),
$$
convergent for small 
$|\gamma|$ for a fixed $\la\in U$ and  
such that $\{\phi_{a_1,...,a_k}(\la, \gm)\}$ is a $\gm$-closed 
multiplicative $k$-form.

\proclaim{Lemma 1.1}

Every closed holomorphic differential $k$-form $C$
defined on an open polydisc is quantizable to 
a holomorphic multiplicative closed $k$-form $\phi(\gamma)$.

\endproclaim

\demo{Proof} Since $U$ is a polydisc, we can find a holomorphic 
$(k-1)$-form $E$ on $U$ such that $dE=C$. Define a multiplicative 
$(k-1)$-form $\theta$ on $U$ by $\theta_{a_1...a_{k-1}}=
e^{- E_{a_1...a_{k-1}}}$. Set $\phi(\gm)=d_\gamma\theta$. 
Since $d_\gamma^2=0$, the form $\phi(\gm)$ is a desired 
multiplicative closed $k$-form. 
$\square$\enddemo

{\bf Remark.} The Taylor expansion of $\phi(\gamma)$ 
in powers of $\gm$ is well defined in $U$, but for each particular
(even very small) nonzero $\gamma$, the form $\phi(\gamma)$ is defined in  a
 smaller open subset $U'(\gamma)\subset U$ which tends to $U$ as $\gamma\to 0$.

Now we introduce gauge transformations of quantum dynamical
R-matrices, \newline
$R \, : \, \h^*  \to \End (V \otimes V)$, of form (1.3.2) with step $\gm$. 
\item{1.4.1.}{  Let $\{\phi_{ab}\}$ be a meromorphic $\gm$-closed 
  multiplicative $2$-form
on $\h^*$.  Set
$$
R(\la) \, \mapsto \, 
\sum_{a = 1}^N \, 
\al_{aa}(\la)\,E_{aa}\otimes E_{aa}\,+
\sum_{a \neq b } \,\phi_{ab}(\la)\, 
\al_{ab}(\la)\,E_{aa}\otimes E_{bb}\,+
\,\sum_{a\neq b } \,\bt_{ab}(\la)\,E_{ba}\otimes E_{ab}.
$$
}
\item{1.4.2.}{ Let the symmetric group 
$S_N$ , the Weyl group of $gl_N$,
act on $\h^*$ and $V$ by permutation of coordinates.
For any permutation $\sigma \in S_N$, set
$$
R(\la) \, \mapsto \, (\sigma\otimes\sigma)\, R(\sigma^{-1} \cdot \lambda)\,
(\sigma^{-1} \otimes \sigma^{-1}) \,.
$$
}
\item{1.4.3.}{ For a nonzero complex number $c$, set
$$
R( \la ) \, \mapsto \, c\, R(  \la )\,.
$$
}
\item{ 1.4.4.}{ For a nonzero complex number $c$
 and an element $\mu \in \h^*$, set
$$
R( \la ) \, \mapsto \, R( c\, \la \, +\,\mu )\,.
$$
}

It is clear that any gauge transformation of types (1.4.2)-(1.4.3)
 transforms a 
quantum dynamical R-matrix with step $\gm$ to
a quantum dynamical R-matrix with step $\gm$.
Any gauge transformation of type (1.4.4) transforms a 
quantum dynamical R-matrix with step $\gm$ to
a quantum dynamical R-matrix with step $ \gm / c$.
In all cases, if the R-matrix is of Hecke type,
then the transformed matrix is of Hecke type.
If the transformation is of type (1.4.3) and the Hecke parameters
of the R-matrix are $p$ and $q$, then the Hecke parameters
of the transformed matrix are $c p$ and $c q$.

\proclaim{Theorem 1.1}
Any gauge transformation of type (1.4.1) transforms a 
quantum dynamical R-matrix with step $\gm$ to
a quantum dynamical R-matrix with step $\gm$.
 If the R-matrix is of Hecke type,
then the transformed matrix is of Hecke type with the same parameters.
\endproclaim

Theorem 1.1 is proved in Section 1.9.

Two R-matrices $R \, : \, \h^*  \to \End (V \otimes V)$ and
$R' \, : \, \h^*  \to \End (V \otimes V)$
 will be called {\it equivalent}
if one of them can be transformed into another by a sequence
of gauge transformations.

\subhead 1.5. Classification of quantum dynamical R-matrices of  Hecke type
with parameters $p, q$ such that $  q=p$
\endsubhead

If Hecke parameters satisfy $p=q$, then the Hecke equation (1.3.6) can be written
as
$$
R^{21}(\la)\,R(\la)\,=\, q^2\,\Id\,.
$$

 Let $X \subset \{1, ..., N\}$ be a subset. Say that $X$ {\it is decomposed
into disjoint
intervals},  $X = X_1 \cup ... \cup X_n$, if
 every $X_k$ has the form $\{ a_k, a_k + 1, ..., b_k\}$ and
$a_{k+1} > b_k $ for $k=1,..., n-1$. 

 A meromorphic function $\mu (\la)$ will be called  
{\it $\gm$-quasiconstant} if $\dl_a \mu = 0$ for all $a$.
Fix a $\gm$-quasiconstant $ \mu : \h^* \to \h^*$ with $\gm=1$. Define
scalar meromorphic $\gm$-quasiconstant functions $\mu_{ab}:\h^*
\to \C$ by $\mu_{ab}(\la)=x_a(\mu(\la))-x_b(\mu(\la))$.
Let $\la_{ab}$ denote $\la_a-\la_b$.

Define $ R_{\cup X_k}:\h^* \to \End (V\otimes V)$
by
$$
R_{\cup X_k}(\la) \, = \, 
\sum_{a, b = 1}^N \, 
E_{aa}\otimes E_{bb}\,+
\sum_{k=1}^n\,
\sum_{a, b \in X_k\, a\neq b} \,
{1\over \la_{ab}-\mu_{ab}(\la)}
\,(\,E_{aa}\otimes E_{bb}\,+\,E_{ba}\otimes E_{ab}\,)\,.
\tag 1.5.1
$$
%where $\la_{ab}$ denotes $\la_a - \la_b$.

\proclaim{Theorem 1.2}
\item{1.} {For every $X \subset \{1, ..., N\}$ , the R-matrix 
$R_{\cup X_k}$ defined by (1.5.1)
is a quantum dynamical
R-matrix of Hecke type with parameters $p=1, \,q=1$ and step $\gm=1$.
}
\item{2.} {Every quantum dynamical R-matrix of Hecke type with parameters
$p, q$, such that $p=q$, is equivalent to one of the matrices (1.5.1).
}
\endproclaim
Theorem 1.2 is proved in Section 1.11.

\subhead 1.6. Classification of quantum dynamical R-matrices of  Hecke type
with parameters $p, q$ such that $q \neq p$
\endsubhead

%where every $X_k$ has the form $\{ a_k, a_k + 1, ..., b_k\}$ and
%$a_{k+1} > b_k + 1$ for $k=1,..., n-1$.
 
Assume that for any $a, b,\, a \neq b,$ a $\gm$-quasiconstant 
$\mu_{ab}:\h^* \to \C$ is given.
We say that this collection of quasiconstants is {\it multiplicative}
if
\item{1.6.1.}{ For any $a,b$, we have 
$$
\mu_{ab}(\la)\, \mu_{ba}(\la)\,=\,1\,.
$$
}
\item{1.6.2.}{For any $a,b,c$, we have
$$
\mu_{ac}(\la)\,=\, \mu_{ab}(\la)\,\mu_{bc}(\la)\,.
%\prod_{i=a}^{b-1}\,\mu_{i,i+1}(\la)\,.
$$
}

Fix a multiplicative family of $\gm$-quasiconstants with $\gm=1$.

Fix a complex number $\epe$ such that $e^\epe \neq 1$.
Let $X \subset \{1, ..., N\}$ be a subset, $X = X_1 \cup ... \cup X_n$
its decomposition  into disjoint intervals.

For any $a,b \in \{ 1, ..., N\}$, $a\neq b$, we shall introduce  functions
$\al_{ab}, \bt_{ab} \,:\, \h^* \to \C$. We shall introduce functions $\bt_{ab}$
and then set $\al_{ab} = e^\epe + \bt_{ab}$.

If  $a,b \in X_k$ for some k, then we set
$$
\bt_{ab}(\la) \,=\, { e^\epe\,-\,1 \over \,\mu_{ab}(\la)
e^{\epe \la_{ab}} \,-\,1 }\,.
\tag 1.6.3
$$
Otherwise we 
set $\bt_{ab}(\la) = 0$ , if $a<b$, and $\bt_{ab}(\la) = 1-e^\epe$,
if $a>b$.

Define $ R_{\cup X_k}:\h^* \to \End (V\otimes V)$
by
$$
R_{\cup X_k,\epe}(\la) \, = \, 
\sum_{a = 1}^N \, 
E_{aa}\otimes E_{aa}\,+
\sum_{ a\neq b} \,
\al_{ab}(\la)\,E_{aa}\otimes E_{bb}\,+\,\sum_{ a\neq b} \,
\bt_{ab}(\la)\,E_{ba}\otimes E_{ab}\,.
\tag 1.6.4
$$

\proclaim{Theorem 1.3}
\item{1.} {For every $X \subset \{1, ..., N\}$ , the R-matrix 
$R_{\cup X_k,\epe}$ defined by (1.6.4)
 is a quantum dynamical
R-matrix of Hecke type with parameters $p=1,\, q=e^\epe$ and step $\gm=1$.
}
\item{2.} {Every quantum dynamical R-matrix of Hecke type with parameters
$p, q$ such that $ q\neq p$ is equivalent to one of the matrices (1.6.4).
}
\endproclaim
Theorem 1.3 is proved in Section 1.12.

\subhead 1.7.
Quantization of classical dynamical r-matrices of $gl_N$ type
\endsubhead

Let $V$ be the $N$ dimensional $\h$-module considered in Section 1.3.
Let $ r:\h^* \to \End (V\otimes V)$ be a zero weight meromorphic function
satisfying CDYB (1.2.2).  Assume 
 that $r$ satisfies  {\it
the unitarity condition, }
$$
r(\la) + r^{21}(\la)\, = \,\epe \,P \,+\,\dl\, \Id
\tag 1.7.1
$$
for some constants $\epe, \dl \in \C$ and all $\la$.
The constant $\epe$ is called {\it the coupling constant},
the constant $\dl$ is called {\it the secondary coupling constant.}
The zero weight condition implies that $r$ has the form
$$
r(\la)\,=\, \sum_{a,b = 1}^N \,\al_{ab}(\la)\,E_{aa}\otimes E_{bb}\,+
\,\sum_{a\neq b } \,\bt_{ab}(\la)\,E_{ab}\otimes E_{ba}\,.
\tag 1.7.2
$$

We recall a classification of such r-matrices. First we introduce
gauge transformations of classical dynamical r-matrices.

\item{1.7.3.}{  Let $\psi \, = \, \sum_{a,b} \psi_{ab}(\la) dx_a \wedge dx_b$ be
a closed meromorphic differential 
2-form on $\h^*$ \newline ( and the notion of a closed 
differential  form has the standard meaning).
Set
$$
r(\la) \, \mapsto \, r(\la) \,+\,
\sum_{a \neq b}^N  \psi_{ab}(\la)\, E_{aa} \otimes E_{bb} \,.
$$
}
\item{1.7.4.}{ For $ \mu \in \h^*$, set 
$$
r(\la) \, \mapsto \,r(\la + \mu)\,.
$$
}
\item{1.7.5.}{ Let the symmetric group 
$S_N$ 
act on $\h^*$ and $V$ by permutation of coordinates.
For any permutation $\sigma \in S_N$, set
$$
r(\la) \, \mapsto \, (\sigma\otimes\sigma)\, 
r(\sigma^{-1} \cdot \lambda)\,
(\sigma^{-1} \otimes \sigma^{-1}) \,.
$$
}
\item{1.7.6.}{ For a nonzero complex number $c$, set 
$$
r(\la) \, \mapsto \,  c\,r(c \la ).
$$
}
\item{1.7.7.}{ For a nonzero complex number $c$, set 
$$
r(\la) \, \mapsto \,  r( \la )\,+\, c\,\Id.
$$
}

%Notice that the first three transformations do not change the coupling
%constant and the fourth transformation multiplies it by $c$.

Any gauge transformation transforms a classical dynamical r-matrix 
to a classical dynamical r-matrix [EV]. 
Two classical dynamical r-matrices $r(\la)$ and  $r'(\la)$ 
will be called {\it equivalent}
if one of them can be transformed into another by a sequence
of gauge transformations.

The gauge transformations of quantum dynamical R-matrices described in
Section 1.4 are analogs of
gauge transformations  of classical dynamical r-matrices.

\subhead Classification of r-matrices with zero coupling constant, $\epe = 0$
\endsubhead

Let $X \subset \{1, ..., N\}$ be a subset, $X = X_1 \cup ... \cup X_n$
its decomposition  into disjoint intervals. 
%where every $X_k$ has the form $\{ a_k, a_k + 1, ..., b_k\}$ and
%$a_{k+1} > b_k + 1$ for $k=1,..., n-1$.

Define
a map  $ r:\h^* \to \End (V\otimes V)$
by
$$
r_{\cup X_k}(\la) \, = \,   \sum_{k=1}^n \sum_{a,b \in X_k\, a\neq b}\,
{1\over \la_{ba}}\,
E_{ba} \otimes  E_{ab}\,.
\tag 1.7.8
$$

\proclaim {Theorem 1.4  }

\item{1.}{ For any $X$ and its decomposition  $X = X_1 \cup ... \cup X_n$
into disjoint intervals,  the function $r_{\cup X_k}$ defined by (1.7.8)
is a classical dynamical r-matrix with zero coupling constant.
}
\item {2.}{ Any classical dynamical r-matrix 
 $ r:\h^* \to \End (V\otimes V)$
with  zero
coupling constant is equivalent to one of the matrices (1.7.8).
}
\endproclaim
Theorem 1.4 follows from [EV].

\subhead Classification of r-matrices with nonzero coupling constant,
 $\epe \neq 0$
\endsubhead

Let $X \subset \{1, ..., N\}$ be a subset, $X = X_1 \cup ... \cup X_n$
its decomposition  into disjoint intervals. 
%where every $X_k$ has the form $\{ a_k, a_k + 1, ..., b_k\}$ and
%$a_{k+1} > b_k + 1$ for $k=1,..., n-1$.
 
For any $a,b \in \{ 1, ..., N\}$, $a\neq b$, we introduce  functions
$ \bt_{ab} \,:\, \h^* \to \C$. 
If  $a,b \in X_k$ for some k, then we set
$$
\bt_{ab} (\la) \,=\, 
 \coth \, ( \, \la_{ba})\, .
$$
Otherwise we set $\bt_{ab} (\la)\,=\,- 1$, if $a<b$ , 
and $\bt_{ab}(\la) \,=\, \,1$, if $a>b$.

Define $ r_{\cup X_k}:\h^* \to \End (V\otimes V)$
by
$$
r_{\cup X_k}(\la) \, = \, 
P \,+\,
\sum_{ a\neq b} \,
\bt_{ab}(\la)\,E_{ba}\otimes E_{ab}\,.
\tag 1.7.9
$$

\proclaim{Theorem 1.5  }
\item{1.} {For every $X \subset \{1, ..., N\}$ 
and its decomposition  $X = X_1 \cup ... \cup X_n$
into disjoint intervals,   the function $r_{\cup X_k}$ defined by (1.7.9)
 is a classical dynamical
r-matrix with nonzero coupling constant $\epe=2$ and the secondary
coupling constant $\dl =0$.
}
\item{2.} {Every classical dynamical r-matrix $ r:\h^* \to \End (V\otimes V)$
with nonzero coupling constant  is equivalent to one of the matrices (1.7.9).
}
\endproclaim
Theorem 1.5 follows from [EV].

\proclaim{Theorem 1.6}
\item{1.} {Every classical dynamical r-matrix $r$ with zero coupling constant, 
holomorphic on an open polydisc $U\subset \h^*$,
can be quantized to a quantum dynamical R-matrix 
$R_\gamma$ on $U$, of Hecke type with parameters
$p, q$ such that $p=q$.
}
\item{2.} {Every classical dynamical r-matrix $r$ 
with nonzero coupling constant,
holomorphic on an open polydisc $U\subset \h^*$,
can be
quantized to a quantum dynamical R-matrix $R_\gamma$ on $U$, 
of Hecke type with parameters
$p, q$ such that $p\neq q$.
}

\endproclaim

\demo{Proof}
The R-matrix
$$
R_{\cup X_k}(\la, \gm) \, = \, 
\sum_{a, b = 1}^N \, 
E_{aa}\otimes E_{bb}\,+
\sum_{k=1}^n\,
\sum_{a, b \in X_k\, a\neq b} \, 
{\gm \over \la_{ab}}\,(\,E_{aa}\otimes E_{bb}\,+\,E_{ba}\otimes E_{ab}\,)
$$
is a quantum dynamical R-matrix of Hecke type with parameters $p=q=1$ and step $\gm$.
Its quasiclassical limit is
$$
r'(\la) \, = \, 
\sum_{k=1}^n\,
\sum_{a, b \in X_k\, a\neq b} \, 
{- 1 \over \la_{ab}}\,(\,E_{aa}\otimes E_{bb}\,+\,E_{ab}\otimes E_{ab}\,).
$$
Making the gauge transformation (1.7.3) corresponding to the closed form
\newline
$
 \sum_k\sum_{a,b\in X_k,a<b}\,  \la_{ab}^{-1}\, dx_a \wedge dx_b\,,
$
we get the r-matrix $r_{\cup X_k}$ 
 defined by (1.7.8) with $\mu_{ab}=0$ for all $a,b$.
This remark and Lemma 1.1 easily imply the first statement of the Theorem.
The second statement is proved analogously. $\square$

\enddemo

\subhead 1.8. Quantum dynamical Yang-Baxter equation in coordinates
\endsubhead 

Consider a quantum dynamical R-matrix $R(\la)$ of  form
(1.3.2). Assume that the matrix is of Hecke type, with
 step $\gm =1$ and Hecke parameters
$p=1$ and $q$. Any R-matrix can be reduced to such an R-matrix
 by  gauge transformations of 
types (1.4.3) and (1.4.4).

%Applying  gauge transformations of types (1.4.4) and (1.4.5)
%we can make the step $\gm =1$ and the Hecke parameters $p=1$
%and $q$. We will keep these assumptions in this section.

The Hecke property implies that $\al_{aa}=1$
and hence the matrix has the form 
$$
R(\la)\,=\, \sum_{a = 1}^N \,E_{aa}\otimes E_{aa}\,+
\sum_{a\neq b } \,\al_{ab}(\la)\,E_{aa}\otimes E_{bb}\,+
\,\sum_{a\neq b } \,\bt_{ab}(\la)\,E_{ba}\otimes E_{ab}\,.
\tag 1.8.1
$$
The Hecke property also implies that for every $a,c \in
\{1, ..., N\},\, a\neq c$, we have
$$
\align
\bt_{ac}(\la )\,+ \bt_{ca}(\la )\, =\, 1\,-\,q,
\tag 1.8.2
\\
\bt_{ac}(\la)\,\bt_{ca}(\la)\,-\,
\al_{ac}(\la)\,\al_{ca}(\la)\,=\, - q,
\tag 1.8.3
\endalign
$$
this is the trace and the determinant of $R^\vee$ restricted to $V_{ac}$.

Applying both sides of the QDYB equation (1.1.3) to a basis vector $v_a\otimes
v_a \otimes v_c \in V^{\otimes 3},\, a\neq c,$ we get  equations
$$
\align
\al_{ca}(\la-\om_a )\, \bt_{ac}(\la )\,\al_{ac}(\la-\om_a)\,+\,
\bt_{ac}(\la-\om_a )^2\,=\,\bt_{ac}(\la-\om_a )\,,
\tag 1.8.4
\\
\bt_{ca}(\la-\om_a)\,\bt_{ac}(\la)\,\al_{ac}(\la-\om_a)\,+\,
\al_{ac}(\la-\om_a)\,\bt_{ac}(\la-\om_a)\,=\,
\tag 1.8.5
\\
\bt_{ac}(\la)\,\al_{ac}(\la-\om_a)\,.
\endalign
$$
Applying both sides of the QDYB equation (1.1.3) to a basis vector $v_a\otimes
v_b \otimes v_c \in V^{\otimes 3}$ with pairwise distinct
$a,b,c$ we get  equations
$$
\align
\al_{ab}(\la-\om_c)\,\al_{ac}(\la)\,\al_{bc}(\la-&\om_a)\,=\,
\al_{bc}(\la)\,\al_{ac}(\la-\om_b)\,\al_{ab}(\la)\,,
\tag 1.8.6
\\
\al_{ac}(\la-\om_b)\,\al_{ab}(\la)\,\bt_{bc}(\la-&\om_a)\,=\,
\bt_{bc}(\la)\,\al_{ac}(\la-\om_b)\,\al_{ab}(\la)\,,
\tag 1.8.7
\\
\bt_{ab}(\la-\om_c)\,\al_{ac}(\la)\,\al_{bc}(\la-&\om_a)\,=\,
\al_{ac}(\la)\,\al_{bc}(\la-\om_a)\,\bt_{ab}(\la)\,,
\tag 1.8.8
\\
\bt_{cb}(\la-\om_a)\,\bt_{ac}(\la)\,\al_{bc}(\la-&\om_a)\,+\,
\al_{bc}(\la-\om_a)\,\bt_{ab}(\la)\,\bt_{bc}(\la-\om_a)\,=\,
\tag 1.8.9
\\
&\,\bt_{ac}(\la)\,\al_{bc}(\la-\om_a)\,\bt_{ab}(\la)\,,
\\
\endalign
$$
$$
\align
\al_{cb}(\la-\om_a)\,\bt_{ac}(\la)\,\al_{bc}(\la-\om_a)\,+\,
\bt_{bc}(\la-\om_a)\,\bt_{ab}(\la)\,\bt_{bc}(\la-\om_a)\,=\,
\tag 1.8.10
\\
\al_{ba}(\la)\,\bt_{ac}(\la-\om_b)\,\al_{ab}(\la)\,+\,
\bt_{ab}(\la)\,\bt_{bc}(\la-\om_a)\,\bt_{ab}(\la)\,,
\endalign
$$
$$
\align
\bt_{ac}(\la-\om_b)&\,\al_{ab}(\la)\,\bt_{bc}(\la-\om_a)\,=
\tag 1.8.11
\\
&\bt_{ba}(\la)\,\bt_{ac}(\la-\om_b)\,\al_{ab}(\la)\,+\,
\al_{ab}(\la)\,\bt_{bc}(\la-\om_a)\,\bt_{ab}(\la)\,.
\endalign
$$

\proclaim{Lemma 1.2}
For any $a,c, \,a\neq c$, the functions $\al_{ac}(\la)$
and $q+\bt_{ac}(\la)$ are not identically equal to zero.

\endproclaim
\demo{Proof}
If $\al_{ac} \equiv 0$, then equations $(1.8.2)_{ac},\,
(1.8.3)_{ac},\,(1.8.4)_{ac},$ and $(1.8.4)_{ca}$ give
 a  contradiction. Thus, $\al_{ac}$ and $\al_{ca}$ are not
identically equal to zero. Equations 
$(1.8.2)_{ac},\,
(1.8.3)_{ac}$ imply 
$$
\al_{ac}(\la)\,\al_{ca}(\la)\,=\,(q\,+ \,\bt_{ac}(\la))\,
\,(q\,+\,\bt_{ca}(\la))\,.
\tag 1.8.12
$$
The Lemma is proved. $\square$

\enddemo

\subhead 1.9. Proof of Theorem 1.1
\endsubhead

Let $\{\phi_{ab}\}$ be a  $\gm$-closed multiplicative 
2-form on $\h^*$. It is easy to see that equations (1.8.2)-(1.8.11)
are invariant with respect to the gauge transformation (1.4.1).
This proves Theorem 1.1. $\square$

\subhead 1.10. Relation $\al_{ac}\,=\, q\,+\,\bt_{ac}$
\endsubhead

Consider a quantum dynamical R-matrix $R(\la)$ of  form
(1.3.2). Assume that the matrix is of Hecke type with step 
 $\gm =1$ and Hecke parameters
$p=1$ and $q$. For any $a,c, \, a\neq c,$ set
$$
\phi_{ac}(\la)\,=\,{q\,+\,\bt_{ac}(\la) \over \al_{ac}(\la)}\,.
\tag 1.10.1
$$ 
\proclaim{Lemma 1.3}
The collection of functions $\phi\,=\,
\{\phi_{ac}\}$ is a $\gm$-closed multiplicative
2-form with $\gm=1$.

\endproclaim

\proclaim{Corollary 1.1} Apply to the R-matrix $R(\la)$ the gauge 
transformation (1.4.1) corresponding to the multiplicative 2-form
$\phi^{-1}$. Then the coefficients of the transformed matrix satisfy
the equation
$$
\al_{ac}\,=\, q\,+\,\bt_{ac}
\tag 1.10.2
$$
for all $a, c$.

\endproclaim

\demo{Proof of Lemma 1.3} Equation $\phi_{ac} \phi_{ca}=1$ follows from
(1.8.12). Equation $d_\gm \phi = 0$ is a direct corollary of
(1.8.6) and (1.8.7).
$\square$\enddemo

\subhead 1.11. Proof of Theorem 1.2
\endsubhead

Let $R(\la)$ be a quantum dynamical R-matrix
of Hecke type with paramaters $p, q$ such that $p=q$.
Using gauge transformations  (1.4.3) and (1.4.4) we can make
step $\gm=1$ and $p=q=1$. By Lemma 1.3 we may assume that
$\al_{ac}(\la)= 1+ \bt_{ac}(\la)$ for all $a\neq c$.
By (1.8.2) we have $\bt_{ac}(\la )=- \bt_{ca}(\la )$ for
all $ a \neq c$.

Fix $a,c,\, a\neq c$, and solve  equations 
 $(1.8.4)_{ac}, (1.8.5)_{ac}, (1.8.4)_{ca}, (1.8.5)_{ca}$.

\proclaim{Lemma 1.4}
Any solution $\bt_{ac}(\la),\,
\bt_{ca}(\la)$ 
of  equations $(1.8.4)_{ac}, (1.8.5)_{ac},$ 
$ (1.8.4)_{ca},$ $ (1.8.5)_{ca}$
has one of the following two forms.
\item{1.}{$\bt_{ac}=\bt_{ca}=0$.
}
\item{2.}{
$$
\bt_{ac}(\la)\,=\,{1\over \la_{ac}-\mu_{ac}},\qquad
\bt_{ca}(\la)\,=\,{1\over \la_{ca}-\mu_{ca}}
\tag 1.11.1
$$
where $\mu_{ac}=-\mu_{ca}$  and 
$\mu_{ac}(\la)$ is a meromorphic function periodic with respect to shifts of
$\la$ by $\om_a$ and $\om_c$,
$\mu_{ac}(\la - \om_a)\,=\,\mu_{ac}(\la - \om_c)\,=
\,\mu_{ac}(\la )$.
}
\endproclaim

\demo{Proof} It is easy to see that $\bt_{ac}(\la)=
\bt_{ca}(\la) \equiv 0$ is a solution. 
Now assume that $\bt_{ac}=-\bt_{ca}\neq 0$. Then
 $(1.8.5)_{ac}$ gives
$$
{1\over \bt_{ac}(\la)}\,+
{1\over \bt_{ac}(\la-\om_a)}\,=\,1\,,
$$
and $(1.8.5)_{ca}$ gives
$$
{1\over \bt_{ac}(\la)}\,+
{1\over \bt_{ac}(\la-\om_c)}\,=\,-1\,.
$$
Let $\mu_{ac}(\la)=\la_{ac}-1/\bt_{ac}(\la)$. Then
$\mu_{ac}(\la-\om_a)=\mu_{ac}(\la)$ and
$\mu_{ac}(\la-\om_c)=\mu_{ac}(\la)$.
Hence
$$
\bt_{ac}(\la)\,=\,{1\over \la_{ac}-\mu_{ac}}
$$
where $\mu_{ac}(\la)$ is a meromorphic function periodic in $\om_a$ and $\om_c$.
Similarly,
$$
\bt_{ca}(\la)\,=\,{1\over \la_{ca}-\mu_{ca}}
$$
where $\mu_{ca}(\la)$ is a function periodic in $\om_a$ and $\om_c$.
We have $\mu_{ac}=-\mu_{ca}$ since
 $\bt_{ac}=-\bt_{ca}$. It is easy to see that these functions
$\bt_{ac}$ and $\bt_{ca}$  solve equations
$(1.8.4)_{ac}$ and $(1.8.4)_{ca}$.
The Lemma is proved. $\square$

\enddemo

Equation (1.8.7) shows that the function $\bt_{ac}(\la)$
and hence the function  $\mu_{ac}(\la)$ is periodic with respects to 
shifts of $\la$ by $\om_b$ for any $b$ different from $a$ and $c$.

Consider equation $(1.8.9)_{abc}$ on functions 
$\bt_{ab}(\la),$ $\bt_{bc}(\la)$, $\bt_{ac}(\la)$. It is easy to see that
if one of these three
functions is identically equal to zero, then there is another
function in this triple which is  identically equal to zero.

Introduce a relation on the set $\{1, ..., N \}$. For any
$a \in \{1, ..., N \}$, let $a$ be related to $a$. For any
$a,b \in \{1, ..., N \}, a\neq b,$ let $a$ be related
to $b$ if the function
$\bt_{ab}(\la)$ is not identically equal to zero. It is easy to 
see that this is an equivalence relation. 

Let $Y \subset \{1, ..., N \}$ be the union of all the equivalence classes
containing more than one element. Let $Y = Y_1 \cup ... \cup Y_n$ be
its decomposition  into equivalence classes.

If pairwise distinct $a,b,c \in \{1, ..., N \}$ do not belong to the same equivalence
class, then at least two of the three functions
$\bt_{ab}(\la),$ $\bt_{bc}(\la)$, $\bt_{ac}(\la)$ are identically equal to zero.
Hence this triple of functions satisfies equation $(1.8.9)_{abc}$. 
If all three elements $a, \,b,\,c$  belong to the same equivalence class, then
equation $(1.8.9)_{abc}$  takes the form
$$
\,{1\over \la_{cb}-\mu_{cb}}\,{1\over \la_{ac}-\mu_{ac}}\,+\,
\,{1\over \la_{ab}-\mu_{ab}}\,{1\over \la_{bc}-\mu_{bc}}
\,=\,
\,{1\over \la_{ac}-\mu_{ac}}\,{1\over \la_{ab}-\mu_{ab}}\,.
$$
This implies that $\mu_{ac}(\la)=\mu_{ab}(\la)+
\mu_{bc}(\la).$ Therefore there exists a 1-quasiconstant meromorphic map
$\mu : \h^* \to \h^*$ such that $\mu_{ac}(\la)=x_a(\mu(\la))-
x_c(\mu(\la))$ for all $a, c$ such that $\mu_{ac}(\la)$ is not identically equal to zero.
It is easy to see that if the functions $\mu_{ab}(\la)$ have this property then equations
(1.8.8) and (1.8.10) are also satisfied.

Let $\sigma$ be a permutation of $\{1, ..., N \}$ which transforms
the set $Y$ and the decomposition $Y = Y_1 \cup ... \cup Y_n$ into a set
$X\subset \{1,..., N\}$ and its decomposition into disjoint intervals
$X = X_1 \cup ... \cup X_n$. Apply to the R-matrix $R(\la)$ the gauge transformation
(1.4.2) corresponding to the permutation $\sigma$. Then the transformed R-matrix will
have form (1.5.1) corresponding to the constructed decomposition
$X = X_1 \cup ... \cup X_n$.
Theorem 1.2 is proved. $\square$

\subhead 1.12. Proof of Theorem 1.3
\endsubhead

Let $R(\la)$ be a quantum dynamical R-matrix 
of Hecke type with paramaters $p, q$ such that $p\neq q$.
Using gauge transformations  (1.4.3) and (1.4.4) we can make
step $\gm=1$ and $p=1$.
Fix a number $\epe$ such that $q=e^\epe$.

By Lemma 1.3 we may assume that
$\al_{ac}(\la)= q + \bt_{ac}(\la)$ for all $a\neq c$.
By (1.8.2) we have $\bt_{ca}(\la ) = 1-q - \bt_{ac}(\la )$ for
all $ a \neq c$.

Fix $a,c,\, a\neq c$, and solve  equations 
$(1.8.4)_{ac}, (1.8.5)_{ac},$ $ (1.8.4)_{ca},$ $ (1.8.5)_{ca}$.

\proclaim{Lemma 1.5}
Any solution $\bt_{ac}(\la),\,
\bt_{ca}(\la)$ 
of  equations $(1.8.4)_{ac}, (1.8.5)_{ac},$ $ (1.8.4)_{ca}, $ $(1.8.5)_{ca}$
has one of the following two forms.
\item{1.}{$\bt_{ac}=0, \bt_{ca}= 1-q$ or
$\bt_{ca}=0, \bt_{ac}= 1-q$.
}
\item{2.}{
$$
\bt_{ac}(\la) \,=\, { e^\epe\,-\,1 \over \,\mu_{ac}(\la)
e^{\epe \la_{ac}} \,-\,1 }\,,\qquad
\bt_{ca}(\la) \,=\, { e^\epe\,-\,1 \over \,\mu_{ca}(\la)
e^{\epe \la_{ca}} \,-\,1 }\,
\tag 1.12.1
$$
 where $\mu_{ac}(\la) \mu_{ca}(\la) = 1$  and 
$\mu_{ac}(\la)$ is a meromorphic function periodic with respect to shifts of
$\la$ by $\om_a$ and $\om_c$,
$\mu_{ac}(\la - \om_a)\,=\,\mu_{ac}(\la - \om_c)\,=
\,\mu_{ac}(\la )$.
}

\endproclaim

%The first type of solutions will be called {\it degenerate.}

\demo{Proof} Equation $(1.8.4)_{ac}$ can be written
in the form
$$
(q+\bt_{ac}(\la-\om_a))\,(1-\bt_{ac}(\la-\om_a))\,\bt_{ac}(\la)\,=\,
(1-\bt_{ac}(\la-\om_a))\,\bt_{ac}(\la-\om_a)\,.
$$
Hence $\bt_{ac}(\la)\equiv 1$ or
$$
(q+\bt_{ac}(\la-\om_a))\,\bt_{ac}(\la)\,=\,\bt_{ac}(\la-\om_a)\,.
\tag 1.12.2
$$
The function $\bt_{ac}(\la)$ cannot be identically equal to 1.
In fact, if  $\bt_{ac}(\la)\equiv 1$, then equation $(1.8.4)_{ca}$ 
gives $ 0=-q(1+q)$ which is impossible since we always assume that 
$-q \neq p$.

Equation $(1.12.2)_{ac}$ has constant solutions $\bt_{ac}(\la)=0$ or
$\bt_{ac}(\la)=1-q$ which correspond to the first statement of the Lemma.
Now assume that $\bt_{ac}(\la)$ is not constant. Introduce a new meromorphic
 function
$y_{ac}(\la)= (\bt_{ac}(\la) + q - 1)/ \bt_{ac}(\la)$. It is easy to see that 
$y_{ac}(\la)\, y_{ca}(\la)\,=\,1$. Now equations 
$(1.12.2)_{ac}$, $(1.12.2)_{ca}$ 
can be written as
$$
y_{ac}(\la)\,=\,q\,y_{ac}(\la-\om_a), \qquad
y_{ac}(\la)\,=\,q^{-1}\,y_{ac}(\la-\om_c).
\tag 1.12.3
$$
Set $\mu_{ac}(\la) = y_{ac}(\la)\,e^{- \epe \la_{ac}}$. Then
the function $\mu_{ac}(\la)$ is periodic with respect to shifts of $\la$
by $\om_a$ and $\om_c$. We have 
$\mu_{ac}(\la)\,\mu_{ca}(\la)\,=\,1$. 
Returning to functions $\bt_{ac}(\la)$ and $\bt_{ca}(\la)$  we get the second type
of solutions. The Lemma is proved.  $\square$

\enddemo

Equation (1.8.7) shows that the function $\bt_{ac}(\la)$
and hence the function  $\mu_{ac}(\la)$ is periodic with respects to 
shifts of $\la$ by $\om_b$ for any $b$ different from $a$ and $c$.

If the function $\bt_{ac}(\la)$
has form (1.12.1), then we say that the function
$\mu_{ac}(\la)$ is {\it  finite}. If $\bt_{ac}(\la)=1-q$, then we say that
$\mu_{ac}(\la)=0$. If $\bt_{ac}(\la)=0$, then we say that
$\mu_{ac}(\la)=\infty$. 
 If $\mu_{ac}(\la)=0$, then $\mu_{ca}(\la)=\infty$.
 If $\mu_{ac}(\la)=\infty$, then $\mu_{ca}(\la)=0$.
% If $\mu_{ac}(\la)=0$ or $\mu_{ac}(\la)=\infty $,
%then we say that $\mu_{ac}(\la)$   is {\it degenerate}.

For  pairwise distinct $a, b, c$, we shall say that the equation
$$
\mu_{ab}(\la)\,\mu_{bc}(\la)\,=\,
\mu_{ac}(\la)\,
\tag 1.12.4
$$
holds if one of the following four conditions is satisfied.
\item{1.12.5.}{ All three functions $\mu_{ab}(\la),$ $\mu_{bc}(\la),$
$\mu_{ac}(\la),$ are finite, and satisfy (1.12.4).
}
\item{1.12.6.}{ $\mu_{ac}(\la)=\infty$ and at least one of the functions
$\mu_{ab}(\la)$, $\mu_{bc}(\la)$
is equal to $\infty$.
}
\item{1.12.7.}{$\mu_{ac}(\la)=0$ and at least one of the functions
$\mu_{ab}(\la)$, $\mu_{bc}(\la)$
is equal to $0$.
}
\item{1.12.8.}{ $\mu_{ac}(\la)$ is finite, one of the functions
$\mu_{ab}(\la)$, $\mu_{bc}(\la)$ is equal to zero and the other 
is equal to infinity.
}

\proclaim{Lemma 1.6}
 For any pairwise distinct $a, b, c$, equation (1.12.4) holds.

\endproclaim

The Lemma easily follows from equation (1.8.9).

Introduce
$$
Y=\{ (a,b)\,|\, (a,b) \in \{1,...,N\}, a\neq b,\, \mu_{ab}= \infty\},
\tag 1.12.9
%\\
%Y'=\{ (a,b)\,|\, (a,b) \in \{1,...,N\}, a\neq b,\, \mu_{ab}= 0\}.
%\endalign
$$
Then
\item{1.12.10.}{ If $(a,b) \in Y$ and $(b,c)\in Y$, then $(a,c)\in Y$.
}
\item{1.12.11.}{ If $(a,b)$ belongs to $Y$, then $(b,a)$
does not belong to $Y$.}

By Theorem 3.11 in [EV], there exists a permutation $\sigma$ of
numbers 
$\{1,...,N\}$ such that for the new order on $\{1,...,N\}$,
if $(a,b)\in Y$, then $a<b$. 
Apply to the R-matrix $R(\la)$
the gauge transformation (1.4.2) corresponding
to the permutation $\sigma$.
Then the set $Y$ defined by (1.12.9) for the transformed R-matrix
is such that if $(a,b)\in Y$, then $a<b$. From now on we denote
 by $R(\la)$ the transformed matrix.

Let $Z \,=\, \{(a,b) \,|\,a <b\}\,-\,Y$.
% and $U\,= \{(a,a+1)\,|\, (a,a+1) \in Z\}$.

\proclaim{Lemma 1.7}
\item{1.}{If $(a,b)$ belongs to $Z$, 
then all pairs $(c,c+1),\,
c=a,a+1,..., b-1,$
belong to $Z$.
}
\item{2.}{If for some $a,b,\,a<b,$ all pairs
 $(c,c+1)$ for $c=a,a+1,..., b-1$
belong to $Z$, then $(a,b)$ belongs to $Z$.
}
\endproclaim

Lemma 1.7 is a special case of Lemma 3.13 in [EV].

Consider the subset $X\subset \{1,...,N\}$ of all $a$
such that there exists $b$ with the property that
 $(a,b)$ or $(b,a)$ belongs to $Z$.

Introduce a relation on the set $X$. For any
$a \in X$, let $a$ be related to $a$. For any
$a,b \in X , a< b,$ let $a$ be related
to $b$ if $(a,b) \in Z$. Lemma 1.7 
implies that this relation is an
 equivalence relation. Let $X=X_1\cup ... \cup X_n$
be the decomposition of $X$ into equivalence classes.
Lemma 1.7 implies that
$X=X_1\cup ... \cup X_n$ is a decomposition into a union of 
disjoint intervals.
It is easy to see that the R-matrix $R(\la)$ has form 
(1.6.4) for the constructed decomposition 
$X=X_1\cup ... \cup X_n$. Theorem 1.3 is proved.
$\square$

\subhead 1.13. Quantum dynamical R-matrices as an extrapolation of constant
quantum R-matrices
\endsubhead 

Consider the vector representation $V$ of the quantum group $U_q(gl_N)$.
Then its R-matrix $\Cal R \in \End (V\otimes V)$ has the form,
$$
\Cal R \, = \, 
\sum_{a = 1}^N \, 
E_{aa}\otimes E_{aa}\,+
\sum_{ a\neq b} \,
\al_{ab}\,E_{aa}\otimes E_{bb}\,+\,\sum_{ a\neq b} \,
\bt_{ab}\,E_{ba}\otimes E_{ab}\,
\tag 1.13.1
$$
where the numbers $\al_{ab}, \bt_{ab}$ are defined as follows:
$\al_{ab} = q,\, \bt_{ab} = 0$ if $a<b$ and
$\al_{ab} = 1,\, \bt_{ab} = 1-q$ if $a>b$.
The matrix $\Cal R$ is a constant solution of the quantum dynamical
Yang-Baxter equation
(1.1.3). 

For any permutation $\sigma$ of numbers $\{1,...,N\}$ we construct
a new constant solution, $\Cal R_\sigma$,
 of the quantum Yang-Baxter equation. $\Cal R_\sigma$ has form (1.13.1)
where the numbers $\al_{ab}, \bt_{ab}$ are defined by the rule:
$\al_{ab} = q,\, \bt_{ab} = 0$ if $\sigma(a)<\sigma(b)$ and
$\al_{ab} = 1,\, \bt_{ab} = 1-q$ if $\sigma(a)>\sigma(b)$.

Fix a complex number $\epe$ such that $e^\epe = q$. Consider the matrix
$$
R (\la) \, = \, 
\sum_{a = 1}^N \, 
E_{aa}\otimes E_{aa}\,+
\sum_{ a\neq b} \,
\al_{ab}(\la)\,E_{aa}\otimes E_{bb}\,+\,\sum_{ a\neq b} \,
\bt_{ab}(\la)\,E_{ba}\otimes E_{ab}\,
\tag 1.13.2
$$
where the functions $\al_{ac}(\la)$ and $\bt_{ac}(\la)$ are defined by
$$
\bt_{ab}(\la) \,=\, { e^\epe\,-\,1 \over 
e^{\epe \la_{ab}} \,-\,1 }\,, \qquad
\al_{ab} = e^\epe + \bt_{ab}\,.
$$
The matrix $R(\la)$ is the R-matrix of form (1.6.4) corresponding
to data $X=X_1=\{1,...,N\}$.

The R-matrix $R(\la)$ extrapolates 
the constant R-matrices $\{\Cal R_\sigma \}$ in the following sence.
Let $\rho=(N/2, (N-2)/2,...,-N/2) \in \h^*$. Let $\sigma(\rho)$
be the vector obtained from $\rho$ by permutation of 
coordinates by $\sigma$. Then
$$
\text{lim}_{t\to +\infty}\, R(\,{}\,{\,t\,\over \,\epe\,}\,{}\, \sigma (\rho){})\,=\,
\Cal R_\sigma\,.
\tag 1.13.3
$$

\head 2. Quantum Dynamical R-matrices with Spectral Parameter
\endhead

\subhead 2.1. Definition
\endsubhead

Let $\h$ be an abelian finite dimensional Lie algebra.
Let $V_i, \, i=1, 2, 3,$ be finite dimensional
diagonalizable  $\h$-modules,
$$
R_{V_iV_j}\,:\, \C\times \h^* \,\to\, \End (V_i \otimes V_j), \qquad 1\leq i < j \leq 3,
$$
meromorphic functions, $\gm $ a nonzero complex number.
The  equation in $\End (V_1 \otimes V_2 \otimes V_3)$,
$$
\align
R_{V_1V_2}^{12}(z_1 - & z_2, \la - \gm h^{(3)})\,
R_{V_1V_3}^{13}(z_1 - z_3, \la )\,
R_{V_2V_3}^{23}(z_2-z_3, \la - \gm h^{(1)})\,
\tag 2.1.1
\\
&=\,R_{V_2V_3}^{23}(z_2 - z_3, \la )\,
R_{V_1V_3}^{13}(z_1- z_3, \la - \gm h^{(2)})\,
R_{V_1V_2}^{12}(z_1, - z_2, \la)\,
\endalign
$$
is called  {\it the quantum dynamical Yang-Baxter
equation with spectral parameter and step} $ \gm$ (QDYB equation).
In what follows we will use a notation $z_{ij}= z_i-z_j$.

A function
$R_{V_iV_j} \, : \C \times \h^* \to \End (V_i \otimes V_j)$ is called  {\it a
function of zero
weight} if
$$
[R_{V_iV_j}(z, \la), h \otimes 1 + 1 \otimes h ] \, = \, 0
\tag 2.1.2
$$
for all $h \in \h$, $z \in \C,\, \la \in \h^*$.
A solution $\{ R_{V_iV_j} \}_{1\leq i< j\leq 3}$  of the
QDYB equation (2.1.1) is called
a solution of zero weight if
each of the functions is of zero weight.

If all the spaces $V_i$ are equal to a space $V$,  then we consider the QDYB equation on
one function $R \, : \h^* \to \End (V \otimes V)$,
$$
\align
R^{12}(z_{12},\,\la - \gm h^{(3)})\,
R^{13}(z_{13},\,\la )\, &
R^{23}(z_{23},\,\la - \gm h^{(1)})\,
\tag 2.1.3
\\
=
R^{23}(z_{23},\,\la )\,&
R^{13}(z_{13},\,\la - \gm h^{(2)})\,
R^{12}(z_{12},\,\la)\,.
\endalign
$$
 A zero weight function $R$ satisfying
the QDYB equation (2.1.3) is called  {\it a quantum dynamical R-matrix
with spectral parameter}. An R-matrix is called {\it unitary},
if it satisfies  the {\it unitarity condition}
$$
R(z, \la)\,R^{21}(-z, \la)\,=\,1\,.
\tag 2.1.4
$$

\subhead 2.2. Quantization and quasiclassical limit
\endsubhead

Let $x_1,...,x_N$ be a basis in $\h$. The basis defines a linear system of coordinates
on $\h^*$. For any $\la \in \h^*$, set $\la_i = x_i (\la)$, $i=1,..., N$.

Let $R_\gm: \C \times \h^* \to \End (V\otimes V)$
be a smooth family of solutions to equations (2.1.3) and (2.1.4) with step $\gm$
such that
$$
R_\gm(z, \la)\,=\, 1\,-\,\gm\,r(\la)\,+\, O(\gm^2).
\tag 2.2.1
$$
Then the function $r : \C\times \h^* \to \End (V\otimes V)$
satisfies the zero weight condition
$$
[r(z,\la), h \otimes 1 + 1 \otimes h ] \, = \, 0
\tag 2.2.2
$$
for all $h \in \h$, $z \in \C,\, \la \in \h^*$, the unitarity
condition 
$$
r(z,\la)\,+\,r^{21}(-z,\la)\,=\,0\,
\tag 2.2.3
$$
and  the classical dynamical Yang-Baxter equation  with spectral parameter (CDYB),
$$
\align
& \sum _{i=1}^N \, x_i^{(1)} {\partial r^{23}
\over \partial x_i}(z_{23}, \la) \, + \, 
\sum_{i=1}^N \, x_i^{(2)} {\partial r^{31}
\over \partial x_i}(z_{31}, \la)
\, + \, \sum _{i=1}^N \, x_i^{(3)} {\partial r^{12 }
\over \partial x_i}(z_{12}, \la)\, +
\tag 2.2.4
\\
& [r^{12}(z_{12}, \la),r^{13}(z_{13}, \la)]\,+ \, 
[r^{12}(z_{12}, \la),r^{23}(z_{23}, \la)]\,
+\,[r^{13}(z_{13}, \la),r^{23}(z_{23}, \la)] \,= \,0\, .
\endalign
$$

A function $r(z,\la)$ with properties (2.2.2)-(2.2.4) is called a
{\it classical dynamical
r-matrix with spectral parameter.}

The function $r$ in (2.2.1) is called  {\it the quasiclassical limit of} $R$,
and the function $R$ is called  {\it a quantization of} $r$.

Let $U\subset \h^*$ be an open set, and let
$R : \C\times U\to \End (V\o V)$ be a 
zero weight meromorphic function on $\C \times U$. We will say that 
$R$ is a quantum dynamical R-matrix 
with spectral parameter on $\C\times U$ if the QDYB equation with spectral paramater 
is satisfied for $R$ whenever it makes sense.

A classical dynamical r-matrix $r(z,\la)$ with spectral parameter 
on $\C\times U$ is called {\it quantizable}
if there exists a power series in $\gm$,
$$
R_\gm(z,\la)\,=\, 1\,-\,\gm\,r(z,\la)\,+\, \sum_{n=2}^\infty \gm^n r_n(z,\la)
\tag 2.2.5
$$
convergent for small 
$|\gamma|$ for any fixed $(z,\la)\in \C\times U$,  
such that $R_\gm(z,\la)$ is a quantum dynamical R-matrix on $\C\times U$ with 
spectral parameter and step $\gamma$.

\subhead 2.3. R-matrices of  $gl_N$ type
\endsubhead

Let $\h$ be an abelian Lie algebra of dimension $N$. Let $V$ be a diagonalizable
$\h$-module of the same dimension such that
its weights  $\omega_1, ... , \omega_N$
form a basis in $\h^*$. Let $x_1,...,x_N$ be the dual basis of $\h$.
Let $v_1, ..., v_N$ be an eigenbasis for $\h$ in $V$ such that $x_i v_j = \dl_{ij} v_j$.
Then the $\h$-module $V\otimes V$ has the weight decomposition,
$$
V\otimes V \,=\, \oplus_{a=1}^N V_{aa} \, \oplus
\oplus_{a < b} V_{ab}\,,
\tag 2.3.1
$$
where $V_{aa}\, =\, \C\,v_a\otimes v_a$ and $V_{ab}\, =\, \C\,v_a\otimes v_b \oplus
\C\, v_b \otimes v_a $ .

%We shall use the basis  $E_{ij}$ in $ \End (V \otimes V)$ by $E_{ij}v_k = \dl_{jk} v_i$.

A quantum dynamical R-matrix with spectral parameter,  $R : \C \times
\h^* \to \End (V\otimes V)$,
for these $\h$ and $V$ will be called an R-matrix of $gl_N$ type.

The zero weight condition implies that the R-matrix preserves the weight decomposition (2.3.1)
and has the form
$$
R(z,\la)\,=\, \sum_{a,b = 1}^N \,\al_{ab}(z,\la)\,E_{aa}\otimes E_{bb}\,+
\,\sum_{a\neq b } \,\bt_{ab}(z,\la)\,E_{ba}\otimes E_{ab}\,
\tag 2.3.2
$$
where $\al_{ab}, \bt_{ab} \,:\,\C\times \h^* \to \C$
are suitable meromorphic functions.

\subhead 2.4. Gauge transformations 
\endsubhead

Fix a nonzero complex number $\gm$.
Let $\psi :\h^* \to \C$ be a function. 
For any $a, b =1,...,N$, set 
$$
\align
\partial_a \psi (\la)\,& =\, \psi (\la) - \psi (\la - \om_a), 
\\
L_{ab} \psi(\la)\, &
=\,\partial_a \psi (\la) - \partial_b \psi (\la - \om_a)\,=\,
\psi (\la) - 2\psi (\la - \om_a) + \psi(\la- \om_a - \om_b ).
\endalign
$$

Introduce gauge transformations of 
quantum dynamical R-matrices,  $R: \C\times\h^* \to \End (V\otimes V)$,
of type (2.3.2) with step $\gm$.
\item{ 2.4.1}{
Let $\psi$ be a meromorphic function on $\h^*$. Set
$$
\align
R(z,\la) \, \mapsto \,
& % \sum_{a = 1}^N \,
%\al_{aa}(z,\la)\,E_{aa}\otimes E_{aa}\,+
\\
&\sum_{a , b =1 }^N \,e^{z \partial_{a}\partial_b \psi(\la)}\,
\al_{ab}(z,\la)\,E_{aa}\otimes E_{bb}\,+
\,\sum_{a\neq b } \,e^{z L_{ab}\psi(\la)}\,
\bt_{ab}(z,\la)\,E_{ba}\otimes E_{ab}.
\endalign
$$
}
\item{ 2.4.2}{
Let $\{\phi_{ab}\}$ be a meromorphic $\gm$-closed multiplicative 2-form on $\h^*$. Set
$$
\align
R(z,\la) \, \mapsto \,
& \sum_{a = 1}^N \,
\al_{aa}(z,\la)\,E_{aa}\otimes E_{aa}\,+
\\
&\sum_{a \neq b } \, \phi_{ab}(\la)\,
\al_{ab}(z,\la)\,E_{aa}\otimes E_{bb}\,+
\,\sum_{a\neq b } \,
\bt_{ab}(z,\la)\,E_{ba}\otimes E_{ab}.
\endalign
$$
}
\item{2.4.3}{ Let the symmetric group
$S_N$ 
act on $\h^*$ and $V$ by permutation of coordinates.
For any permutation $\sigma \in S_N$, set
$$
R(z, \la) \, \mapsto \, (\sigma\otimes\sigma)\, R(z,
\sigma^{-1} \cdot \lambda)\,
(\sigma^{-1} \otimes \sigma^{-1}) \,.
$$
}
\item{2.4.4}{ For a nonzero holomorphic scalar function $c(z)$, set
$$
R(z, \la ) \, \mapsto \, c(z)\, R(z,  \la )\,.
$$
}
\item{ 2.4.5}{ For nonzero complex number $b, c$  and an element $\mu \in \h^*$, set
$$
R(z, \la ) \, \mapsto \, R(bz,\, c\la + \mu )\,.
$$
}

It is clear that any gauge transformation of type (2.4.3)
 transforms a (unitary)
quantum dynamical R-matrix with spectral parameter and step $\gm$ to
a (unitary) quantum dynamical R-matrix with spectral parameter and step $\gm$.
Any gauge transformation of type (2.4.4)
 transforms a 
quantum dynamical R-matrix with spectral parameter and step $\gm$ to
a  quantum dynamical R-matrix with spectral parameter and step $\gm$.
If in addition we have $c(z)c(z^{-1})=1$, then the
gauge transformation of types (2.4.4)
 transforms a unitary
quantum dynamical R-matrix with spectral parameter and step $\gm$ to
a unitary quantum dynamical R-matrix with spectral parameter and step $\gm$.
Any gauge transformation of type (2.4.5) transforms a (unitary)
quantum dynamical R-matrix with spectral parameter and step $\gm$ to
a (unitary)
 quantum dynamical R-matrix with spectral parameter and step $ \gm / c$.

\proclaim{Theorem 2.1}
Any gauge transformation of type (2.4.1) or (2.4.2) transforms a
quantum dynamical R-matrix with spectral parameter and step $\gm$ to
a quantum dynamical R-matrix with spectral parameter and step $\gm$.
Moreover, if the initial quantum dynamical R-matrix is unitary,
then the transformed R-matrix is unitary.
\endproclaim

Theorem 2.1 is analogous to Theorem 1.1 and is also
proved by direct verification. Namely, in order to prove Theorem 2.1
it is enough to
write  the QDYB equation (2.1.3) in coordinates, as it was done for equation
(1.1.3) in Section 1.8, and then check that if functions $\al_{ab}(z,\al)$
and $\bt_{ab}(z,\al)$ form a solution of the coordinate equations, then
the transformed functions also form a solution.

Two R-matrices $R \, : \, \C\times\h^*  \to \End (V \otimes V)$ and
$R' \, : \, \C\times \h^*  \to \End (V \otimes V)$
 will be called {\it equivalent}
if one of them can be transformed into another by a sequence
of gauge transformations.

\subhead 2.5. Examples 
\endsubhead

{\bf The elliptic R-matrix.}

Fix a point $\tau$ in the upper half plane and a 
complex number $\gamma$. Let
$$
\theta(z, \tau)=-\!\sum_{j\in\Z+\frac12}e^{\pi ij^2\tau+2\pi ij(z+\frac12)}
$$
be Jacobi's first theta function.

Let $\h$ be the Cartan subalgebra of $gl_N$. It is the abelian
Lie algebra of diagonal  complex $N\times N$ matrices with the standard
basis $x_i
=\text{diag}(0,\dots,0,1_i,0,\dots,0)$, 
$i=1,\dots, N$. Its dual space $\h^*$ has the dual basis  
$\omega_i$.

The vector representation of $gl_N$ is $V=\C^N$
with the standard basis $v_1,\dots,v_N$, $x_iv_j=\dl_{ij}v_j$.

Let $R^{ell}_{\gm,\tau} (z,\lambda)\in\End(V\otimes V)$
be the $R$-matrix of the elliptic quantum
group $E_{\tau,\gamma/2}(sl_N)$, [F1-2, FV2]. It is a function of the spectral
parameter $z\in\C$ and an additional variable
$\lambda=(\lambda_1,\dots,\lambda_N)
\in\h^*$. It is a solution of the CDYB equation (2.1.3) and satisfies 
the unitarity condition
(2.1.4) [F1-2]. The formula for $R^{ell}_{\gm,\tau}$ is
$$
R^{ell}_{\gm,\tau} (z,\lambda)=\sum_{a=1}^NE_{aa}\otimes E_{aa}
+\sum_{a\neq b}\alpha(z,\lambda_{ab}) E_{aa}\otimes E_{bb}
+
\sum_{a\neq b}\beta(z,\lambda_{ab}) E_{ba}\otimes E_{ab},
\tag 2.5.1
$$
where $\la_{ab}=\la_a-\la_b$ and the functions $\alpha,\beta$ are  ratios of theta functions:
$$
\alpha(z,\lambda)=\frac{\theta(\lambda+\gamma,\tau)\theta(z,\tau)}
{\theta(\lambda ,\tau)\theta(z-\gamma,\tau)}\,,
\qquad
\beta(z,\lambda)=\frac{\theta(z-\lambda,\tau)\theta(\gamma ,\tau)}
{\theta(z-\gamma,\tau)\theta(\lambda,\tau)}\,.
\tag 2.5.2
$$
%The R-matrix $R^{ell}_\gm (z,\lambda)$ is invariant under the symmetric group
%$S_N$ (the Weyl group of $sl_N$),
%in the sense that for any permutation $\sigma$
%$$
%R^{ell}(z, \la) \, \mapsto \, (\sigma\otimes\sigma)\, R^{ell}(\sigma^{-1} \cdot \lambda)\,
%(\sigma^{-1} \otimes \sigma^{-1}) \,.
%$$
%where $S_N$ acts linearly on $\h$ and $V$ by permutation
%of coordinates.

{\bf Trigonometric R-matrices}.
 
Let $X \subset \{1, ..., N\}$ be a subset, $X = X_1 \cup ... \cup X_n$
its decomposition  into disjoint intervals.

For any $a,b \in \{ 1, ..., N\}$, $a\neq b$, we introduce  functions
$\al_{ab}, \bt_{ab} \,:\, \C\times \h^* \to \C$. 

If  $a,b \in X_k$ for some k, then we set
$$
\al_{ab}(z,\la)\,=\,{\sin (\la_{ab} + \gm)\, \sin (z) \over
\sin (\la_{ab})\, \sin (z-\gm)}, \qquad
\bt_{ab}(z,\la)\,=\,{\sin (z-\la_{ab})\, \sin (\gm) \over
\sin (\la_{ab})\, \sin (z-\gm)}.
\tag 2.5.3
$$
Otherwise we set
$$
\al_{ab}(z,\la)\,=\, e^{-i\gm}\, { \sin (z) \over
 \sin (z-\gm)}, \qquad
\bt_{ab}(z,\la)\,=\,-\,e^{iz}\, { \sin (\gm) \over
 \sin (z-\gm)}
\tag 2.5.4
$$
if $a<b$, and 
$$
\al_{ab}(z,\la)\,=\, e^{i\gm}\, { \sin (z) \over
 \sin (z-\gm)}, \qquad
\bt_{ab}(z,\la)\,=\,-\,e^{-iz}\, { \sin (\gm) \over
 \sin (z-\gm)}
\tag 2.5.5
$$
if $a>b$.

Define a function $ R^{trig}_{\cup X_k, \gm}:\C\times \h^* \to \End (V\otimes V)$
by
$$
R^{trig}_{\cup X_k, \gm}(z,\la) \, = \, 
\sum_{a = 1}^N \, 
E_{aa}\otimes E_{aa}\,+
\sum_{ a\neq b} \,
\al_{ab}(\la)\,E_{aa}\otimes E_{bb}\,+\,\sum_{ a\neq b} \,
\bt_{ab}(\la)\,E_{ba}\otimes E_{ab}\,
\tag 2.5.6
$$
where $\al_{ab}$ and $\bt_{ab}$ are defined by (2.5.3) - (2.5.5).

{\bf Rational R-matrices}.

Let $X \subset \{1, ..., N\}$ be a subset, $X = X_1 \cup ... \cup X_n$
its decomposition  into disjoint intervals.

For any $a,b \in \{ 1, ..., N\}$, $a\neq b$, we shall introduce  functions
$\al_{ab}, \bt_{ab} \,:\, \C\times \h^* \to \C$. 

If  $a,b \in X_k$ for some k, then we set
$$
\al_{ab}(z,\la)\,=\,{(\la_{ab} + \gm)\, z \over \la_{ab} \,(z-\gm)}, \qquad
\bt_{ab}(z,\la)\,=\,{(z- \la_{ab})\,\gm \over
\la_{ab}\,  (z-\gm)}.
\tag 2.5.7
$$
Otherwise we set
$$
\al_{ab}(z,\la)\,=\,  { z \over
  z-\gm }, \qquad
\bt_{ab}(z,\la)\,=\,-\, { \gm \over
 z-\gm}.
\tag 2.5.8
$$

Define a function $ R^{rat}_{\cup X_k,\gm}:\C~\times \h^* \to \End (V\otimes V)$
by
$$
R^{rat}_{\cup X_k, \gm}(z,\la) \, = \, 
\sum_{a = 1}^N \, 
E_{aa}\otimes E_{aa}\,+
\sum_{ a\neq b} \,
\al_{ab}(\la)\,E_{aa}\otimes E_{bb}\,+\,\sum_{ a\neq b} \,
\bt_{ab}(\la)\,E_{ba}\otimes E_{ab}\,.
\tag 2.5.9
$$
where $\al_{ab}$ and $\bt_{ab}$ are defined by (2.5.7) - (2.5.8).

\proclaim {Theorem 2.2}
For any subset $X\subset \{1,...,N\}$ and its 
decomposition $X = X_1 \cup ... \cup X_n$ into disjoint intervals, the 
functions
$ R^{trig}_{\cup X_k,\gm}$ and $ R^{rat}_{\cup X_k,\gm}$ are zero weight solutions
of the QDYB equation (2.1.3) satisfying 
 the unitarity condition (2.1.4).

\endproclaim

\demo{Proof} According to [F1-2] the elliptic R-matrix $R^{ell}_{\gm,\tau}$ is a
zero weight solution of the QDYB equation (2.1.3) satisfying 
 the unitarity condition (2.1.4).

If $q=e^{2\pi i \tau} \to 0$, then  $\theta (z) \sim 2 q^{1/8} \sin (\pi z)$.

These two facts show that the R-matrix $R^0(z,\la)$ of the form (2.3.2), with
$$
\al_{ab}(z,\la)\,=\,{\sin (\la_{ab} + \gm)\, \sin (z) \over
\sin (\la_{ab})\, \sin (z-\gm)}, \qquad
\bt_{ab}(z,\la)\,=\,{\sin (z-\la_{ab})\, \sin (\gm) \over
\sin (\la_{ab})\, \sin (z-\gm)}
$$
for all $a\neq b$ and $\al_{aa} \equiv 1$ for all $a$,
is a zero weight solution of the QDYB equation (2.1.3) satisfying 
the unitarity condition (2.1.4).

For any fixed $d \in \h^*$, the R-matrix $R^0(z,\la + d)$ is also
a zero weight solution of the QDYB equation (2.1.3) satisfying 
the unitarity condition (2.1.4).

Fix a subset $X\subset \{1,...,N\}$ and its 
decomposition $X = X_1 \cup ... \cup X_n$ into disjoint intervals. 
It is easy to see that there exists a sequence of elements
$d_i \in \h^*,\, i=1,2,...$  such that
 the R-matrix $ R^0(z, \la + d_i)$ has a limit when $i$
tends to infinity, and this limit is equal to $R^{trig}_{\cup X_k, \gm}(z,\la)$.
This observation shows that $R^{trig}_{\cup X_k, \gm}(z,\la)$ 
is a zero weight solution of the QDYB equation (2.1.3) satisfying 
the unitarity condition (2.1.4).

Rescale the R-matrix $R^{trig}_{\cup X_k, \gm}(z,\la)$ and consider
a matrix $R_\epe (z,\la)=R^{trig}_{\cup X_k, \epe \gm}(\epe z,\epe \la)$ where
$\epe$ is a new parameter. Let  $\gm, z, \la$ be fixed and let $\epe$ tends to
0. Then the limit
of $R_\epe (z, \la)$ is equal to $R^{rat}_{\cup X_k, \gm}(z,\la)$.
Hence, $R^{rat}_{\cup X_k, \gm}(z,\la)$
is a zero weight solution of the QDYB equation (2.1.3) satisfying 
the unitarity condition (2.1.4). Theorem 2.2 is proved.

$\square$\enddemo

\subhead 2.6.
Quantization of classical dynamical r-matrices of $gl_N$ type with spectral parameter
\endsubhead

Let $V$ be the $N$ dimensional $\h$-module considered in Section 2.3.
Let $ r: \C\times\h^* \to \End (V\otimes V)$ be a zero weight meromorphic function
satisfying CDYB (2.2.4) and  the unitarity condition (2.2.3). 

The zero weight condition implies that $r$ has the form
$$
r(z,\la)\,=\, \sum_{a,b = 1}^N \,\al_{ab}(z,\la)\,E_{aa}\otimes E_{bb}\,+
\,\sum_{a\neq b } \,\bt_{ab}(z, \la)\,E_{ab}\otimes E_{ba}\,.
\tag 2.6.1
$$

Assume that the function $r$ satisfies also the {\it residue condition}
$$
{\text{Res}}_{z=0} \, r(\la, z)\, = \, \epe \,P\, +\dl\, \Id\,.
$$
Here $P \in \End (V\otimes V)$ is the permutation of factors
and $\Id \in \End (V\otimes V)$ is the identity operator.
The complex numbers
$\epe$ and $\dl$ are called  {\it the coupling constant}
and  {\it the secondary coupling constant}, respectively. 
We always assume that the coupling constant $\epe$ is not equal to zero.

We recall a classification of such r-matrices. First we introduce
gauge transformations of classical dynamical r-matrices with spectral parameter.
\item{2.6.2}{  Let $\psi \, = \, \sum_{a,b} \psi_{ab}(\la) dx_a \wedge dx_b$ be
a closed meromorphic differential
2-form on $\h^*$.
Set
$$
r(z, \la) \, \mapsto \, r(z, \la) \,+\,
\sum_{a \neq b}  \psi_{ab}(\la)\, E_{aa} \otimes E_{bb} \,.
$$
}
\item{2.6.3}{For a holomorphic function
$\psi \, : \, \h^* \, \to \, \C$, set
$$
\align
r(z ,\la) \, \mapsto \,\sum_{a,b=1}^N  \,(\al_{ab}(z, \la)\,+ \,& z\,
{\partial^2 \psi \over \partial x_a \, \partial x_b}
(\la))\,
E_{aa} \otimes E_{bb} \, +
\,
\\
& \sum _{a\neq b} \,\bt_{ab} (z,\la)\,
e^{z ({\partial \psi \over \partial x_a }(\la)
- {\partial \psi \over \partial x_b }(\la))  }\,
E_{ab}\otimes E_{ba} \,.
\endalign
$$
}
\item{2.6.4}{ For $ \mu \in \h^*$, set
$$
r(z,\la) \, \mapsto \,r(z,\la + \mu)\,.
$$
}
\item{2.6.5}{ Let the symmetric group
$S_N$
act on $\h^*$ and $V$ by permutation of coordinates.
For any permutation $\sigma \in S_N$, set
$$
r(z,\la) \, \mapsto \, (\sigma\otimes\sigma)\,
r(z,\sigma^{-1} \cdot \lambda)\,
(\sigma^{-1} \otimes \sigma^{-1}) \,.
$$
}
\item{2.6.6}{ For a nonzero complex number $c$, set
$$
r(z,\la) \, \mapsto \,  c\,r(z, c \la ).
$$
}
\item{2.6.7}{ For an odd scalar meromorphic function $f(z)$,
$f(z)+f(-z)=0$, set 
$$
r(z, \la) \, \mapsto \,  r(z, \la )\,+\, f(z)\,\Id.
$$
}

%Notice that the first three transformations do not change the coupling
%constant and the fourth transformation multiplies it by $c$.

Any gauge transformation transforms a classical dynamical r-matrix
with spectral parameter
to a classical dynamical r-matrix with spectral parameter [EV].
Two classical dynamical r-matrices $r(z,\la)$ and  $r'(z,\la)$
will be called {\it equivalent}
if one of them can be transformed into another by a sequence
of gauge transformations.

The gauge transformations of quantum dynamical R-matrices with spectral
parameter described in Section 2.4 are analogs of the
gauge transformations  of classical dynamical r-matrices with spectral parameter.

\subhead Classification of the classical dynamical r-matrices with spectral
parameter
\endsubhead

\subhead
The elliptic r-matrix
\endsubhead

Fix a point $\tau$ in the upper half plane.
Introduce the functions
$$
\sigma_w(z) =
{\theta(w-z, \tau)\,\theta'(0, \tau)
\over
\theta(w, \tau)\,\theta(z, \tau)},
\qquad
\rho (z) =
{\theta'(z, \tau)
\over
\theta(z, \tau)},
$$
where $\theta'(z,\tau)=\frac{\d \theta(z,\tau)}{\d z}$. 
Set
$$
r^{ell}_\tau (z,\la)\,=\, \rho(z)\sum_{a=1}^N E_{aa}\otimes E_{aa}\,
+ \,\sum_{a\neq b } \sigma_{\la_{ba}}(z) E_{ab} \otimes E_{ba}\,.
\tag 2.6.8
$$
 For every $\tau \in \C$, Im $\tau > 0$, the
function $r^{ell}_\tau (z,\la)$ is a classical dynamical r-matrix
with spectral parameter $z$, coupling constant $\epe= 1$ and
secondary constant $\dl = 0$, [FW].

{\bf Trigonometric r-matrices}.

Let $X \subset \{1, ..., N\}$ be a subset, $X = X_1 \cup ... \cup X_n$
its decomposition  into disjoint intervals.

For any $a,b \in \{ 1, ..., N\}$, $a\neq b$, we introduce a  function
$\bt_{ab}\,:\, \C\oplus \h^* \to \C$.

If  $a,b \in X_k$ for some k, then we set
$$
\bt_{ab}(z,\la)\,=\,-\,{\sin (\la_{ab} + z) \over
\sin (\la_{ab})\, \sin (z)}.
$$
Otherwise we set
$$
\bt_{ab}(z,\la)\,=\, {e^{-iz} \over  \sin (z)}, \qquad
\text{for} \, a<b, \qquad
\bt_{ab}(z,\la)\,=\,{e^{iz} \over  \sin (z)} \qquad
\text{for} \, a>b.
$$

We introduce a trigonometric r-matrix 
$ r^{trig}_{\cup X_k, \gm}:\C\oplus  \h^* \to \End (V\otimes V)$
by
$$
r^{trig}_{\cup X_k}(z,\la)\,=\, \text{cotan}\,( z) \,\sum_{a=1}^N E_{aa}\otimes E_{aa}\,
+ \, \sum _{a \neq b} \, \bt_{ab} (z,\la)\,
\,E_{ab} \otimes  E_{ba},
\tag 2.6.9
$$
where $\text{cotan}\, (z)\,=\, \text{cos}\, (z)\,/ \text{sin} \,(z).$

\subhead  Rational r-matrices
\endsubhead

Let $X \subset \{1, ..., N\}$ be a subset, $X = X_1 \cup ... \cup X_n$
its decomposition  into disjoint intervals.
Set
$$
r^{rat}_{\cup X_k}(z,\la) \, = \,  { P  \over z} \,
+ \sum _{k=1}^n \sum_{a,b \in X_k,\, a\neq b}  \,
{1\over \la_{ab} } \,
E_{ab} \otimes  E_{ba}\, .
\tag 2.6.10
$$

\proclaim {Theorem 2.3 }

\item{1.} {For every subset $X \subset \{1, ..., N\}$ and its decomposition
 $X = X_1 \cup ... \cup X_n$  into disjoint intervals, the matrices
$r^{trig}_{\cup X_k}$ and $r^{rat}_{\cup X_k}$ are classical dynamical
r-matrices with spectral parameter.
}
\item{2.} {Every classical dynamical r-matrix $r:\C\times \h^* \to \End (V\otimes V)$
with nonzero coupling constant is equivalent to one of the matrices (2.6.8)-(2.6.10).
}

\endproclaim

Theorem 2.3 follows from [EV].

\proclaim {Theorem 2.4 }

Let $R(z,\la)$ be 
a unitary classical dynamical r-matrix $r$ with spectral parameter and nonzero
coupling constant, 
meromorphic on $\C\times U$, where $U$ is an open polydisc.
Assume that for any $\la\in U$ there exists $z\in\C$ 
such that $r$ is holomorphic at $(\la,z)$. 
Then $r$ 
can be quantized to a unitary
quantum dynamical R-matrix $R_\gamma$ on $\C\times U$ 
of $gl_N$ type.
Moreover, if a classical  dynamical r-matrix with spectral parameter and nonzero
coupling constant is equivalent to the elliptic  r-matrix (2.6.8)
(resp., a trigonometric r-matrix (2.6.9) or a rational r-matrix (2.6.10)),
then it has a quantization equivalent to the elliptic R-matrix (2.5.1)
(resp., a trigonometric R-matrix (2.5.6) or a rational R-matrix (2.5.9)).
\endproclaim

{\bf Remark.} It follows from \cite{EV} that the 
holomorphicity assumption of Theorem 2.4 holds for any meromorphic
unitary classical dynamical r-matrix with a nonzero coupling constant, 
as long as $U\subset Y(r)\subset \h^*$, where $Y(r)$ is
a dense open set. Thus, this assumption does not impose any 
significant restriction. 

\demo{Proof} We shall prove that if a classical dynamical r-matrix is equivalent to
the elliptic r-matrix (2.6.8), then it is quantizable to a quantum dynamical
R-matrix equivalent to the elliptic R-matrix (2.5.1). The other statements of the Theorem
are proved similarly.

Compute the quasiclassical limit of $R^{ell}_{\gm,\tau} (z,\la)$. For the functions
$\al(z,\la,\gm)$ and $\bt(z,\la,\gm)$ defined in (2.5.2), we have
$$
\text{lim}_{\gm \to 0}{\al(z,\la,\gm)- 1\over\gm}=
{\theta '(\la)\over \theta (\la)}+
{\theta '(z)\over \theta (z)},
\qquad
\text{lim}_{\gm \to 0}{\,\bt(z,\la,\gm)\,\over\gm}=
{\theta '(0) \theta (z-\la)\over \theta (\la)\theta (z)}.
$$
Hence
$$
R^{ell}_{\gm,\tau} (z,\la)\,=\, 1\, - \,\gm \, r(z,\la) \,+\, O(\gm^2)
$$
where 
$$
\align
r(z,\la)\,& =\,-\, \sum_{a\neq b}\,({\theta '(\la_{ab})\over \theta (\la_{ab})}+
{\theta '(z)\over \theta (z)})\,
E_{aa}\otimes E_{bb}\,-\,
\sum_{a\neq b} \, {\theta '(0) \theta (z-\la_{ab})\over \theta (\la_{ab})\theta (z)}
\, E_{ba}\otimes E_{ab}
\\
& =\,-\, \sum_{a\neq b}\,({\theta '(\la_{ab})\over \theta (\la_{ab})}+
{\theta '(z)\over \theta (z)})\,
E_{aa}\otimes E_{bb}\,+\,
\sum_{a\neq b} \, \sigma_{\la_{ba}}(z) 
\, E_{ab}\otimes E_{ba}.
\endalign
$$
Now applying to the r-matrix $r(z,\la)$
the transformation (2.6.1) coresponding to the closed differential 2-form
$$
\sum_{a\neq b}\, {\theta '(\la_{ab})\over \theta (\la_{ab})}\,
dx_a\wedge dx_b
$$
and then applying to the result the transformation (2.6.6) corresponding 
to the function $f(z)= \theta '(z)/ \theta (z)$ we get the matrix $r^{ell}_\tau
(z,\la)$
defined by (2.6.7). 
This remark and Lemma 1.1 easily imply the statement of the Theorem
concerning the elliptic r-matrix. Theorem 2.4 is proved.

$\square$\enddemo

{\bf Remark.} The elliptic 
quantum dynamical R-matrix (2.5.1) was invented by G. Felder [F1-2]
as a quantization of the classical dynamical r-matrix (2.6.8).

\subhead 2.7. Formal dynamical R-matrices and gauge fixing conditions
\endsubhead

Let $R_\gm(z,\la)=1-\gamma r(z,\la)+\sum_{n\ge 2}\gamma^nr_n(z,\la)$
be a power series in $\la$ and $\gm$, 
whose coefficients are meromorphic functions of $z$, 
taking values in $\End(V\o V)$. The series $R_\gm$ is called 
a formal quantum dynamical R-matrix
of $gl_N$ type with spectral parameter and step $\gm$ 
if it is of zero
weight and satisfies the quantum dynamical Yang-Baxter equation. 
In addition, $R_\gm$  is called unitary if it satisfies the unitarity condition
(2.1.4). In this section for brevity we will 
refer to formal quantum dynamical R-matrices of $gl_N$ type
with spectral parameter and step $\gm$ as ``formal dynamical R-matrices''. 
As we know, any such R-matrix has form (2.3.2).  

The theory of formal dynamical R-matrices is completely 
analogous to the theory of analytic dynamical 
R-matrices. In particular, one can define 
formal classical dynamical r-matrices
and formal gauge transformations in an obvious way. 
If $R_\gamma=1-\gamma r+...$ is a (unitary) formal dynamical R-matrix,
then $r$ is a (unitary) formal 
dynamical r-matrix.

An example of a formal dynamical R-matrix is 
the Taylor expansion of an analytic dynamical R-matrix 
$R_\gamma(z,\la)$ at a point $\gm=0,\la=\la_0$, such that $R$ is regular at
this point for  generic values of $z$. 

\proclaim{Proposition 2.1} Let $R_\gamma=1-\gamma r+...$ be a unitary 
formal dynamical R-matrix, and
$z_0\in \C$ a point where $R_\gm$ is regular. 
 Let $\al_{ab}$, $\beta_{ab}$ be the matrix coefficients 
of $R_\gm$, see (2.3.2). Then $R_\gamma$ can be transformed, by a sequence of
formal gauge transformations, 
to a unitary formal dynamical R-matrix satisfying 
the following conditions:

1) for every $a,b$, 
the ratio $\frac{\alpha_{ab}(z,\la-\gm\omega_c)}{\alpha_{ab}(z,\la)}$
is independent of $z$; 

2) for every $a<b$, \ \  $\al_{ab}(z_0,\la)=1$; 

3) the coefficient $\al_{11}(z,\la)$ is independent of $z$.
\endproclaim

{\it Proof.} 
The QDYB equation with spectral parameter implies the equation 
$$
\al_{ab}(u,\la-\gm\om_c)\,\al_{ac}(u+v,\la)\,\al_{bc}(v,\la-\gm\om_a)\,=\,
\al_{bc}(v,\la)\,\al_{ac}(u+v,\la-\gm\om_b)\,\al_{ab}(u,\la)\,
\tag 2.7.1
$$
for any $a,b,c$. Therefore, we have 
$$
\frac{\al_{ab}(u,\la-\gm\om_c)}{\al_{ab}(u,\la)}=H_{abc}(\la)e^{D_{abc}(\la)u},
\tag 2.7.2
$$
for suitable power series $H_{abc}(\la),\,D_{abc}(\la)$.

\proclaim{Lemma 2.1} There exists a formal power series $\psi (\la)$ 
such that $D_{abc}=\d_a \d_b\d_c\psi$. 
\endproclaim

\demo{Proof} From (2.7.1) it follows that 
$D_{abc}$ is symmetric. From (2.7.2) it follows that $\d_dD_{abc}$ is 
symmetric. The rest of the proof of the Lemma
follows from the basic theory of difference equations
with infinitesimal shift.
$\square$\enddemo

\proclaim{Corollary 2.1} Performing a gauge transformation (2.4.1), 
we can arrange $D=0$, i.e. condition 1. 
\endproclaim

   From now on we assume that $D=0$, i.e. 
$$
\frac{\al_{ab}(u,\la-\gm\om_c)}{\al_{ab}(u,\la)}=H_{abc}(\la).
\tag 2.7.3
$$
This implies that $\al_{ab}(u,\la)=\al^1_{ab}(u)\al^2_{ab}(\la)$, where 
$\al^i_{ab}$ are new functions. 

Consider the multiplicative 2-form 
$\phi$ defined by $\phi_{ab}(\la)=\al_{ab}(z_0,\la)$, $a<b$. It follows 
from (2.7.1) that $d_\gm \phi=0$. Therefore, by a gauge transformation 
of type (2.4.2) we can arrange $\phi=1$, i.e. condition 2. 

It remains to arrange condition 3. By (2.7.3), $\al_{11}(z,\la)=
f(z)g(\la)$ for a suitable formal power series $g(\la)$ and a
meromorphic function $f(z)$ 
such that $f(z)f(-z)=1$.
Applying transformation (2.4.4) with $c(z)=1/f(z)$, we get condition 3. 
The Proposition is proved.
$\square$

We will call conditions 1-3 the gauge fixing conditions. 

\subhead 2.8. Classification of unitary formal dynamical R-matrices
with elliptic quasiclassical limit  \endsubhead

We will say that a formal classical dynamical r-matrix
$r$ is of elliptic, trigonometric, or rational type
if it is gauge equivalent (by formal gauge transformations) to
an r-matrix of the form (2.6.8), (2.6.9),(2.6.10), respectively, 
expanded near a point $\la_0\in \h^*$. 
It follows from \cite{EV} that any formal classical 
dynamical r-matrix satisfying the residue condition with 
coupling constant $\epe\ne 0$ is of elliptic, trigonometric, or rational type.

\proclaim{Theorem 2.5} Let $R_\gamma=1-\gm r+O(\gm^2)$ be a unitary formal 
dynamical R-matrix whose quasiclassical limit $r$ is of the elliptic type.
Then there exist a point $\la_0\in \h^*$ and a power series
$\tau (\gm) = \tau_0+O(\gm)\in \C[[\gm]]$, $\text{Im}(\tau_0)>0$
such that the R-matrix $R_\gm$ can be transformed, by a sequence
of formal gauge transformations, into the Taylor series of
$R^{ell}_{\gm, \tau (\gm)}(z, \la-\la_0)$ where
$R^{ell}_{\gm, \tau }(z, \la)$ is the elliptic R-matrix (2.5.1).

\endproclaim

The proof of this Theorem occupies the next section. 

\subhead 2.9. Proof of Theorem 2.5 \endsubhead

Let $X^0$ be the space of unitary formal classical dynamical 
 r-matrices with spectral parameter 
and a nonzero coupling constant. Let 
$X_*^0$ be the subset of elements of $X^0$ which satisfy the  
following gauge fixing conditions: 
 
1c) $\frac{\d}{\d \la_c}\alpha_{ab}(z,\la)$ is independent of $z$;

2c) $\al_{ab}(z_0,\la)=0$, $a<b$;  

3c) $\al_{11}(z,\la)$
is independent of $z$

(these conditions are quasiclassical analogues of conditions 1-3 above).

According to the results of \cite{EV}, 
the space $X_*^0$ is a connected, finite-dimensional complex manifold 
(with singularities),
and any element of $X^0$ is gauge equivalent to an element of $X_*^0$. 
(i.e. $X_*^0$ is a ``cross-section''). Moreover, since $r\in X_*^0$ is 
of elliptic type, the manifold $X_*^0$ is smooth at $r$. 

Let $X$ be the space of unitary formal quantum dynamical 
 R-matrices with spectral parameter, and $X_*$ 
the subset of elements of $X$ satisfying the gauge fixing conditions 
1-3. 

As we have shown in Section 2.7, 
we can assume that our family $R_\gm$ is in 
$X_*$. In this case, $r\in X_*^0$. 

Now let us prove the statement of the Theorem 
modulo $\gamma^{m+1}$ by induction in $m$. 

For $m=1$, the Theorem is a tautology. 
Suppose we know the Theorem for $m=k\ge 2$, and want to prove it for 
$m=k+1$. 

We have a polynomial $R_k=1-\gm r+...+\gm^k r_k$ which satisfies the 
condition $R_k\in X_*$ modulo $\gamma^{k+1}$. 
We know that $R_k$ satisfies the conclusion of Theorem 2.5 modulo 
$\gamma^{k+1}$, i.e. is
of the form (2.5.1) modulo $\gamma^{k+1}$. 

Consider any extension of 
this polynomial to order $k+1$: $R_{k+1}=R_k+\gamma^{k+1}r_{k+1}$. 
The condition that $R_{k+1}\in X_*$ modulo $\gamma^{k+1}$ can be expressed as  
a nonhomogeneous linear equation with respect to $r_{k+1}$ having the form
$A\,r_{k+1}=s_{k+1}(r_k,...,r_2,r)$, where 
$A$ is a linear operator. 

The obvious, but crucial observation now is the following. 

\proclaim{Lemma 2.2} $\text{Ker}\, A=T_rX_*^0$, where $T_rX_*^0$ 
denotes the tangent space at the point $r$. 
\endproclaim

\demo{Proof} Indeed, it is easy to see by an explicit calculation that
the linear homogeneous equation $A\rho=0$ is nothing else 
but the equation for a tangent vector to $X_0^*$ at the point $r$. 
$\square$\enddemo

\proclaim{Corollary 2.2} 
The dimension of the space of solutions of $Ar_{k+1}=s_{k+1}$ is less than 
or equal to $K=\text{dim}(X_*^0)$.  
\endproclaim 

However, by Theorem 2.4, we already have a family 
of elements of $X_*$ with $K$ parameters -- the quantizations of 
elements of $X^0_*$. Therefore, using 
dimension arguments, we obtain that if $r_{k+1}$ satisfies $Ar_{k+1}=s_{k+1}$, 
then $R_{k+1}(\gamma)$ has to be in this 
K-parametric family, which completes the induction step. 

The Theorem is proved. 

$\square$

{\bf Remark.} If $r$ is not elliptic but rational 
or trigonometric, the result of Theorem 2.5 can be generalized, 
in the sense that formal dynamical R-matrices 
$R_\gamma=1-\gm r+...$ with rational or trigonometric $r$
can be explicitly classified up to gauge transformations   
by the same method as above. However,  
both the statement and the proof in this case are more delicate, as 
the manifold $X_0^*$ may now be singular at $r$, and it is necessary 
to describe carefully these singularities.
For simplicity one should first consider the case dim  
$V \, =2$, and then generalize to an arbitrary dimension. 
We are not giving this argument here. 

\head{3. Quantum Dynamical R-matrices and monoidal categories}\endhead

\vskip .05in

Let us briefly recall some standard notions 
of the category theory \cite{Mac, Kass}. 

Recall that a morphism $a:F\to G$ of two functors 
from a category $\Cal C$ to 
a category $\Cal C'$
is a choice of a morphism $a_X: F(X)\to G(X)$ for any object $X$ in 
$\Cal C$, such that for any two objects 
$X,Y\in \Cal C$ and any 
morphism $g:X\to Y$ we have $a_Y\circ F(g)=G(g)\circ a_X$.
An endomorphism of a functor is just a morphism of this functor into itself. 
 
Recall that a {\it monoidal category} 
is a category $\Cal C$ with a bifunctor $\o: \Cal C\times \Cal C\to 
\Cal C$ (i.e. a functor with respect to each factor), 
 called the tensor product, and an isomorphism of functors
$\Phi: (*\o *)\o *\to *\o (*\o *)$, called the associativity 
isomorphism, such that $\Phi$ satisfies 
the pentagon relation, and there exists a unit object $\bold 1\in \C$ 
with certain properties. A {\it braided monoidal category}
 is a monoidal category 
with a functorial isomorphism $\beta: \o\to \o^{op}$ called the
commutativity isomorphism, which satisfies the hexagon relations.
A braided category is called {\it symmetric} if $\beta^2=1$.  
A monoidal category will be called a {\it tensor category}
 if it has an additive 
structure $\oplus$, such that $\o$ is distributive with respect to $\oplus$.  
 
\subhead 3.1. The category of $\h$-vector spaces \endsubhead

Let $\h$ be a finite-dimensional commutative Lie algebra over $\C$. 
Let $M_{\h^*}$ denote the field of meromorphic functions on $\h^*$. 
Fix a complex number $\gamma$. 

Let $\V_\h$ denote the category whose objects are 
diagonalizable $\h$-modules,
and morphisms are defined by $\text{Hom}_{\V_\h}(X,Y)=
\text{Hom}_\h(X,Y\o_\C M_{\h^*})$. 

 Let $W\o *$ be the functor 
of multiplication by $W$. 
For any $W\in V_\h$ and $f\in \End_{\Vh}(W)$, 
define $f(*-\gamma h^{(2)})\in \End(W\o *)$ by the formula
$$
f_V(\la-\gamma h^{(2)})(w\o v)=f_V(\la-\gamma \mu)w\o v,\tag 3.1.1
$$
for any $v\in V$ of weight $\mu$ (cf. Section 1.1). 

Define a bifunctor $\bo: \V_\h\times \V_\h\to \V_\h$ as follows. 
For any $X,Y\in \V_\h$, define $X\bo Y$ to be the usual tensor product 
$X\o Y$. For any two morphisms $f:X\to X'$, $g:Y\to Y'$ define
the morphism $f\bo g: X\o Y\to X'\o Y'$ by the formula
$$
f\bo g(\la)=f^{(1)}(\la -\gamma h^{(2)})(1\o g(\la)).\tag 3.1.2
$$

It is easy to see that the category $\V_\h$ equipped with 
the bifunctor $\bo$ is a tensor category (cf. \cite{Mac}).  
Indeed, the functors $*\bo (*\bo *)$ and 
$(*\bo *)\bo *$ are equal, so $\bo$ is associative. 
Moreover, the object $\bold 1=\C$
(the trivial $\h$-module), satisfies the condition 
$\bold 1=\bold 1\bo \bold 1$,
and the functors $X\to \bold 1\bo X$, $X\to X\bo \bold 1$ are autoequivalences
of $\V_\h$, so $\bold 1$ is an identity object in $\V_\h$.

We will call this monoidal category the category of 
$\h$-vector spaces. 
If $\h=0$, the category $\V_\h$ coincides with the category
of complex vector spaces. 

If $\gamma=0$, the category $\Vh$ is equivalent, as a tensor category,
 to the category of diagonalizable $\h$-modules, 
with scalars extended from $\C$ to $M_{\h^*}$.
This case is not very interesting, so from now on we will assume that 
$\gamma\ne 0$.  

The category $\Cal V_\h$ depends on $\gamma$, but the categories 
with different nonzero $\gamma$ are obviously equivalent. 
We will suppress the dependence of $\Cal V_\h$ on $\gamma$ in the notation. 
 
{\bf Remark.} 
It is clear that for any two objects $X,Y\in\Vh$ the permutation operator
$\sigma_{XY}: X\bo Y\to Y\bo X$ is an isomorphism in $\Vh$. However, if 
$\h\ne 0$, then
this isomorphism is not functorial in $X$ and $Y$. In fact, it is quite
easy to see that if $\h\ne 0$, 
there is no functorial isomorphism between $X\bo Y$ and
$Y\bo X$: such an isomorphism would have to conjugate 
$f^{(1)}(\la-\gamma h^{(2)})(1\o g(\la))$ into 
$g^{(1)}(\la-\gamma h^{(2)})(1\o f(\la))$ 
for any $f,g$, which is impossible, since there is no relation between 
$f(\la)$ and $f(\la-\gamma\mu)$ 
for a generic function $f$.  Thus, the category $\Vh$ is a tensor category 
which in general does not admit a braided structure.

\subhead 3.2. Dynamical quantum R-matrices and tensor functors\endsubhead
 
It is known from the theory of quantum groups that 
if we are given a braided monoidal category $\B$,
a symmetric tensor category $\V$, and
a tensor functor $F:\B\to \V$, then for any object $X\in \B$
we can construct an element $R(\Cal B,F,X)\in \text{Aut}_\V(F(X)\o F(X))$ 
which satisfies the quantum Yang-Baxter equation, by the formula
$$
R(\Cal B,F,X)=F(\beta_{XX})P,\tag 3.2.1 
$$
 where 
$$
\beta_{XY}: X\o Y\to Y\o X 
$$ is
the braiding in $\B$, and $P$ is the permutation. 
For brevity we will write $R(\Cal B,F,X)$ as $R_X$. 

Suppose now that we are given a braided monoidal category $\B$ and
a tensor functor $F:\B\to \V_\h$. 
Observe that formula (3.2.1)
makes sense in this situation. However, since 
$\sigma_{XY}$ is not a functorial isomorphism, we should not expect 
$R_X$ to be a solution to the quantum 
Yang-Baxter equation. Still, it turns out that $R_X$ 
satisfies a modified version of the quantum Yang-Baxter equation, namely, 
the quantum dynamical Yang-Baxter equation (1.1.3). 

\proclaim{Theorem 3.1} The element $R_X$ 
satisfies the quantum dynamical Yang-Baxter equation (1.1.3)
in $\End_\Vh(F(X)^{\bo 3})$.
\endproclaim

\demo{Proof} We start with the braid relation
$$
(\beta\o 1)(1\o\beta)(\beta\o 1)=(1\o\beta)(\beta\o 1)(1\o\beta).\tag 3.2.2
$$
Applying the functor $F$ to (3.2.2), and using the definition of
the tensor product of morphisms in $\Vh$, we get (1.1.3).
$\square$\enddemo
 
\vskip .05in

\subhead 3.3. Representations of a quantum dynamical R-matrix\endsubhead

The notions discussed in this section were introduced in \cite{F1,F2,FV1}.

Let $R:\h^*\to \End(V\o V)$ be a quantum dynamical R-matrix
(see Chapter 1).

\proclaim{Definition} A representation of $R$ 
is an object $W\in \Vh$ endowed with
 an invertible morphism $L\in \End_\Vh(V\bo W)$, called the L-operator,  
such that 
$$
\align
R^{12}(\la - \gm h^{(3)})\,
L^{13}(\la )\, &
L^{23}(\la - \gm h^{(1)})\,
\tag 3.3.1
\\
=
L^{23}(\la )\,&
L^{13}(\la - \gm h^{(2)})\,
R^{12}(\la)\,.
\endalign
$$
in $\End_\Vh(V\bo V\bo W)$.
\endproclaim

{\bf Examples.} 1. The trivial representation: 
$W=\C$,  
$L=\text{Id}$. 

2. The basic representation: $W=V$, $L=R$.

Let $(W,L)$ be a representation of $R$. Let
$A\in \text{Aut}_\Vh(W)$. 
Let $L^A(\la):=(1\o A(\la)^{-1})L(\la)(1\o A(\la-\gm h^{(1)}))$.  

\proclaim{Lemma 3.1} $(W,L^A)$ is a representation of $R$.
\endproclaim

\demo{Proof} Straightforward.
$\square$\enddemo

Let $(W,L_W)$ and $(U,L_U)$ be representations of $R$. 
 
\proclaim{Definition} A morphism $A\in\text{Hom}_\Vh(W,U)$
is called an $R$-morphism if 
$$
(1\o A(\la))L_W(\la)=L_U(\la)(1\o A(\la-\gm h^{(1)})),\tag 3.3.2
$$
\endproclaim

Denote the space of $R$-morphisms from $W$ to $U$ by $\text{Hom}_R(W,U)$. 

It is clear that the composition of two R-morphisms is again an R-morphism. 
Thus, representations of $R$ form a category, which we denote by $Rep(R)$. 
This category
is additive, with the obvious notion of direct sum. 

\proclaim{Definition} The tensor product of $W$ and $U$ is the pair
$(W\o U, L_{W\o U})$, where
$$
L_{W\o U}(\la)=L^{12}_W(\la-\gm h^{(3)})L^{13}_U(\la).\tag 3.3.3
$$
\endproclaim

\proclaim{Lemma 3.2} $(W\o U, L_{W\o U})$ is a representation of $R$. 
\endproclaim
 
\demo{Proof} Straightforward.
$\square$\enddemo

It is clear that $(W\o U)\o X=W\o (U\o X)$.

\proclaim{Lemma 3.3} If $W,W',U,U'$ are representations of $R$ and
$f,g$ are $R$-morphisms then $f\bo g$ is an $R$-morphism.
\endproclaim

\demo{Proof} Straightforward.
$\square$\enddemo

Thus, we have equipped the category $Rep(R)$  
with a structure of a tensor category. 
Moreover, the forgetful functor $F: Rep(R)\to\Vh$ is naturally a tensor
functor. 

Theorem 3.1 shows that any pair $(\B,\ F:\B\to\Vh)$ defines
a system of quantum dynamical R-matrices. It turns out that conversely,
any quantum dynamical R-matrix $R$ defines $\B$, $F$, and $X$, such 
that $R=R(\B,F,X)$. The construction of
$\B,F,X$ is parallel to the case of usual R-matrices ($\h=0$),
where it is well known.

Namely, let $\B$ be the subcategory of $Rep(R)$ whose objects are 
tensor powers of $V$, and morphisms are the same as in 
$\text{Rep} (R)$. It is clearly a monoidal category.
Define a braiding $\beta$ on $\B$ by $\beta_{VV}=
RP$ (from the dynamical Yang-Baxter 
equation it follows that this is a morphism of representations
$V\o V\to V\o V$). 
It is easy to check using the hexagon axioms for the braiding 
that there exists  a unique braiding on 
$\B$ with such $\beta_{VV}$. 

Let $F:\B\to \Vh$ be the forgetful functor. We assign the pair $(\B,F)$ 
to $R$. It is clear that $R=R(\B,F,X)$ if we take $X=V$.

\subhead 3.4. Dual representations\endsubhead

It is useful to define the notion of the 
left and right dual representations. 

\proclaim{Definition} Let $(W,L_W)$ be a representation of $R$. 
The right dual representation to $W$ is the pair $(W^*,L_{W^*})$, where 
$W^*$ denotes the $\h$-graded dual of $W$, and  
$$
L_{W^*}(\la)=L_W^{-1}(\la+\gamma h^{(2)})^{t_2},\tag 3.4.1
$$
provided that the r.h.s. of (3.4.1) is invertible
(here $t_2$ denotes dualization in the second component).
The left dual representation to $W$ is the pair $({ }^*W,L_{{ }^*W})$, where 
${ }^*W=W^*$, and 
$$
L_{^*W}(\la)=L_W^{t_2}(\la-\gamma h^{(2)})^{-1},\tag 3.4.2
$$
provided that the r.h.s. of (3.4.2) is well defined. 
\endproclaim 

{\bf Remark 1.} Here $L_W^{-1}(\la+\gm h^{(2)})^{t_2}$ denotes the result
of three operations applied successively to $L_W$: inversion,  
shifting of the argument, and dualization in the second component. Similarly,
 $L_W^{t_2}(\la-\gm h^{(2)})^{-1}$ denotes the result
of three operations applied successively to $L_W$: 
dualization in the second component, shifting of the argument, and inversion. 
   
{\bf Remark 2.} 
We do not define the representation 
$W^*$ if $L_{W^*}$ is not invertible, and do not define the representation 
${}^*W$ if $L_W^{t_2}$ is not invertible.

\proclaim{Lemma 3.4} The right dual 
representation $(W^*,L_{W^*})$ and the left dual representation 
$({ }^*W,L_{{ }^*W})$
are representations 
of $R$, and if $W$ has finite dimensional weight subspaces then ${ }^*(W^*)=
({ }^*W)^*=W$. 
\endproclaim

\demo{Proof} The Lemma can be checked by a direct calculation. 
It also follows from Propositions 4.1 and 4.4 below. 
$\square$\enddemo

\proclaim{Lemma 3.5} If $A:W_1\to W_2$ is a homomorphism 
of representations of $R$, then the linear map $A^*(\la):=
A(\la+\gm h^{(1)})^t=A^t(\la-\gm h^{(1)})$ is a
homomorphism of representations $W_2^*\to W_1^*$, and is a homomorphism of 
represenations ${ }^*W_2\to { }^*W_1$, when these representations are defined. 
\endproclaim 

\demo{Proof} The Lemma can be checked by a direct calculation. 
It also follows from Propositions 4.1 and 4.4.
$\square$\enddemo

{\bf Remark.} It is easy to show that 
for two finite dimensional representations $W_1,W_2$ of $R$, the representation
$(W_1\o W_2)^*$ is naturally isomorphic to $W_2^*\o W_1^*$, 
and similarly for the left dual, if the corresponding dual representations 
are defined. 

\head 4. $\h$-Hopf algebroids and their dynamical representations\endhead

In this Chapter we will define the notion of an
$\h$-bialgebroid, 
and give the simplest nontrivial examples --
dynamical quantum groups associated to
quantum dynamical R-matrices from Chapter 1. 
We will generalize this material in 
the next chapter. 

\subhead 4.1. $\h$-bialgebroids\endsubhead

Let $\h$ be a finite dimensional commutative Lie algebra over $\C$, and 
$\gm$ a nozero complex number. Recall that $M_{\h^*}$ denotes the field 
of meromorphic functions on $\h^*$. 

\proclaim{Definition} An $\h$-algebra  
with step $\gm$ is an associative algebra 
$A$ over $\C$ with $1$, endowed with an $\h^*$-bigrading 
$A=\oplus_{\al,\beta\in \h^*}A_{\al\beta}$
(called the weight decomposition), and two algebra embeddings
$\mu_l,\mu_r:M_{\h^*}\to A_{00}$ (the left and the right moment maps), such  
that for any $a\in A_{\al\beta}$ and $f\in M_{\h^*}$, we have 
$$
\mu_l(f(\la))a=a\mu_l(f(\la+\gm\al)),\quad
\mu_r(f(\la))a=a\mu_r(f(\la+\gm\beta)).\tag 4.1.1
$$
\endproclaim

A morphism $\phi:A\to B$ of two $\h$-algebras 
is an algebra homomorphism, preserving 
the moment maps. By (4.1.1), such a homomorphism 
also preserves the weight decomposition. 

Let $A,B$ be two $\h$-algebras   with step $\gamma$, 
and $\mu_l^A,\mu_r^A,\mu_l^B,\mu_r^B$ their moment maps. 
Define their ``matrix tensor product'', 
$A\wo B$, which is also an $\h$-algebra. 

\proclaim{Definition} Let
$$
(A\wo B)_{\al\delta}:=\oplus_{\beta}A_{\al\beta}
\o_{M_{\h^*}} B_{\beta\delta},\tag 4.1.2
$$
 where $\o_{M_{\h^*}}$ means 
the usual tensor product 
modulo the relation $\mu_r^A(f)a\o b=a\o \mu_l^B(f)b$, for any
$a\in A,b\in B, f\in M_{\h^*}$.  
\endproclaim

Introduce a multiplication in $A\wo B$ by the rule $(a\o b)(a'\o b')=
aa'\o bb'$. It is easy to show that 
this product is well defined (cf. Proposition 5.1). Define 
the moment maps for $A\wo B$ by 
$\mu_l^{A\wo B}(f)=\mu_l^A(f)\o 1$,  
$\mu_r^{A\wo B}(f)=1\o \mu_r^A(f)$. 
It is easy to check that this makes $A\wo B$ into an $\h$-algebra. 
It is clear that 
$\wo$ is functorial with respect to both factors, 
and $(A\wo B)\wo C=A\wo (B\wo C)$. However,  
$A\wo B$ is not, in general, isomorphic to $B\wo A$. 
 
{\bf Remark.} The name ``matrix tensor product'' is used because 
formula (4.1.2) reminds the matrix multiplication.

\proclaim{Definition} A coproduct on an $\h$-algebra $A$ 
is a homomorphism of $\h$-algebras 
$\Delta: A\to A\wo A$.
\endproclaim

Let $D_\h$ be the algebra of difference operators 
$M_{\h^*}\to M_{\h^*}$, i.e. 
operators of the form $\sum_{i=1}^nf_i(\la)T_{\beta_i}$, where
$f_i\in M_{\h^*}$, and 
for $\beta\in \h^*$ we denote by $T_\beta$ the field automorphism 
of $M_{\h^*}$ given by $(T_\beta f)(\la)=f(\la+\gm\beta)$. 

The algebra $D_\h$ is the simplest
nontrivial example of an $\h$-algebra. Indeed
if we define the weight decomposition by $D_\h=\oplus (D_\h)_{\al\beta}$, 
where $(D_\h)_{\al\beta}=0$ if $\al\ne\beta$, and 
$(D_\h)_{\al\al}=\{f(\la)T_\al^{-1}: f\in M_{\h^*}\}$, and the moment maps
$\mu_l=\mu_r:M_{\h^*}\to (D_\h)_{00}$ 
to be the tautological isomorphism, 
then $D_\h$ becomes an $\h$-algebra.  

\proclaim{Lemma 4.1} For any $\h$-algebra $A$, the algebras
$A\wo D_\h$ and $D_\h \wo A$ are canonically isomorphic to $A$. 
\endproclaim

\demo{Proof} Straightforward.
$\square$\enddemo

Lemma 4.1 shows that the category of $\h$-algebras 
equipped with the product $\wo$ is a 
monoidal category, where the unit object is $D_\h$. 

\proclaim{Definition} A counit on an $\h$-algebra $A$ is
a homomorphism of $\h$-algebras 
$\eps: A\to  D_\h$. 
\endproclaim

\proclaim{Definition} An $\h$-bialgebroid is an 
$\h$-algebra $A$ equipped with a coassociative coproduct 
$\Delta$ (i.e. such that $(\Delta\o \Id_A)\circ \Delta=
(\Id_A\o \Delta)\circ \Delta$), and a counit
$\eps$ such that $(\eps\o \Id_A)\circ \Delta=
(\Id_A\o \eps)\circ \Delta=\Id_A$. 
\endproclaim

The property of the counit in the definition makes sense because 
of Lemma 4.1. 

\subhead 4.2. Dynamical representations of $\h$-bialgebroids
\endsubhead

Let $W$ be a diagonalizable $\h$-module, and let $D^\al_{\h, W}\subset
\text{Hom}_\C(W,W\o D_\h)$ be the space of all difference 
operators on $\h^*$ 
with coefficients in $\End_\C(W)$, which have weight $\al$ with 
respect to the action of $\h$ in $W$. 

Consider the algebra $D_{\h, W}=\oplus_\al D_{\h, W}^\al$.  
This algebra has a weight decomposition 
$D_{\h,W}=\oplus_{\al,\beta} (D_{\h,W})_{\al\beta}$ defined as follows: 
if $g\in \text{Hom}_\C(W,W\o M_{\h^*})$ is
an operator of weight $\beta-\al$ 
then $gT_{\beta}^{-1}\in (D_{\h,W})_{\al\beta}$.

Define the moment maps $\mu_l,\mu_r: M_{\h^*}\to (D_{\h, W})_{00}$
by the formulas $\mu_r(f(\la))=f(\la)$,
$\mu_l(f(\la))=f(\la-\gm h)$.  

\proclaim{Lemma 4.2} The algebra $D_{\h, W}$ equipped with this weight 
decomposition and these moment maps is an $\h$-algebra. 
\endproclaim

\demo{Proof} Straightforward.$\square$\enddemo

\proclaim{Lemma 4.3} There is a natural embedding of $\h$-algebras
$\theta_{WU}:D_{\h, W}\wo D_{\h, U}\to D_{\h, W\o U}$, 
given by the formula $fT_\beta\o gT_\delta\to (f\bo g)T_\delta$,
where $\bo$ is defined in Chapter 3, and 
$f\in \text{Hom}(W,W\o M_{\h^*})$.  This embedding is 
 an isomorphism if $W,U$ are finite-dimensional. 
\endproclaim

\demo{Proof} We have to show that the map $\theta_{WU}$ is well 
defined, and is an embedding.
We also have to show that $\theta_{WU}$ is a homomorphism of $\h$-algebras,
which is an isomorphism in the finite-dimensional case. 

The fact that $\theta_{WU}$ is well defined follows from the 
identity $\phi(\la)f\bo g=f\bo \phi(\la-\gm h)g$, for any  
function $\phi\in M_{\h^*}$. The injectivity of $\theta_{WU}$, and 
its surjectivity in the finite dimensional case are straightforward. 

It remains to show that $\theta_{WU}$ is a homomorphism 
of $\h$-algebras. It is obvious that $\theta_{WU}$ preserves the moment 
maps, so it remains to show that it is multiplicative. 
We have
$$
\gather
\theta_{WU}((f(\la)T_\beta^{-1}\o g(\la)T_\delta^{-1})
(f'(\la)T_{\beta'}^{-1}\o g'(\la)T_{\delta'}^{-1}))=\\
\theta_{WU}(f(\la)f'(\la-\gm\beta)T_{\beta+\beta'}^{-1}\o 
g(\la)g'(\la-\gm\delta)T_{\delta+\delta'}^{-1})=\\
f^{(1)}(\la-\gm h^{(2)})f^{'(1)}(\la-\gm h^{(2)}-\gm\beta)
(1\o g(\la)g'(\la-\gm\delta))T_{\delta+\delta'}^{-1}=\\
f^{(1)}(\la-\gm h^{(2)})(1\o g(\la))f^{'(1)}(\la-\gm h^{(2)}-\gm\delta)
(1\o g'(\la-\gm\delta))T_{\delta+\delta'}^{-1}=\\
f^{(1)}(\la-\gm h^{(2)})(1\o g(\la))T_\delta^{-1}
f^{'(1)}(\la-\gm h^{(2)})
(1\o g'(\la))T_{\delta'}^{-1}=\\
\theta_{WU}(f(\la)T_\beta^{-1}\o g(\la)T_\delta^{-1})\theta_{WU}
(f'(\la)T_{\beta'}^{-1}\o g'(\la)T_{\delta'}^{-1}).\tag 4.2.1
\endgather
$$
The Lemma is proved.
$\square$\enddemo

\proclaim{Definition} A dynamical representation of an 
$\h$-algebra $A$ is a diagonalizable $\h$-module 
$W$ endowed with a homomorphism of $\h$-algebras 
$\pi_W:A\to D_{\h, W}$. 
A homomorphism of dynamical representations 
$\phi: W_1\to W_2$ is an element of $\text{Hom}_{\C}(W_1,W_2\o M_{\h^*})$
such that $\phi\circ \pi_{W_1}(x)=
\pi_{W_2}(x)\circ \phi$ for all $x\in A$.  
\endproclaim

{\bf Example.} If $A$ has a counit, then 
it has the trivial representation: $W=\C$, $\pi=\eps$. 

Suppose now that $A$ is an $\h$-bialgebroid. Then, if 
$W$ and $U$ are two dynamical representations of $A$, the $\h$-module 
$W\o U$ also has a natural structure of a dynamical representation, 
defined by $\pi_{W\o U}(x)=\theta_{WU}\circ (\pi_W\o \pi_U)\circ \Delta(x)$. 

It is easy to show that if $f:W_1\to W_2$ and $g:U_1\to U_2$ 
are homomorphisms of dynamical representations, then 
$f\bo g$ is a homomorphism $W_1\o U_1\to W_2\o U_2$
(where $\bo$ is defined in Chapter 3). 
This gives a rule of tensoring morphisms. 
Thus, dynamical representations of $A$ form a monoidal category
$\text{Rep}(A)$, 
whose identity object is the trivial representation. 

Moreover, 
 the category $\text{Rep}(A)$ is equipped with a natural 
tensor functor $\text{Rep} (A)\to \Vh$ to the category of 
$\h$-vector spaces -- the forgetful functor.

\subhead 4.3. $\h$-Hopf algebroids and dual representations\endsubhead

Let us introduce the notion of an antipode
on an $\h$-bialgebroid.

Let $A$ be an $\h$-algebra. A linear map $S:A\to A$ is called an
antiautomorphism of $\h$-algebras if it is an antiautomorphism 
of algebras and $\mu_r\circ S=\mu_l,\mu_l\circ S=\mu_r$. 
 From these conditions it follows that $S(A_{\al\beta})=A_{-\beta,-\alpha}$. 

Let $A$ be an $\h$-bialgebroid, and let $\Delta$, $\eps$
 be the coproduct and counit of $A$. For $a\in A$, let 
$$
\Delta(a)=\sum_i a^1_i\o a^2_i.\tag 4.3.1
$$ 

\proclaim{Definition} An antipode on the $\h$-bialgebroid $A$ 
is an antiautomorphism of $\h$-algebras $S:A\to A$ such that 
for any $a\in A$ and any presentation (4.3.1) of $\Delta(a)$, one has
$$
\sum_i a_i^1S(a_i^2)=\mu_l(\eps(a)1),\ 
\sum_i S(a_i^1)a_i^2=\mu_r(\eps(a)1),\tag 4.3.2
$$ 
where $\eps(a)1\in M_{\h^*}$ is the result of application of the difference 
operator $\eps(a)$ to the constant function $1$.  
\endproclaim

{\bf Remark.} It is easy to see that  
$\sum_i a_i^1S(a_i^2)$ and $\sum_i S(a_i^1)a_i^2$ depends only on $a$ and 
not on the choice of the presentation (4.3.1). 

\proclaim{Definition} An $\h$-bialgebroid with an antipode is called an
$\h$-Hopf algebroid. 
\endproclaim

{\bf Remark.} If $\h=0$, the notions of an $\h$-algebra, $\h$-bialgebroid, 
$\h$-Hopf algebroid coincide with the notions of an algebra, 
bialgebra, and Hopf algebra, respectively. 

For any $\h$-Hopf algebroid $A$, the category $\text{Rep}(A)$ has 
the following natural notion of the left and right dual representation. 

If $(W, \pi_W)$ is a dynamical representation of an $\h$-algebra $A$, 
we denote by $\pi_W^0:A\to \text{Hom}(W,W\o M_{\h^*})$ 
the map defined by $\pi_W^0(x)w=
\pi_W(x)w$, $w\in W$ 
 (the difference operator $\pi_W(x)$ restricted to the constant functions).
It is clear that $\pi_W$ is completely determined by $\pi_W^0$. 

\proclaim{Definition} Let $(W,\pi_W)$ be a dynamical representation of $A$. 
Then the right dual representation to $W$ is $(W^*,\pi_{W^*})$, 
where $W^*$ is the $\h$-graded dual to $W$, and 
$$
\pi_{W^*}^0(x)(\la)=\pi_W^0(S(x))(\la+\gm h-\gm \alpha)^t\tag 4.3.3
$$
for $x\in A_{\al\beta}$, where $t$ denotes dualization. 
The left dual representation to $W$ is $({ }^*W,\pi_{{ }^*W})$, 
where ${ }^*W=W^*$, and 
$$
\pi_{{ }^*W}^0(x)(\la)=
\pi_W^0(S^{-1}(x))(\la+\gm h-\gm \alpha)^t\tag 4.3.4
$$
for $x\in A_{\al\beta}$. 
\endproclaim

\proclaim{Proposition 4.1} Formulas (4.3.3) and (4.3.4) define 
dynamical representations of $A$. 
Moreover, if $A(\la):W_1\to W_2$ is a morphism 
of dynamical representations, then $A^*(\la):=A(\la+\gm h)^t$ 
defines a morphism $W_2^*\to W_1^*$ and ${ }^*W_2\to { }^*W_1$. 
\endproclaim

\demo{Proof} Let $x\in A_{\al_x\beta_x}$, $y\in A_{\al_y\beta_y}$.
Then $\pi_W^0(xy)(\la)=\pi_W^0(x)(\la)\pi_W^0(y)(\la-\gm \beta_x)$
by the definition of a dynamical representation. Therefore, we have 
$$
\gather
\pi_{W^*}^0(xy)(\la)^t=\pi_W^0(S(xy))(\la+\gm h-\gm \al_x-\gm \al_y)
=\pi_W^0(S(y)S(x))(\la+\gm h-\gm \al_x-\gm \al_y)=\\
\pi_W^0(S(y))(\la+\gm h-\gm \al_x-\gm \al_y+\gm \al_{S(x)}-\gm
\beta_{S(x)})
\pi_W^0(S(x))(\la+\gm h-\gm \al_x-\gm \al_y-\beta_{S(y)})=\\ 
\pi_W^0(S(y))(\la+\gm h-\gm \al_y-\gm \beta_x)
\pi_W^0(S(x))(\la+\gm h-\gm \al_x).\tag 4.3.5\endgather
$$
Dualizing (4.3.5), we get
$$
\gather
\pi_{W^*}^0(xy)(\la)=\pi_W^0(S(x))(\la+\gm h-\gm \al_x)^t
\pi_W^0(S(y))(\la+\gm h-\gm \al_y-\gm \beta_x)^t=\\
\pi_{W^*}^0(x)(\la)\pi_{W^*}^0(y)(\la-\gm \beta_x),\tag 4.3.6
\endgather
$$
which implies the first statement of the Proposition for $W^*$.
The proof for ${}^*W$ is obtained by replacing $S$ by $S^{-1}$. 

Let us prove the second statement. The intertwining property of $A(\la)$ 
can be written as
$$
A(\la)\pi_W^0(x)(\la)=\pi_W^0(x)(\la)A(\la-\gm \beta_x).\tag 4.3.7
$$
Replacing $x$ with $S(x)$ and shifting the arguments, we get
$$
\gather
A(\la+\gm h-\gm \beta_x)
\pi_W^0(S(x))(\la+\gm h-\gm \al_x)=\\
\pi_W^0(S(x))(\la+\gm h-\gm \al_x)A(\la+\gm h-\gm\al_x-
\gm\beta_S(x)).\tag 4.3.8\endgather
$$
Dualizing (4.3.8) and using the identity $\beta_{S(x)}+\al_x=0$, we get 
the second statement of the Proposition. 
The Proposition is proved.
$\square$\enddemo

\subhead 4.4. $\h$-bialgebroids associated to a 
 function $R:\h^*\to \End(V\o V)$\endsubhead

Let $\h$ be a finite dimensional commutative Lie algebra, 
and $V=\oplus_{\al\in \h^*} V_\al$ 
a finite dimensional diagonalizable $\h$-module. 
Let $R(\la)$ be a meromorphic function $\h^*\to \End(V\o V)$, 
such that $R(\la)$ is invertible for a generic $\la$.   
Using $R$, we will now define an $\h$-bialgebroid $A_R$ 
which we call the {\it dynamical quantum group} corresponding to $R$. 
This construction is analogous to the Faddeev-Reshetikhin-Sklyanin-Takhtajan
construction of the quantum function algebra on $GL_N$. 

As an algebra, $A_R$ by definition is generated by two copies of  
$M_{\h^*}$ (embedded as subalgebras) and certain new generators, 
which are matrix elements of the operators
$L^{\pm 1}\in \text{End}(V)\o A_R$.
 We denote the elements of the first copy of 
$M_{\h^*}$ as $f(\la^1)$ and of the second copy as $f(\la^2)$, 
where $f\in M_{\h^*}$. We denote by $(L^{\pm 1})_{\al\beta}$ 
the weight components of $L^{\pm 1}$ with respect to the natural 
$\h$-bigrading on $\End (V)$, so 
that $L^{\pm 1}=(L^{\pm 1}_{\al\beta})$, where
$L^{\pm 1}_{\al\beta}\in \text{Hom}_\C(V_\beta,V_\al)\o A_R$.    

Then the defining relations for $A_R$ are:
 
$$
f(\la^1)L_{\al\beta}=L_{\al\beta}f(\la^1+\gm\al);\
f(\la^2)L_{\al\beta}=L_{\al\beta}f(\la^2+\gm\beta);
[f(\la^1),g(\la^2)]=0;\tag 4.4.1
$$
$$
LL^{-1}=L^{-1}L=1;\tag 4.4.2
$$
and the dynamical Yang-Baxter relation 
$$
R^{12}(\la^1)L^{13}L^{23}=:L^{23}L^{13}R^{12}(\la^2):,\tag 4.4.3
$$
Here the :: sign (``normal ordering'')
means that the matrix elements of $L$ 
should be put on the right of the matrix elements of $R$. Thus, if 
$\{v_a\}$ is a homogeneous basis of $V$, and $L=\sum E_{ab}\o L_{ab}$, 
$R(\la)(v_a\o v_b)=\sum R^{ab}_{cd}(\la)v_c\o v_d$, then (4.4.3) has the form
$$
\sum R^{xy}_{ac}(\la^1)L_{xb}L_{yd}=\sum R_{xy}^{bd}(\la^2)
L_{cy}L_{ax},\tag 4.4.4
$$ 
where we sum over repeated indices. 

More precisely, the algebra $A_R$ is, by definition, the quotient of the 
algebra $\tilde A$ freely generated by $M_{\h^*}\o M_{\h^*}$ and 
elements $L_{ab},(L^{-1})_{ab}$, $a,b=1,...,\text{dim}V$, 
by the ideal defined by relations (4.4.1)-(4.4.3). 

Introduce the moment maps for $A_R$ by 
$\mu_l(f)=f(\la^1)$, $\mu_r(f)=f(\la^2)$, and the weight 
decomposition by $f(\la^1),f(\la^2)\in (A_R)_{00}$, $L_{\al\beta}\in 
\text{Hom}_\C(V_\beta,V_\al)\o (A_R)_{\al\beta}$.
It is clear that $A_R$ equipped with such 
structures is an $\h$-algebra. 

Now define the coproduct on $A_R$, $\Delta: A_R\to A_R\wo A_R$, by
the usual Lie-theoretic formulas
$$
\Delta(L)=L^{12}L^{13},\Delta(L^{-1})=(L^{-1})^{13}(L^{-1})^{12} \tag 4.4.5
$$
(here $\Delta$ is applied to the second component of $L,L^{-1}$). 

\proclaim{Proposition 4.2} $\Delta$ 
extends to a well defined homomorphism 
$A\to A\wo A$. 
\endproclaim

\demo{Proof} From (4.4.5) we get
$$
\Delta(L_{\al\beta})=\sum_{\theta}L_{\al\theta}^{12} 
L_{\theta\beta}^{13}.\tag 4.4.6
$$
So it remains to show that the defining relations 
of $A_R$ are invariant 
under $\Delta$. The invariance of relations (4.4.1) follows directly 
from (4.4.6). Relation (4.4.2) is obviously invariant.
To check the invariance of relation (4.4.3), we have to 
show that 
$$
R^{12}(\la_1^1)L^{13}L^{14}L^{23}L^{24}=:L^{23}L^{24}L^{13}L^{14}R^{12}
(\la_2^2): .\tag 4.4.7
$$
(the subscripts $1,2$ under $\la$ indicate that the corresponding 
functions are taken from the first and the second components of $A_R$ in the 
product $A_R\wo A_R$; and, as before, the :: sign indicates 
that the functions of $\la^i$ are written 
on the left from the L-operators).  

We have
$$
\gather
R^{12}(\la_1^1)L^{13}L^{14}L^{23}L^{24}=R^{12}(\la_1^1)L^{13}L^{23}L^{14}L^{24}
=:L^{13}L^{23}R^{12}(\la_1^2):L^{14}L^{24}=\\ 
L^{13}L^{23}R^{12}(\la_2^1)L^{14}L^{24}=L^{23}L^{13}:L^{24}L^{14}R^{12}
(\la_2^2):=:L^{23}L^{24}L^{13}L^{14}R^{12}
(\la_2^2): .\tag 4.4.8\endgather
$$
(We replaced $\la_1^2$ by $\la_2^1$ in the middle of (4.4.8) 
since $A_R\wo A_R$ is by definition inside of the tensor product
$A_R\o_{M_{\h^*}}A_R$, where $M_{\h_*}$ is mapped into the first 
component of $A_R$ by $\mu_r$ and into the second by $\mu_l$, acting 
from the left). The Proposition is proved. 
$\square$\enddemo

Now define the counit on the algebra $A_R$. Recall that 
the counit
has to be an algebra homomorphism $\epsilon: A_R\to D_\h$. 

Define the counit by the formula
$$
\epsilon(L_{\al\beta})=
\delta_{\al\beta} \Id_{V_\al} \o T_\al^{-1},
\epsilon((L^{-1})_{\al\beta})=
\delta_{\al\beta} \Id_{V_\al} \o T_\al,
\tag 4.4.9
$$
where $\Id_{V_\al}: V_\al\to V_\al$ is the identity operator. 
 
We need to check that the counit is well defined, i.e. that the defining 
relations are annihilated by it. For relations (4.4.1),(4.4.2) it is 
obvious. 
Relation (4.4.3) reduces to checking that
$$
(\sum_{\al,\beta}R^{12}(\la) (\Id_{V_\al}\o \Id_{V_\beta}))\o T_{\alpha+\beta}^{-1}=
(\sum_{\al,\beta} (\Id_{V_\al}\o \Id_{V_\beta})R^{12}(\la))\o T_{\alpha+\beta}^{-1},
\tag 4.4.10
$$
which holds because $R$ has zero weight. 

\proclaim{Proposition 4.3}
The counit axiom
$(\Id\o \epsilon)\circ \Delta=(\epsilon \o \Id)\circ \Delta
=\Id$ is satisfied for $A_R$.
\endproclaim

\demo{Proof} We need to check the relations on $L$. 
These relations follow from the fact
that the elements $T_\al^{-1}\o L_{\al\beta}$,
$L_{\al\beta}\o T_\beta^{-1}$ are mapped to $L_{\al\beta}$ 
under the natural isomorphisms $D_\h\wo A_R\to A_R$, 
$A_R\wo D_\h\to A_R$.  
$\square$\enddemo

Thus, $A_R$ is an $H$-biequivariant 
bialgebroid.  
We will call it the {\it dynamical quantum group} corresponding to 
the function $R$. 

It is also possible to consider the algebra generated by 
$f(\la_1),f(\la_2)$, $L$ (without $L^{-1}$). Denote this algebra 
by $\bar A_R$. The algebra $\bar A_R$ is an $\h$-bialgebroid, 
which is naturally mapped to $A_R$. 

{\bf Remark.} The algebra $\bar A_R$ was introduced in \cite{FV1}
under the name of ``the operator algebra''. 

\subhead 4.5. The antipode on $A_R$\endsubhead

Let $A,B$ be algebras with $1$. For $X\in B\o A$, 
define $i(X)$ to be the inverse of $X$, and $i_*(X)$
to be the inverse of $X$ in the algebra $B\o A_{op}$, 
where $A_{op}$ is $A$ with the reversed order of multiplication. 
Clearly, $i^2=i_*^2=\Id$. 

Let $I$ be the group freely generated by $i,i^*$ with relations 
$i^2=i_*^2=\Id$. We will say that an element $X$ is 
{\it strongly invertible} if for any $g\in I$ 
the element $g(X)$ is well defined. 

\proclaim{Definition} An invertible, 
weight zero matrix function $R$ is said to be rigid if
the element $L\in \End(V)\o A_R$ is strongly invertible.
\endproclaim

\proclaim{Proposition 4.4} $R$ is rigid if and only if 
$A_R$ admits an antipode $S$ such that $S(L)=L^{-1}$.
In this case, $S^{2n}(L)=(i^*i)^n(L),S^{2n+1}(L)=i(i^*i)^n(L)$. 
In particular, $S(L^{-1})=i^*i(L)$.  
\endproclaim

\demo{Proof} Suppose that $R$ is rigid. Extend the definition of the 
antipode by $S(L^{-1})=i_*(L^{-1})=i_*i(L)$. 
It is easy to see that the relations of $A_R$ are preserved, so 
this indeed defines an antihomomorphism $S:A\to A$. Moreover, 
$S$ is an isomorphism: the inverse is given by $S^{-1}(L^{-1})=L$, 
$S^{-1}(L)=i_*(L)$. 

Now suppose that $S$ is defined. Then it is easy to check that 
$(i^*i)^n(L)=S^{2n}(L)$, $i(i^*i)^n(L)=S^{2n+1}(L)$, $n\in \Bbb Z$. 
This defines $g(L)$ for all $g\in I$. 
The proposition is proved. 
$\square$\enddemo

{\bf Remark 1.} 
The proposition shows that for rigidity of $R$, it is sufficient
that $i_*(L)$ and $i_*(L^{-1})$ be defined. 
 
{\bf Remark 2.} Observe that in general $S^2\ne 1$. 

Thus, if $R$ is rigid then $A_R$ is an $\h$-Hopf algebroid. 

\subhead 4.6. Representation theory of $A_R$\endsubhead

Now consider the representation theory of $A_R$.  
As was pointed out in \cite{FV1}, the category
$\text{Rep}(A_R)$ of dynamical representations of 
$A_R$ is tautologically isomorphic to the category $\text{Rep}(R)$ 
of representations of $R$.  

\proclaim{Proposition 4.4} The tensor categories 
$Rep(A_R)$ and $Rep(R)$ are equivalent.
\endproclaim

\demo{Proof}
Define the functor $\Gamma: Rep(A_R)\to Rep(R)$ to be the 
identity at the level of vector spaces, and set
$$
L_{\Gamma(W)}=\pi_W^0(L). \tag 4.6.1
$$
Define the functor $\Gamma^{-1}: Rep(R)\to Rep(A_R)$ by
$$
\pi_{\Gamma^{-1}(W)}^0(L)=L_W.\tag 4.6.2
$$
These functors preserve tensor structure, and are obviously inverse to each 
other. The Proposition is proved.
$\square$\enddemo

It is easy to see that 
the functor $\Gamma$ commutes with the duality functors.  
Therefore, if $R$ is rigid, then the representations $W^*$,${}^*W$ of $R$ 
are well defined for any $W$, and the category $Rep_f(R)$
of finite-dimensional representations of $R$ (= the category $Rep_f(A_R)$
of finite dimensional dynamical representations of $A_R$) 
is a rigid tensor category\cite{DM}. This explains our use of the word 
``rigid''.  

Although $A_R$ is an $\h$-Hopf algebroid 
for any rigid zero weight function $R$, it does not always 
have nice properties. For a generic $R$, this algebra will be very small 
and will not have interesting dynamical representations. 
However if $R$ is a dynamical quantum  
R-matrix, then the category $\text{Rep}(R)$ is nontrivial 
(it contains the basic representation defined in Chapter 3), so by 
Proposition 4.4 the category $\text{Rep}(A_R)$ 
 is also nontrivial. Thus, algebras $A_R$ with $R$ being a dynamical 
quantum R-matrix form a good class of $\h$-Hopf algebroids. 
 From now on we will only consider $A_R$ for $R$ being a dynamical quantum
R-matrix. 

\subhead 4.7. Sufficient conditions for rigidity\endsubhead

Unfortunately, the definition of rigidity cannot be effectively checked, 
since it depends on the properties of the algebra $A_R$, about whose 
structure we do not know very much. 
Therefore, we would like to find some effective sufficient 
conditions of rigidity. 

For any function $X:\h^*\to \End(V\o V)$, define the function 
$\tilde X: \h^*\to \End(V\o V)$ as follows. Suppose that 
for $v,w\in V$ one has $X(\la)(v\o w)=\sum_i f_i(\la)v_i\o w_i$, where 
$f_i\in M_{\h^*}$ and $w_i$ are homogeneous. Then set
$\tilde X(\la)(v\o w)=\sum_i f_i(\la\,+\,\gm\, wt(w_i))v_i\o w_i$, where
$wt(w_i)$ denotes the weight of $w_i$. 

Let $R$ be a dynamical quantum R-matrix with step $\gm$.
Assume that $i_*(\tilde R)$ is defined.

Let us write $\tilde R$ in the form $\tilde R=\sum a_i\o b_i$, 
and $i_*(\tilde R)$ in the form
$i_*(\tilde R)=\sum c_i\o d_i$.

Define the operators 
$Q=\sum d_ic_i, Q'=\sum c_id_i: \h^*\to\End(V)$.
These operators are of weight zero with respect to $\h$, since $R$ is of 
weight zero. 

\proclaim{Proposition 4.5} Suppose $R$ is such that
$i_*(\tilde R)$ is defined, and $R$ satisfies the following conditions:

(i) The operator $Q$ is invertible 
for a generic $\la$. 

(ii) The operator $Q'$ is invertible 
for a generic $\la$. 

Then $R$ is rigid, and 
$$
i_*(L^{-1})=S^2(L)=:(Q(\la^1)\o 1)L(Q^{-1}(\la^2)\o 1):=
:(Q'(\la^1+\gm h)^{-1}\o 1)L(Q'(\la^2+\gm h)\o 1):.\tag 4.7.1
$$ 
\endproclaim

{\bf Remark.} It is clear that (i) and (ii) are satisfied for $R=1$ and
are open conditions. Therefore, 
Proposition 4.5 shows that if $R_\gm$ is a continuous family of 
quantum dynamical R-matrices with step $\gm$ such that $R_0=1$, then 
$R_\gm$ is rigid for small $\gm$. 

\demo{Proof} First of all, let us deduce a commutation relation 
between $L$ and $L^{-1}$. 

Multiplying the dynamical Yang-Baxter equation 
by $(L^{-1})^{23}$ on the right, we get
$$
R^{12}(\la^1)L^{13}=:L^{23}L^{13}R^{12}(\la^2)(L^{23})^{-1}:,\tag 4.7.2
$$
Let $\{v_a\}$ be an $\h$-homogeneous basis of $V$, and 
$L=\sum E_{ab}\o L_{ab}$. Denote by $\omega_a$ the weight of $v_a$. 
Then we have 
$$
\gather
:L^{23}L^{13}R^{12}(\la^2):(L^{23})^{-1}=
\sum E_{ab}^{(2)}:L_{ab}^{(3)}L^{13}R^{12}(\la^2)(L^{23})^{-1}:=\\
\sum E_{ab}^{(2)}L_{ab}^{(3)}:L^{13}R^{12}((\la+\gm\omega_b)^2)(L^{23})^{-1}:=\\
\sum E_{ab}^{(2)}L_{ab}^{(3)}:L^{13}\tilde 
R^{12}(\la^2)(L^{23})^{-1}:=L^{23}:L^{13}
\tilde R^{12}(\la^2):(L^{23})^{-1}.\tag 4.7.3
\endgather
$$
Therefore, multiplying (4.7.2) on the left by $(L^{23})^{-1}$ we get
$$
(L^{23})^{-1}:
R^{12}(\la^1)L^{13}:=:L^{13}\tilde R^{12}(\la^2)(L^{23})^{-1}:,\tag 4.7.4
$$
Transforming the left hand side of this equation similarly to (4.7.3), 
we arrive at the equation
$$
:(L^{23})^{-1}\tilde
R^{12}(\la^1)L^{13}:=:L^{13}\tilde R^{12}(\la^2)(L^{23})^{-1}:,\tag 4.7.5
$$
which is the desired commutation relation. 

Now, using property (i), define 
$$
T=:(Q(\la^1)\o 1)L(Q^{-1}(\la^2)\o 1):\in \End(V)\o A_R.\tag 4.7.6
$$
Let * denote the product in the algebra 
$\End(V)\o (A_R)_{op}$. Let us compute the product $L^{-1}*T$. 

Set $L^{-1}=\sum E_{ab}\o (L^{-1})_{ab}$. Then we get
$$
L^{-1}*T=\sum (E_{pq}Q(\la^2)E_{rs}Q^{-1}(\la^1)\o 1)(1\o L_{rs}(L^{-1})_{pq}).
\tag 4.7.7
$$
Using (4.7.5), we can rewrite (4.7.7) in the form
$$
L^{-1}*T=(\sum d_i(\la^2)E_{rs}b_j(\la^1)Q(\la^1)a_j(\la^1)E_{pq}c_i(\la^2)
Q^{-1}(\la^2)\o 1)(1\o L_{rs}(L^{-1})_{pq}).\tag 4.7.8
$$
Using the definition of $Q$, we have
$$
\sum b_iQa_i=1.\tag 4.7.9
$$
Substituting (4.7.9) into (4.7.8), we get $L^{-1}*T=1$. 

Now, using property (ii), define 
$$
T'=
:(Q'(\la^1+\gm h)^{-1}\o 1)L(Q'(\la^2+\gm h)\o 1):.\tag 4.7.10
$$
Then, analogously to the above, we get
$T'*L^{-1}=1$. Thus, $T=T'=i_*(L^{-1})$. 

It is easy to see that 
$$
i_*(L)=:(Q^{-1}(\la^2)\o 1)L(Q(\la^1)\o 1):\tag 4.7.11
$$
Thus, $R$ is rigid. $\square$\enddemo

Now we will show that any rigid quantum dynamical R-matrix 
satisfies a certain crossing symmetry condition. 

For an invertible zero weight function $X(\la)\in \text{End}(V\o V)$, 
set 
$$
\tau(X)(\la)=X^{-1}(\la+\gm h^{(2)})^{t_2}.\tag 4.7.12
$$

\proclaim{Corollary 4.1}
Let $R$ be a rigid quantum dynamical R-matrix on $V$. 
Then $\tau(R)$ is invertible, and $R$
satisfies the crossing symmetry condition
$$
\tau^2(R)=(Q(\la-\gm h^{(2)})\o 1)R(\la)(Q^{-1}(\la)\o 1).\tag 4.7.13
$$
\endproclaim

\demo{Proof} It is clear that $\tau^2(R)=L_{V^{**}}$, where 
$V$ is the basic representation of $R$. Therefore, using (4.7.1) in the basic  
representation, we get (4.7.12). 
$\square$\enddemo

\subhead 4.8. Dynamical quantum groups associated to 
dynamical R-matrices of $gl_N$ type\endsubhead

Now suppose that $R$ is a dynamical R-matrix of $gl_N$-type. 
Then it has form (1.3.2), and we can write the defining relations
for $A_R$ more explicitly. Since all weight subspaces of $V$ are 
1-dimensional, we have $(L^{\pm 1})_{\al\beta}\in A$. For brevity we write 
$(L^{\pm 1})_{ab}$ for $(L^{\pm 1})_{\om_a\om_b}$. 
Thus, we have
$L^{\pm 1}=\sum E_{ab}\o (L^{\pm 1})_{ab}$.  

In this notation, the defining relations for $A_R$ look like
$$
\gather
LL^{-1}=L^{-1}L=1,\\
f(\la^1)L_{bc}=L_{bc}f(\la^1+\gamma\om_b),\
f(\la^2)L_{bc}=L_{bc}f(\la^2+\gamma\om_c), \\
L_{as}L_{at}=\frac{\alpha_{st}(\la^2)}{1-\beta_{ts}(\la^2)}L_{at}L_{as},
s\ne t,\\
L_{bs}L_{as}=\frac{\alpha_{ab}(\la^1)}{1-\beta_{ab}(\la^1)}L_{as}L_{bs},
a\ne b,\\
\al_{ab}(\la_1)L_{as}L_{bt}-\al_{st}(\la_2)L_{bt}L_{as}=
(\beta_{ts}(\la_2)-\beta_{ab}(\la_1))L_{bs}L_{at}, a\ne b, s\ne t, \tag 
4.8.1 \endgather
$$
where $\al_{ab}$, $\beta_{ab}$ are the functions from (1.3.2). 

{\bf Remark.} It is also possible to define dynamical quantum groups 
associated with dynamical R-matrices with spectral parameter. 
It is done analogously to the above, and we will do it in detail 
in a forthcoming paper. For example, if $R(z,\la)$ is 
a quantum dynamical R-matrix with spectral parameter of elliptic type 
(i.e. of the form (2.5.1)), we will get the elliptic 
quantum group defined
in \cite{F1,F2,FV1,FV2}. Relations (4.8.1) 
(for dynamical R-matrices of $gl_N$ Hecke type) can be obtained 
as a limiting case of the defining relations for the elliptic 
quantum group. 

\subhead 4.9. Rigidity of the 
rational and the trigonometric dynamical R-matrix\endsubhead
 
Consider the trigonometric dynamical R-matrix 
$R(\la)$ defined by (1.6.4), with $X=\{1,...,N\}$, and $\mu_{ab}=1$. 

\proclaim{Proposition 4.6} $R(\la)$ is rigid, and 
the matrices $Q,Q'$ are given by the formulas
$$
\gather
Q=\text{diag}(Q_1,...,Q_N),\ Q'=\text{diag}(Q_1',...,Q_N'),
\\
Q_a(\la)=
\prod_{i\ne a}
\frac{q^{1+\la_i}-q^{\la_a}}
{q^{\la_i}
-q^{\la_a}},\\
Q_a'(\la)=qQ_a^{-1}(\la), \tag 4.9.1\endgather
$$
where $q=e^\epe$.
\endproclaim

\demo{Proof}
First of all, it is not hard to show by a direct computation that the 
matrix $i_*(\tilde R)$ is defined. So it remains to show that the elements 
$Q,Q'$ are invertible.  

Let $P(\la)=Q'(\la+\gm h)$. 
Let $P_i,Q_i$ be the 
diagonal entries of $P,Q$. As we know, 
these entries are defined by the following systems of linear equations:
$$
\gather
Q_a+\sum_{b\ne a}\beta_{ab}(\la+\gm \om_a)Q_b=1,\\ 
P_a+\sum_{b\ne a}\beta_{ba}(\la+\gm \om_b)P_b=1.\tag 4.9.2  
\endgather
$$ 
 The explicit form of the systems (4.9.2) is 
$$
\gather 
Q_a+\sum_{b\ne a}\frac{q-1}{q^{1+\la_a-\la_b}-1}Q_b=1,\\
P_a+\sum_{b\ne a}\frac{q-1}{q^{1+\la_b-\la_a}-1}P_b=1.
\endgather
$$ 
Thus, if one of these systems is nondegenerate (which we show below) then 
$Q(\la)=P(-\la)$. 

 From now on we consider only the first system. 
Note that it can be conveniently written as
$$
\sum_{b}\frac{q-1}{q^{1+\la_a-\la_b}-1}Q_b=1\tag 4.9.3
$$
Define $X_b=q^{\la_b}Q_b$. Then (4.9.3) can be written as
$$
\sum_{b}\frac{1}{[1+\la_a]-[\la_b]}X_b=1,\tag 4.9.4
$$
where $[x]=\frac{q^x-1}{q-1}$. 
Thus, the vector $X$ is defined by $X=C^{-1}\bold 1$, 
where $C_{ab}=\frac{1}{[1+\la_a]-[\la_b]}$, and $\bold 1$ is the vector 
whose components are all equal to 1. 

To invert the matrix $C$, we use the well known combinatorial 
identity (which is 
called the ``Bose-Fermi correspondence'' in physics): 
$$
\text{det}(\frac{1}{x_i-y_j})=\frac{\prod_{i<j}(x_i-x_j)\prod_{i<j}(y_i-y_j)}
{\prod_{i,j}(x_i-y_j)}.\tag 4.9.5
$$
Applying this identity to $x_i=[1+\la_i],y_i=[\la_i]$, and using the
usual rule of inverting matrices, we get
$$
(C^{-1})_{ab}=\frac{\prod_{(i,j):i=b\text{ or }j=a} (x_i-y_j)}
{\prod_{j\ne b}(x_b-x_j)\prod_{i\ne a}(y_i
-y_a)}.
\tag 4.9.6
$$
In particular, 
$$
X_a=\sum_b (C^{-1})_{ab}=\frac{\prod_i(x_i-y_a)}
{\prod_{i\ne a}(y_i
-y_a)}
\sum_b\frac{\prod_{j\ne a}
(x_b-y_j)}{\prod_{j\ne b}(x_b-x_j)}.\tag 4.9.7
$$

{\bf Claim.} 
$$
\sum_b\frac{\prod_{j\ne a}
(x_b-y_j)}{\prod_{j\ne b}(x_b-x_j)}=1.\tag 4.9.8
$$

{\it Proof of the claim.} Consider the expression on the l.h.s. of
(4.9.8) as a rational function of $z=x_a$ for fixed $x_b$, $b\ne  a$. 
This function has no more than simple poles at 
$x_b$, $b\ne a$, and no other singularities; it equals 1 at 
infinity. Thus, it suffices to show that its residues vanish, 
which is obvious: only two terms contribute to each residue, ant these two 
terms cancel each other. 

Thus, we get:
$$
Q_a(\la)=q^{-\la_a}\frac{\prod_i([1+\la_i]-[\la_a])}
{\prod_{i\ne a}([\la_i]
-[\la_a])},\tag 4.9.9
$$
i.e. 
$$
\gather
Q_a(\la)=
\prod_{i\ne a}
\frac{q^{1+\la_i}-q^{\la_a}}
{q^{\la_i}
-q^{\la_a}},\\
P_a(\la)=
\prod_{i\ne a}
\frac{q^{1-\la_i}-q^{-\la_a}}
{q^{-\la_i}
-q^{-\la_a}}.\tag 4.9.10
\endgather
$$
Therefore, 
$$
Q_a'(\la)=P_a(\la-\omega_a)=
\prod_{i\ne a}
\frac{q^{1-\la_i}-q^{1-\la_a}}{q^{-\la_i}
-q^{1-\la_a}}
=\prod_{i\ne a}
\frac{q^{1+\la_i}-q^{1+\la_a}}{q^{1+\la_i}
-q^{\la_a}}=qQ_a^{-1}(\la). \tag 4.9.11
$$

Thus, $R$ is rigid, and $Q,Q'$ are given by formula (4.9.1). 
The proposition is proved. 
$\square$\enddemo 

An analogous theorem holds for the rational dynamical R-matrix
(1.5.1) (with $X=\{1,...,N\}$ and 
$\mu_{ab}=0$). The formulas for $Q,Q'$ for such $R$
are obtained from (4.9.1) as $q\to 1$. 

It is easy to show that the property of rigidity is preserved by gauge 
transformations, so we get

\proclaim{Corollary 4.2}
Any quantum dynamical R-matrix $R$ of $gl_N$ Hecke type is rigid. 
\endproclaim

Clearly, the elements $Q,Q'$ for any such $R$ can be easily 
computed from (4.9.1).
 
\head 5. $H$-biequivariant Hopf algebroids\endhead

In this chapter we generalize the notions of an $\h$-algebra, 
$\h$-bialgebroid, $\h$-Hopf algebroid to the case when 
the Lie algebra $\h$ is not necessarily commutative, and define 
quantum counterparts of the quasiclassical 
notions introduced in Chapters 1-2 of \cite{EV}. 

We will define the notions of an $H$-biequivariant 
Hopf algebroid and quantum groupoid.
The notion of an $H$-biequivariant quantum groupoid is 
a quantum analogue of the notion of an $H$-biequivariant 
Poisson groupoid, introduced in \cite{EV}. We will also introduce
less general notions of a dynamical quantum groupoid
and Hopf algebroid, which 
are quantum analogues of the notions of a dynamical Poisson groupoid
and Hopf algebroid. 

In this chapter we will work 
mostly in the setting of perturbation theory. That is, quantum objects 
will be defined over $k[[\hbar]]$, where $k$ is some field, 
and give classical objects modulo $\hbar$ and quasiclassical 
ones modulo $\hbar^2$. We discuss the relationship between the quasiclassical 
and quantum objects,  and questions regarding quantization.  

\subhead 5.1. Quantization of Poisson algebras\endsubhead

In this section we will remind some well known facts from the theory of 
deformation quantization. 

Let $k$ be a field of characteristic zero. Let $K=k[[\hbar]]$.
By a topologically free $K$-module we mean a 
$K$-module of the form $V[[\hbar]]$, where $V$ is a $k$-vector space. 
All $K$-modules we will use will be topologically free. 
By tensor product of two such modules we will always mean 
completed tensor product over $K$. 
  
Let $A_0$ be a commutative algebra over $k$ with $1$.
Recall that according to Grothendieck, a linear operator 
$D: A_0\to A_0$ is a differential operator of order 
$\le N, N\ge 1$ if for any $a\in A_0$ the operator $f\to D(af)-aDf$
is a differential operator of order $\le N-1$, and a differential 
operator of order $0$ is the operator of multiplication by an element of
$A_0$. If $A_0$ is the algebra of regular functions on a manifold 
(smooth, analytic, algebraic, formal) then ``differential operator
of order $N$'' means what it usually means in geometry.  
    
Let $A_0$ be a Poisson algebra over $k$ with $1$, with Poisson bracket 
$\{,\}$. 
Recall that by a quantization of $A_0$ is meant 
a $K$-module $A=A_0[[\hbar]]$ equipped with a K-linear binary 
operation $*: A\o A\to A$, which defines 
an associative algebra structure on $A$, such that $A/\hbar A=A_0$ 
as an algebra, 
and $\frac{1}{\hbar}
(f*g-g*f)\text{ mod }\hbar=\{f,g\}$, $f,g\in A_0\subset A$. 
In this case $A_0$ is called the quasiclassical limit of $A$. 

Let $f,g\in A_0$. Then 
$$
f*g=fg+\hbar c_1(f,g)+\hbar^2c_2(f,g)+...,\tag 5.1.1 
$$
where $c_i: A_0\o A_0\to A_0$ are linear maps. 
A quantization defined by (5.1.1) is called local if 
$c_i(f,g)$ is a differential operator in $f$ and $g$ for any $i$. 
If $A_0$ is the algebra $\O(X)$ of regular functions on a smooth manifold $X$, 
and $A$ is a local quantization of $A_0$, then $A$ defines 
(by formula (5.1.1)) a quantization
$A_U$ of the algebra $(A_U)_0=\O(U)$ of 
regular functions on any open subset $U$ of 
$X$. In other words, it defines a quantization of the sheaf of regular 
functions.  
This holds also in the holomorphic and algebraic situations, if $X$ is 
affine.  

Let $X$ be a manifold, and let $T^*X$ be its cotangent bundle.
Let $A_0=\O(T^*X)_p$ be the Poisson algebra of regular functions 
on $T^*X$ which are fiberwise polynomial of a uniformly bounded degree. 
This Poisson algebra has a 
distinguished quantization 
$A=\O_q(T^*X)_p$ called the canonical quantization
(q is not a parameter here but the first letter of the word 
``quantum''). Namely,  
$A$ is the algebra of formal series of the form $\sum_{n\ge 0} \hbar^nD_n$, where 
$D_i$ are differential operators on $X$, 
such that $n\ge \text{order}(D_n)$, and $n-\text{order}(D_n)\to +\infty$, as
as $n\to \infty$. It is easy to check that this quantization is local, 
so it defines a quantization $A_U=\O_q(U)$ of the Poisson algebra
$(A_U)_0=\O(U)$ of regular functions on an open subset $U\in T^*X$. 

Let $\g$ be a Lie algebra, and $\g^*$ be its dual space, with the
usual Poisson structure. Consider the Poisson algebra 
$\O(\g^*)_p$ of polynomial functions on $\g^*$. This algebra 
has a distinguished quantization $A=\O_q(\g^*)_p$, called the
geometric quantization. 
Namely,  
$A$ is the algebra of formal series of the form $\sum_{n\ge 0} \hbar^nD_n$, 
where 
$D_i\in U(\g)$, $n\ge \text{order}(D_n)$,
and $n-\text{order}(D_n)\to +\infty$,
$n\to \infty$. It is easy to check that this quantization is local, 
so it defines a quantization $A_U=\O_q(U)$ of the Poisson algebra
$(A_U)_0=\O(U)$ of regular functions on an open subset $U\in \g^*$. 

\subhead 5.2. $H$-biequivariant associative algebras
\endsubhead

In this section we will introduce the notion of an $H$-biequivariant 
associative algebra. This notion is a quantum analogue of 
the notion of an $H$-biequivariant Poisson algebra, 
introduced in a previous paper \cite{EV}.

Let $A$ be an associative algebra
over $K$ with 1, which is commutative mod $\hbar$, 
$H$ a connected affine algebraic group
over $k$, and
$\psi: A\times H\to A$ be a right algebraic
action of $H$ on $A$ by automorphisms, defined over $k$.
This means that $A$, as a representation of $H$, has the form 
$A_0[[\hbar]]$, where $A_0$ is a sum of finite dimensional 
representations of $H$ over $k$. 

Let $\h$ be the Lie algebra of $H$.
Let $U\subset \h^*$ be an $H$-invariant open set. A 
homomorphism $\mu: \O_q(U)\to A$ is called a quantum moment map for $\psi$
if for any linear function
on $U$ given by $a\in \h$ and any $f\in A$ we have  
$$
[\mu(a),f]=\hbar d\psi|_{h=1}(a,f).\tag 5.2.1
$$ 
Here
$d\psi|_{h=1}:\h\times A\to A$ is the differential of $\psi$ at $h=1\in H$. 
Using the Leibnitz identity for the operator $g\to [\mu(g),f]$,
from (5.2.1) one can compute $[\mu(g),f]$ for any rational function $g$. 

For a left action of $H$ a quantum moment map is defined in the same way, 
with the only difference that it is an anti-homomorphism rather than a 
homomorphism. 

\proclaim{Definition} An $H$-biequivariant associative algebra over $U$
is a 5-tuple $(A,l,r,\mu_l,\mu_r)$, where 
$A$ is an associative
algebra with $1$ over $K$, which is commutative mod $\hbar$, 
$l,r$ is
a pair of commuting algebraic 
 actions of $H$ on $A$ (a left action and a right action) 
by algebra automorphisms, defined over $k$, 
and $\mu_l,\mu_r: \O_q(U)\to A$ 
are quantum moment maps for $l$, $r$, such that 

(i) $\mu_l,\mu_r$ are embeddings, and their images commute;

(ii) There exists an $l(H)\times r(H)$-invariant 
k-subspace $A_0^l$ of $A$ such that the multiplication map $\mu_r(\O_q(U))
\o A_0^l
\to A$ is a linear isomorphism; 
there exists an $l(H)\times r(H)$-invariant 
k-subspace $A_0^r$ of $A$ such that the multiplication map $ 
\mu_l(\O_q(U))\o A_0^r\to A$ is a linear isomorphism. 

A morphism of $H$-biequivariant associative algebras over $U$ is a morphism
of algebras which preserves $l,r$ and $\mu_l,\mu_r$. 
\endproclaim

{\bf Remark 1.} From $[l,r]=0$ it follows that $[\mu_l\circ x,
\mu_r\circ y]$ is 
a central element for $x,y\in\h$, but it does not follow 
that this commutator equals $0$. So we require that it is zero by condition 
(i). 

{\bf Remark 2.} Condition (ii) is of technical nature and is not very 
important in the discussion below.

Denote the category of $H$-biequivariant associative algebras over $U$ by 
$\Cal A_U^q$ (q stands for ``quantum'').

For convenience we will write $l(h)a$ as $ha$ and $r(h)a$ as $ah$. 

Let us now describe the monoidal structure on $\Cal A_U^q$.

Let $A,B\in \Cal A_U^q$. Then the group $H$ acts in
$A\o B$ by $\Delta(h)(a\o b)=ah^{-1}\o hb$. We will construct a new 
$H$-biequivariant associative algebra $A\wo B$, which is obtained by 
quantum Hamiltonian reduction of $A\o B$ by the action of $H$. 
 
Denote by $A*B$ the space $A\o_{\O_q(U)}B$, where $\O_q(U)$ is mapped to
$A$ via $\mu_r^A$, and to $B$ via $\mu_l^B$, acting in both algebras 
from the left. Then $A*B$ is the 
quotient of $A\o B$ by the linear span $I$ of elements of the form
$\mu^A_r(f)a\o b-a\o \mu^B_l(f)b$, $f\in \O_q(U)$, $a\in A,b\in B$.  
The space $A*B$ has two commuting actions of $H$ ($l_A\o 1$ and $1\o r_B$). 
But
we cannot claim that $A*B\in \Cal A_U^q$, since 
the algebra structure on $A\o B$ does not, in general, descend to 
$A*B$ ($I$ is only a right ideal and not necessarily 
a left ideal). 

However, the action $\Delta$ of $H$ on $A\o B$ descends to
one on $A*B$, so we can define $A\wo B:=(A*B)^H$, where $H$ acts by $\Delta$. 

\proclaim{Proposition 5.1} The algebra structure on $A\o B$ descends to
one on $A\wo B$.
\endproclaim

\demo{Proof} Let $x,y\in A\wo B$. We can regard $x,y$ as elements
of $A*B$. Choose their liftings $X=\sum a_i\o b_i$, $Y=\sum c_i\o d_i$
into $A\o B$. By definition, $xy$ is the image of $XY$ 
in $A*B$. 

We have to check two things. 

1. That $xy$ is $H$-invariant.

2. That $xy$ does not depend on the choice of liftings $X,Y$.

First we check property 1. Since $x,y$ are $H$-invariant, we have
$$
\gather
\sum [\mu^A_r(z),a_i]\o b_i+\sum a_i\o [\mu^B_l(z),b_i]\in I, \\
\sum [\mu^A_r(z),c_i]\o d_i+\sum c_i\o [\mu^B_l(z),d_i]\in I, z\in \h. 
\tag 5.2.2\endgather
$$
Therefore, since $I$ is a right ideal, 
$$
\sum [\mu^A_r(z),a_ic_j]\o b_id_j+
\sum a_ic_j\o [\mu^B_l(z),b_id_j]\in XI+I.\tag 5.2.3
$$

{\bf Lemma 5.1} If $X$ is $H$-invariant modulo $I$, then $XI\subset I$. 

{\it Proof of the Lemma.} 

Since $\sum c_j\o d_j$ is $H$-invariant modulo $I$, for any
$z\in \h$ we have
$$
\sum c_j\mu^A_r(z)\o d_j-\sum c_j\o d_j\mu^B_l(z)\in I.\tag 5.2.4
$$
Therefore, the same equality holds any rational function $g\in \O_q(U)$
instead of $z$.  
 This proves the Lemma.

The Lemma shows that the RHS 
of (5.2.3) is in $I$, i.e. $xy$ is $H$-invariant. 

Now we check property 2. 
If $X',Y'$ are any other liftings of $x$ and $y$, then 
$X-X'\in I$, and $Y-Y'\in I$. 
So it remains to show that $X(Y-Y')\in I$. 
But this follows from the Lemma. 
$\square$
\enddemo

Thus, we have shown that the product descends to $A\wo B$. 
 The two commuting actions of $H$ on $A\o B$ 
by $(h_1,h_2)(a\o b)=h_1a\o bh_2$, and the corresponding
quantum moment maps descend to $A\wo B$. 
So, in order to check that $A\wo B\in \Cal A_U^q$, it 
suffices to check properties (i) and (ii). 

Using 
properties (i) and (ii) 
of the quantum moment maps $\mu_l^A,\mu_r^A,\mu_l^B,\mu_r^B$,
it is easy to see that $A*B$ is naturally identified with
$\mu_l^A(\O_q(U))\o A_0^r\o B_0^r$, and $A\wo B$ is identified with
$\mu_l^A(\O_q(U))\o (A_0^r\o B_0^r)^H$, where $H$ acts by
$a\o b\to ah^{-1}\o hb$. This implies properties (i) and (ii) 
for the quantum moment map $\mu_l^A\o 1: \O_q(U)\to A\wo B$, corresponding to
the left action of $H$ on $A\wo B$ (with $(A\wo B)_0^r=(A_0^r\o B_0^r)^H$). 
For the quantum moment map $1\o \mu^B_r:\O_q(U)\to A\wo B$ corresponding 
to the right action, these properties are proved 
analogously. 

Thus, $A\wo B\in\Cal A_U^q$. 
It is clear that the assignment $A,B\to A\wo B$ is a bifunctor
$\Cal A_U^q\times \Cal A_U^q\to \Cal A_U^q$. 

Recall \cite{EV} that $(T^*H)_U$ denotes 
the variety of points $(h,p)\in T^*H$ such that $h^{-1}p\in U$. 
Consider the algebra $\O_q((T^*H)_U)$, which is the canonical 
quantization of the standard
symplectic structure on $(T^*H)_U$. 
It is equipped with the standard actions $l,r$ 
of $H$ on left and right given by
$(x,p)\to (h_1xh_2,h_1ph_2)$ (these actions obviously respect the 
quantization). 

Let  
$\mu_{l,r}:\O_q(U)\to \O_q((T^*H)_U)$ be the embeddings, which 
assign to an element of $U(\h)$ the corresponding right-,
respectively left-invariant differential operator on $H$.  
 It is easy to check 
that $\mu_{l,r}$ are quantum moment maps for $l,r$.  

Let $\bold 1=(\O_q((T^*H)_U),l,r,\mu_l,\mu_r)$. It is easy to check that 
we have natural isomorphisms $A\wo \bold 1\equiv A\equiv \bold 1\wo A$.

\proclaim{Proposition 5.2} (i) $(A\wo B)\wo C=A\wo (B\wo C)$. 

(ii) $\bold 1$ is a unit object in $\Cal A_U^q$ with respect to 
$\wo$, and 
$(\Cal A_U^q,\wo,\bold 1)$ is
a monoidal category.
\endproclaim

\demo{Proof} Easy.
$\square$\enddemo

Let $A\in \Cal A_U^q$. Denote by $\bar A$ the new object of $\Cal A_U^q$
obtained as follows: $\bar A$ is $A^{op}$  (the opposite algebra), 
with the left and the right actions of $H$ permuted
(i.e. the left, respectively right, action of $h$ on $\bar A$ is the right, 
respectively left, action of $h^{-1}$ on $A$), and the quantum moment maps
also permuted. We will call $\bar A$ the dual object to $A$. 
By a {\it quasireflection} on $A$ we will mean a morphism 
$i:\bar A\to A$. Note that unlike \cite{EV}, here we do not
require that $i^2=1$. 

Let $A\in \Cal A_U^q$ and $i:\bar A\to A$ be a quasireflection.
Let $\phi^i_+,\phi^i_-: A\o A\to A$ be given by the formulas $\phi^i_+(a\o b)= 
ai(b)$, $\phi^i_-(a\o b)=i(a)b$. It is easy to see that these 
maps descend to linear 
maps $\psi^i_\pm: A\wo A\to A$. 

\subhead 5.3. $H$-biequivariant Hopf algebroids\endsubhead

Now let us define the quantum version of the notion of an
$H$-biequivariant Poisson-Hopf algebroid. 

\proclaim{Definition} Let $A$ be an $H$-biequivariant 
associative algebra.
Then $A$ is called an $H$-biequivariant Hopf algebroid over $U$
if it is equipped with a coassociative $\Cal A_U^q$-morphism   
$\Delta: A\to A\wo A$ called the coproduct, 
a $\Cal A_U^q$-morphism $\epsilon: A\to \bold 1$ called the counit, and 
a quasireflection $S: \bar A\to A$ called the antipode,   
such that 

(i) $(id\bullet \epsilon)\circ \Delta=(\epsilon \bullet id)\circ \Delta
=id$, and 

ii) $\psi_+^S\circ \Delta=\mu_l\circ P\circ \epsilon$,
$\psi_-^S\circ \Delta=\mu_r\circ P\circ \epsilon$, where $P: \bold 1\to
\O_q(U)$ is the map which assigns to a differential operator on $H$
its value at the identity element (which is in $U(\h)$).   
\endproclaim 

The same structure without the antipode will be called an $H$-biequivariant 
bialgebroid. 

If $H=1$, then these notions coincide with notions of a Hopf algebra 
and a bialgebra over $K$. 

{\bf Remark 1.} In the above discussion, $U$ is a Zariski open set. 
If $k=\Bbb R$ or $\Bbb C$, then we can take $U$ to be an open set 
in the usual sense, and define $\O(U)$ to be the algebra of smooth,
respectively analytic, functions on $U$. Then we can repeat sections 
5.2, 5.3, and thus define the notions of
an $H$-biequivariant
associative algebra and Hopf algebroid
over $U$. 
Similarly, one can take $U$ to be the infinitesimal neighborhood of zero 
in $\h^*$ (i.e. $\O(U)=k[[\h]]$). The material of Sections 5.2 and 5.3
can be generalized to this case as well. 

{\bf Remark 2.} In the smooth, 
analytic, and formal case one has to drop the condition 
that $A$ is the sum of finite dimensional representations 
of $H$ (because $\O_q(U)$ does not satisfy this condition). 
One should instead require that $A$ is a representation of $\h$. 
One should also impose the locality condition
for a quantum moment map $\mu$: for any $f\in A$
the operation $g\to [\mu(g),f]$ is local in $g$, in the sense
that $[\mu(g),f]=\sum \mu(D_ig)f_i$, where $f_i\in A$, and $D_i$
are h-adically convergent series of differential operators on $U$.  
Using (5.2.1) and the locality property, one can compute
$[\mu(g),f]$ not only for rational functions $g$ but for arbitrary 
smooth, holomorphic, or formal functions. 

\subhead 5.4. Quantization of $H$-biequivariant Poisson-Hopf algebroids
and Poisson groupoids \endsubhead

Consider the following two settings. 

1. Let $A_0$ be an $H$-biequivariant Poisson algebra
(see Section 2.3 of \cite{EV}). 
Let $A=A_0[[\hbar]]$. Suppose that $A$ is equipped with 
an associative product $*$ in such a way that 
$A$ is a local quantization of $A_0$ as a Poisson algebra, and the 5-tuple
$(A,l,r,\mu_l,\mu_r)$ is an $H$-biequivariant associative algebra 
(where $l,r,\mu_r,\mu_r$ are the $K$-linear extensions
of the structure maps of $A_0$ to $A$). 

2. Assume that in addition 
$A_0$ is an $H$-biequivariant Poisson-Hopf algebroid, 
i.e. it is equipped with maps $\Delta_0,\epsilon_0,S_0$ satisfying
certain axioms (see Section 2.4 of \cite{EV}). Suppose that $A$ 
is as above, and in addition that $A$ is equipped with maps
$\Delta,\epsilon,S$, which make $A$ an $H$-biequivariant 
Hopf algebroid, and equal $\Delta_0,\epsilon_0,S_0$ modulo $\hbar$. 

\proclaim{Definition} In these cases,
$A_0$ is called the quasiclassical limit of $A$, and 
$A$ is called a quantization of $A_0$. 
\endproclaim

If $H=1$, then this definition is the usual definition of 
a quantization of a Poisson and Poisson-Hopf algebra. 

Now consider the geometric version of this definition. 
Let $X$ be an $H$-biequivariant Poisson manifold over $U$. 
Let $A_0=\O(X)$. Then $A_0$ satisfies the axioms of an 
$H$-biequivariant Poisson algebra, except for maybe property (ii). 
The notion of quantization of $A_0$ is defined as above. 
A quantization $A$ of $A_0$ will be called an 
{\it $H$-biequivariant quantum space}. 

If $X$ is in addition an $H$-biequivariant Poisson groupoid, 
then $A_0$ satisfies the axioms of an $H$-biequivariant 
Poisson-Hopf algebroid, except for property (ii) 
and the fact that the coproduct 
$\Delta$ maps $A_0$ to $A_0^2:=\O(X\bullet X)[[\hbar]]$, which is a completion 
of $A_0\wo A_0$, but not to $A_0\wo A_0$ itself. 
(This problem already exists for Lie groups, where 
the coproduct maps $\O(G)$ to $\O(G\times G)$ and not to $\O(G)\times \O(G)$). 
The notion of quantization of $A_0$ is defined as above. 
The quantization is called local if $f*g$ is a bidifferential 
operator of $f,g$ modulo any power of $\hbar$, and $\Delta(f)=D\Delta_0(f)$, 
where $D$ is a differential operator modulo any power of $\hbar$.   
A local quantization $A$ of $A_0$ will be called an $H$-biequivariant quantum 
groupoid. 

Suppose that $X=X(G,H,U)$ is a dynamical Poisson groupoid
(see Chapter 1 of \cite{EV}), and $A_0=\O(X)$ is as 
above. In this case a local quantization $A$ of $A_0$ will be called a 
{\it dynamical quantum groupoid}. 
If the subspace $\O(U)\o \O(G)\o \O(U)[[\hbar]]\subset A$ is 
closed under the product, then it is an  $H$-biequivariant 
Hopf algebroid. Such 
Hopf algebroid is called a {\it dynamical Hopf algebroid}.

Recall that by a preferred quantization of a Poisson Lie group 
is meant a quantization in which the coproduct is undeformed. 
The notion of a preferred quantization of an $H$-biequivariant Poisson 
groupoid or Poisson-Hopf algebroid is defined in the same way. 

\proclaim{Conjecture} (i) Any dynamical Poisson groupoid admits a 
quantization. 

(ii) Any quasitriangular dynamical Poisson groupoid
admits a preferred quantization.
\endproclaim

In the case $H=1$ (Poisson-Lie groups), this conjecture goes back to Drinfeld
and is proved in \cite{EK1,EK2}.
 
\subhead 5.5. The case $H=(\C^*)^N$\endsubhead

In this section we will 
consider the special case when $H=(\C^*)^N$, and 
establish the connection between the constructions of this chapter and 
Chapter 4.  

Let $H=(\C^*)^N$. In this case, the 
main notions of Chapter 5 are simplified:

1. Since $H$ is commutative,
the algebra $\O_q(U)$ is just $\O(U)[[\hbar]]$. 

2. Denote by $P\subset \h^*$ be the lattice of 
characters of $H$ ($P=\Z^n$).
Let $A$ be an $H$-biequivariant associative algebra.
Then the algebra $A$ can be written as $A=\oplus_{\al,\beta\in P}A_{\al\beta}$,
where $A_{\al\beta}$ is the set of elements $a\in A$ such that $h_1ah_2=
\alpha(h_1)\beta(h_2)a$ (the direct sum is understood in the 
$\hbar$-adically complete sense). The images of the maps $\mu_l,\mu_r$ are in 
$A_{00}$. The product $A\wo B$ can be written in the form 
$(A\wo B)_{\alpha\delta}=\oplus_{\beta\in P} 
A_{\alpha\beta}\o_{\O(U)}B_{\beta\delta}$ , where $\O(U)$ is embedded in $A$ 
via $\mu_r^A$ and in $B$ via $\mu_l^B$, and acts form the left
(thus this product is similar to the matrix product). 

3. The algebra $\O_q((T^*H)_U)=\bold 1$ can be written in form 
$\O(U)\times \O(H)[[\hbar]]=\O(U)\o \C[P][[\hbar]]$, where 
the commutation relations between $P$ and $\O(U)$ are given by
$f\chi=\chi f^\chi$, $f\in \O(U)$, $\chi\in P$, where 
$f^\chi(u)=f(u+\hbar\chi)$. 

In particular, in this case we can replace the algebra $\O(U)$ with 
the field $M_{\h^*}$ of meromorphic functions on
$\h^*$, imposing the locality condition 
(see Remark 2, Section 5.3). Then equation (5.2.1)
together with the locality condition implies identities (4.1.1). 

Now
nothing prevents us from setting $\hbar$ to be 
no longer a formal parameter, but a nonzero complex number $\gamma$. 
In this situation, it is easy to see that an 
$H$-biequivariant algebra (bialgebroid, Hopf algebroid) is 
the same as an $\h$-algebra ($\h$-bialgebroid, $\h$-Hopf algebroid)
with weights belonging to $P\subset \h^*$. 
This gives a connection between Chapters 4 and 5. 

\head 6. $\h$ bialgebroids associated to 
quantum dynamical R-matrices of Hecke type.\endhead

\subhead 6.1. The Hecke condition\endsubhead

Let $R:\h^*\to \End(V\o V)$ be a quantum dynamical R-matrix
with step $\gm$.
Consider the $\h$-bialgebroid $\bar A_R$ introduced in Chapter 4. 

It is clear that if $R=1$ and $\gm=0$ 
then $\bar A_R=M_{\h^*}\o M_{\h^*}\o \O(\End (V))$.  
Therefore, for $R\ne 1$
 we want the algebra $\bar A_R$ to look like a quantum deformation
of $M_{\h^*}\o M_{\h^*}\o \O(\End (V))$. 

A natural formalization of this wish is the PBW property, defined below. 

The algebra $\bar A_R$ has a natural $\Z_+$-grading, given by 
$deg(f(\la^i))=0,deg(L_{ab})=1$. Denote by $\bar A_R^n$ 
the degree $n$ component of $\bar A_R$. It is clear that 
$\bar A_R^n$ are $M_{\h^*}\o M_{\h^*}$-modules, 
where the two components of $M_{\h^*}$ act by 
left multiplication by $f(\la^1)$ and $f(\la^2)$.  

\proclaim{Definition} The algebra $\bar A_R$ is said to satisfy 
the Poincare-Birkhoff-Witt (PBW) property if 
the  $M_{\h^*}\o M_{\h^*}$-module $\bar A_R^n$
is isomorphic to the free module $M_{\h^*}\o M_{\h^*}\o S^n\text{End}(V)$.
\endproclaim

For a general dynamical R-matrix, the PBW property is not the case. 
However, the property holds if one imposes an additional 
``Hecke type'' condition on $R$. 

\proclaim{Definition} $R$ is said to be of strong Hecke type 
if 

(i) $R$ satisfies equation (1.3.6) for some nonzero
parameters $p,q\in \C$, 
$p\ne -q$, such that $q/p$ is not a root of unity,
and 

(ii) There exists a continuous family $R(t)$, $t\in [0,1]$, of
quantum dynamical R-matrices with step $\gm(t)$, satisfying (i)
with parameters $p(t)$, $q(t)$, such that 
$R(0)=1,p(0)=q(0)=1,\gm(0)=0$, $R(1)=R,p(1)=p,q(1)=q,\gm(1)=\gm$. 
\endproclaim

{\bf Example.} It is easy to see from the classification 
that all dynamical R-matrices 
of $gl_N$ Hecke type are of strong Hecke type. Thus, for 
dynamical R-matrices of $gl_N$-type, strong Hecke type is the same as 
the Hecke type. 

\proclaim{Theorem 6.1} If $R$ is of strong Hecke type 
then $\bar A_R$ satisfies the PBW property. 
\endproclaim

This theorem explains the meaning of the Hecke type conditions 
introduced in Chapter 1. If $\h=0$, this theorem is well known 
(see \cite{FRT}). 
 
\subhead 6.2. Proof of Theorem 6.1\endsubhead

 Let $\tilde A$ be the algebra 
with the same generators as $\bar A_R$ and the same relations 
except the Yang-Baxter relation. Then, as a vector space, the algebra 
$\tilde A$ has the form $\oplus_{n\ge 0}\tilde A^n$, $\tilde A^n=M_{\h^*}\o 
M_{\h^*}\o (\text{End}(V))^{\o n}$, and $\bar A_R$ is the quotient of
$\tilde A$ by the Yang-Baxter relation. 

Let $H_n(v)$ be the Hecke algebra of type $A_n$ with parameter $v$. 
It is the algebra generated by elements $T_i$, $1\le i\le n-1$, with relations
$$
[T_i,T_j]=0, |i-j|\ge 2;\ T_iT_{i+1}T_i=T_{i+1}T_iT_{i+1};\ 
(T_i-1)(T_i+v)=0\tag 6.2.1
$$
If $v$ 
is not a root of unity of degree $n$, 
this algebra is isomorphic to $\C[S_n]$ and therefore semisimple.

Denote by $R^{ii+1}(\la)$ the operator 
$1^{i-1}\bo R(\la)\bo 1^{n-i-1}:V^{\o n}\to M_{\h^*}\o V^{\o n}$,
where $\bo$ has the meaning defined by (3.1.2). 

If $R$ satisfies condition (i), then we have an action 
of $H_n(v)$, $v=q/p$, on the $M_{\h^*}\o M_{\h^*}$-module  $\tilde A^n$, 
defined by the formula
$$
T_iX= P_{ii+1}:R^{ii+1}(\la^1)XR^{ii+1}(\la^2)^{-1}:
P_{ii+1},\tag 6.2.2
$$  
where $P_{ii+1}$ is the permutation of the $i$-th 
and the $i+1$-st components in the tensor product $V^{\o n}$.
This construction explains the origin of the term ``Hecke type''. 

The Yang-Baxter relation in $A_R$ implies that
the degree $n$ component $\bar A^n_R$ of $\bar A_R$ is isomorphic to 
the space of coinvariants of $T_1,...,T_{n-1}$ 
in $\tilde A^n$. By semisimplicity
of $H_n(v)$, this space is isomorphic 
to the space of vectors in $M_{\h^*}\o M_{\h^*}\o (\End (V))^{\o n}$,
 which are invariant under $T_i$. 

Now recall that $R$ satisfies condition (ii). Let 
$R(t)$ be the corresponding family. Consider the 
corresponding modules $\bar A_{R(t)}^n$. 
Since they can be defined both as coinvariants and invariants, 
their dimensions cannot jump, which implies that $\bar A_{R(0)}^n$
is isomorphic to $\bar A_{R(1)}^n$ as a $M_{\h^*}\o M_{\h^*}$-module.
However, by our assumptions, $\bar A_{R(0)}^n=  
M_{\h^*}\o M_{\h^*}\o S^n\text{End}(V)$, while 
$\bar A_{R(1)}^n=\bar A_R^n$. This proves the Theorem. 
$\square$

\subhead 6.3. Hecke condition and quantization\endsubhead

Theorem 6.1 has the following generalization to the case when 
the step $\gm$ is a formal parameter. 

Let $R_\gm=1-\gm r+\sum \gm^nr_n$ be a formal series whose 
coefficients are meromorphic functions $\h^*\to \End(V\o V)$. 
Suppose that $R$ is a quantum dynamical R-matrix with step $\gm$.
Let $\bar A_{R_\gm}$,
$A_{R_\gm}$ denote the algebras over $K:=\C[[\gm]]$ defined as 
in Chapter 4. 

It is clear that $\bar A_{R_\gm}/\gm\bar A_{R_\gm}=M_{\h^*}\o M_{\h^*}\o 
\O(\text{End}(V))$. Thus the analogue of the PBW property 
for $\bar A_{R_\gm}$ in 
this case is the property that the $K$-module 
$\bar A_{R_\gm}$ is a topologically free module, i.e. provides
a flat deformation of $M_{\h^*}\o M_{\h^*}\o 
\O(\text{End}(V))$.

\proclaim{Theorem 6.2} If $R_\gm$ satisfies the Hecke equation 
(1.3.6) for some $p(\gm)=1+O(\gm)$, $q(\gm)=1+O(\gm)$, then 
$\bar A_{R_\gm}$ is 
a flat deformation of $M_{\h^*}\o M_{\h^*}\o 
\O(\text{End}(V))$.
\endproclaim

\demo{Proof} Analogous to the proof of Theorem 6.1
$\square$\enddemo 

\proclaim{Corollary 6.1} Under the assumption of Theorem 6.2, 
$A_{R_\gm}$ is 
a flat deformation of $M_{\h^*}\o M_{\h^*}\o 
\O(GL(V))$.
\endproclaim

If $R_\gm$ is holomorphic 
in an open set $U\subset \h^*$ then 
 we can define algebras $\bar A_{R_\gm}^U$,
$A_{R_\gm}^U$ in the same way as 
$\bar A_{R_\gm}$, $A_{R_\gm}$, except that 
$M_{\h^*}$ 
is replaced with the algebra of holomorphic functions
$\O(U)$ on $U$. It is clear 
that Theorem 6.2 and Corollary 6.1 are valid for these algebras:

\proclaim{Proposition 6.1}
Under the assumptions of Theorem 6.2, the
algebras $\bar A_{R_\gm}^U$,
$A_{R_\gm}^U$ are topologically free over $K$. 
\endproclaim

Now let $R_\gm:U\to \End(V\o V)[[\gm]]$ 
be a quantum dynamical R-matrix 
holomorphic on $U$ which satisfies the condition of Theorem 6.2.
Let $p(\gm)=1+a\gm+O(\gm^2)$, $q(\gm)=1+b\gm+O(\gm^2)$, $\gm \to 0$. 
Then 
from the quadratic equation for $R^\vee$ we get
the unitarity condition
$$
r^{21}+r=(b-a)P-(b+a),\tag 6.3.1
$$
and from the quantum 
dynamical Yang-Baxter equation for $R$ we get the classical
dynamical Yang-Baxter equation for $r$. Thus, 
according to Chapter 1 of \cite{EV}, $r$ defines a structure 
of a quasitriangular dynamical Poisson 
groupoid on $U\times GL(V)\times U$.
In particular, we have the 
corresponding dynamical Poisson-Hopf algebroid 
$A_r^{0U}=\O(U)\o \O(GL(V))\o \O(U)$
(here $\O(G)$ denotes the algebra of polynomial functions on $G$). 

\proclaim{Theorem 6.3} The dynamical Hopf algebroid
$A_{R_\gm}^U$ is a quantization of the dynamical 
Poisson-Hopf algebroid $A_r^{0U}$. 
\endproclaim

\demo{Proof} Since we know that $A_{R_\gm}^U$ is topologically free, 
the proof is the direct computation
of the quasiclassical limit and then comparison with
Chapter 1 of \cite{EV}.
$\square$\enddemo

Let $G=GL(V)$, $H$ be a maximal torus in
$G$, and $U\subset \h^*$ a polydisc. 
Let $X(G,H,U)$ be the Lie groupoid $U\times G\times U$
with two actions of $H$, defined in Chapter 1 of \cite{EV}. 

\proclaim{Theorem 6.4} 
Any structure of a quasitriangular dynamical Poisson groupoid
on $X(G,H,U)$ admits a preferred quantization. 
\endproclaim

\demo{Proof} The statement follows from Theorem 1.6 and Theorem 6.3. 
$\square$\enddemo

{\bf Remark.} Notice that if $R_\gm$ fails to satisfy the Hecke condition 
modulo $\gm^2$, then the algebra $A_{R_\gm}$ is not topologically free. 
Indeed, in this case $r$ does not satisfy 
the unitarity condition, so according to Chapter 1 of \cite{EV}
the bracket defined by $r$ on $U\times GL(V)\times U$ is not Poisson 
(i.e. does not satisfy the Jacobi identity). This means that 
the corresponding deformation is not flat, since a flat deformation
of a commutative algebra induces a Poisson structure on this algebra. 
Thus, the Hecke condition seems to be intrinsic for good properties of 
the algebra $A_R$.

\vskip10ex 

P.Etingof
Department of Mathematics, Harvard University, Cambridge, MA 02138

{\it E-mail address: etingof\@math.harvard.edu.}

A.Varchenko,
Department of Mathematics, Phillips Hall, University of North Carolina at
Chapel Hill, Chapel Hill, NC 27599-3250, USA

{\it E-mail address: av\@math.unc.edu.}

\end